\documentclass[12pt]{article}
\usepackage{WideMargins2}
\def\Aut{{{\rm Aut}}}

\def\order{{{\rm order}}}
\def\GL{{{\rm GL}}}
\def\PGL{{{\rm PGL}}}

\def\PSL{{{\rm PSL}}}

\def\SL{{{\rm SL}}}
\def\GCD{{{\rm GCD}}}

\def\mult{{{\rm mult}}}
\def\N{{{\rm N}}}

\def\F{{{\mathbb F}}}
\def\Z{{{\mathbb Z}}}

\def\cj#1{{{\overline{#1}}}}

\def\deg{{{\rm deg}}}

\def\textmatrix#1#2#3#4{{ \left({#1 \atop #3}{#2\atop #4}\right)}}
\def\set#1{{\left\{\,#1\,\right\}}}
\let\a=\alpha
\let\b=\beta
\let\g=\gamma
\let\d=\delta

\let\l=\lambda
\let\ep=\varepsilon

\let\x=\times
\let\s=\sigma
\let\z=\zeta
\let\ee = e
\def\ie{{ \it i.e.}}
\def\eg{{ \it e.g.}}

\def\jacobi#1#2{{{\left(\frac{#1}{#2}\right)}}}
\def\choose#1#2{{{\left({#1}\atop {#2}\right)}}}

\def\calC{{{\cal C}}}
\def\calK{{{\cal K}}}
\def\calL{{{\cal L}}}
\def\calO{{{\cal O}}}
\def\calP{{{\cal P}}}

\def\inv{{{\rm inv}}}

\title{\Large\bfseries  Explicit Artin maps into $\PGL_2$}
\author{Antonia W.\ Bluher \\ National Security Agency}
\date{June 2019, revised August 2021}
\begin {document}
\fancytitle
\begin{abstract}
Let $G\subset \PGL_2(\F_q)$ where $q$ is any prime power, and let $Q\in \F_q(x)$ such that $\F_q(x)/\F_q(Q)$ is a Galois extension with group $G$. 
By explicitly computing the Artin map on unramified degree-1 primes in $\F_q(Q)$ for various groups $G$,
interesting new results emerge about finite fields, additive polynomials,
and conjugacy classes of $\PGL_2(\F_q)$.
For example, by taking $G$ to be 
a unipotent group, one obtains
a new characterization for when an additive polynomial splits completely over $\F_q$.
When $G = \PGL_2(\F_q)$, one obtains information about conjugacy classes 
of $\PGL_2(\F_q)$.
When $G$ is the group of order~3 generated 
by $\textmatrix {1\,}{-1}{1\,}{\,\,0}$, one obtains a natural tripartite symbol on $\F_q$
with values in $\Z/3\Z$. 
Some of these results generalize to $\PGL_2(K)$ for arbitrary fields $K$.
\end{abstract}

\section{Introduction} \label{sec:intro}

Let $K$ be a field and $G$ a finite subgroup of $\PGL_2(K)$. It is
well known, and will be proved in Section~\ref{sec:Q}, that there is 
$Q \in K(x)$
such that $K(x)/K(Q(x))$ has Galois group~$G$. Normalize $Q$ so that
$Q(\infty) = \infty$; $Q$ will be called a {\it quotient map} for~$G$. If $\tau \in \cj K\cup\{\infty\}$, then $Q^{-1}(\tau)$ is a $G$-orbit in~$\cj K\cup\{\infty\}$.
If $|Q^{-1}(\tau)|=|G|$, then $\tau$ is said to be {\it regular}. 
Let $\s \in \Aut(\cj K/K)$ such that $\s(\tau)=\tau$.  Then $Q^{-1}(\tau)$ is closed under $\sigma$, and it is a $G$-orbit,
so for each $v \in Q^{-1}(\tau)$ there is $\gamma \in G$ such that $\sigma(v) = \gamma(v)$.
If $\tau$ is regular, then $\gamma$ is uniquely determined by $v$ and~$\sigma$, and the conjugacy class
$\calC_{\gamma,G} = \set{\a \g \a^{-1} : \a \in G}$ is uniquely determined by $\tau$ and $\sigma$.
In that case, define 
$$\inv_Q(\tau,\sigma) = \calC_{\gamma,G}.$$
If $G$ is abelian, so that $\calC_{\gamma,G} = \{\gamma\}$, then we write more simply $\inv_Q(\tau,\sigma)=\gamma$. We abbreviate $\calC_\g=\calC_{\g,G}$ if $G$ is clear from the context.

If $K=\F_q$ and $\sigma(v)=v^q$, where $q$ is any prime power, then we write $\inv_Q(\tau,q)$ instead of $\inv_Q(\tau,\s)$,
or just $\inv(\tau)$ if $Q$ and $q$ are clear from the context.
In that case,
Xander Faber astutely observed that the map $\tau \mapsto \inv(\tau) = \calC_\g$ is essentially the Artin map for the extension 
$\F_q(x)/\F_q\left(Q(x)\right)$.
The connection is as follows. 
Let $\tau \in \F_q$.
The polynomial $Q-\tau\in \F_q(Q)$ corresponds to a degree-1 place $P=(Q-\tau)$ of $\F_q(Q(x))$. This place is unramified in $\F_q(x)$
if and only if $\tau$ has $|G|$ distinct preimages in $\cj \F_q$; 
this coincides with our definition that $\tau$ is regular. 
If $v \in Q^{-1}(\tau)$
and $g$ is its minimal polynomial over $\F_q$, then $g$ corresponds to a place $\calP$ in $\F_q(x)$ that lies over $P$. 
The element $\gamma$ such that $v^q = \gamma(v)$ is the {\it Frobenius automorphism} of $\calP$ for the extension $\F_q(x)/\F_q(Q)$, 
and the map $P \mapsto \calC_\g$
is the {\it Artin map}. For a more general discussion of the Artin map over function fields, see Rosen \cite[Ch.~9]{Rosen}. 
The Artin map is defined in \cite[page~122]{Rosen}. Because of this close relation with the Artin map, we call $\inv_Q(\tau,\s)$ 
the {\it Artin invariant} of $\tau$ with respect to $Q$ and $\s$.

While the existence and general properties of the Artin map are widely known, what has not been previously appreciated is 
how interesting the examples are, even in the genus-0 case, {\it i.e.}, over the rational function field $\F_q(x)$.
Thus, the emphasis in this article is not so much on the existence of $\inv(\tau)$, but rather on the
wealth of arithmetic information that is revealed by specific examples.  Some of this arithmetic information was previously known, and some is new.

\begin{example} \label{example:QRes} 
The simplest example is $G=\left\{\textmatrix 1001,\textmatrix{-1}001\right\}$ and $Q(x)=x^2$ over $\F_q$.
Here one should assume $1\ne -1$, \ie, $q$ is odd.
All $\tau \in \F_q$ are regular except $\tau=0$.  If $v^2=\tau\ne 0$ then 
$$v^q = (v^2)^{(q-1)/2}v = \tau^{(q-1)/2} v = \jacobi \tau q v,$$
where $\jacobi \tau q$ is the quadratic residue symbol: 1 if $\tau$ is a nonzero square, or $-1$ if
$\tau$ is a nonsquare.  Thus, $v^q = \gamma(v)$ with $\gamma = \textmatrix {\jacobi \tau q}  0 0 1$, and $\inv(\tau) = \gamma$.
As can be seen, the Artin invariant for this group is intimately connected with the quadratic residue symbol.
\qed
\end{example}

The next example is related to the theory of Kummer extensions.  For an exposition on Kummer extensions, 
see \cite[\S II.L]{Artin} or \cite[Ch.~4, \S4]{Neukirch}.  If $(n,q)=1$ then let $\mu_n$ denote the group of $n$th roots of unity
in $\cj\F_q$.

\begin{example} \label{example:Kummer} 
Let $q$ be any prime power (possibly even).
Suppose that $n$ divides $q-1$, and let $G = \set{\textmatrix \zeta001 : \zeta^n=1} \subset \PGL_2(\F_q)$.
The quotient map is $Q(x)=x^n$, and every element of $\F_q$ is regular except~0.
Let $\tau \in \F_q^\x$ and $\inv(\tau) = \textmatrix \z 001$.
Then $\tau = v^n \in \F_q^\x \iff v^q = \zeta v$.  
The field extension $\F_q(v)=\F_q(\tau^{1/n})$ is a Kummer extension. The Galois group sends $v$ to $v^q=\inv(\tau) (v) = \z v$. 
To express $\inv(\tau)$ directly in terms of $\tau$, note that $\zeta = v^{q-1} = (v^n)^{(q-1)/n}=\tau^{(q-1)/n}$.
Example~\ref{example:QRes} is the special case $n=2$.
\qed
\end{example}

\begin{example}{\bf (Klein Group)} \label{example:Klein} Let $G = \set{\textmatrix 1001, \textmatrix {-1}001, \textmatrix 0110,\textmatrix 0{-1}10}\subset
\PGL_2(\F_q)$. Here assume
$1 \ne -1$, so $q$ is odd. 
$Q(x) = (x+1/x)^2/4$ is a quotient map.  Let $\tau=Q(v) \in \F_q$.
Then $\tau$ is regular iff $\tau \not \in \{0,1\}$ iff $v^4 \ne 1$.  The following theorem pertains to this example.

\noindent{\bf Theorem.}\ {\it Let $q$ be an
odd prime power. Every element $\tau \in \F_q$ can be written as 
$$\tau = (v+1/v)^2/4,\qquad\text{$v\in\mu_{2(q-1)}\cup\mu_{2(q+1)}$.}$$
Moreover, if $\tau \not\in\{0,1\}$, then $v^{q-AB}=A$, where
$A = \jacobi\tau q$ and $B=\jacobi{\tau-1}q$. }
\medskip

The theorem implies $v^q = A v^{AB} = \gamma(v)$, where 
\begin{equation} \g=\inv(\tau,q) = \textmatrix A001 \textmatrix 0110^{(1-AB)/2}.
\label{eq:Klein}
\end{equation}
In other words, $\inv(\tau,q)$ is given explicitly in terms of $\jacobi \tau q$ and $\jacobi{\tau-1}q$. 
The above theorem is proved in \cite[Theorem~4.1]{Wilson-like} and is used in \cite{Permutation,Wilson-like} 
to obtain a new factorization formula for Dickson and Chebyshev polynomials 
and new theorems in elementary number theory. For a short and self-contained exposition on these topics, see \cite{Structure}.
In fact, the article \cite{Structure} directly motivated the current article.
\qed
\end{example}

The next example (explained in Section~\ref{sec:unipotent}) computes the Artin invariant when $G$ is a unipotent subgroup of 
$\PGL_2(\F_q)$, and leads to unexpected results about additive polynomials that have all their roots in the ground field.  

\begin{example} \label{example:1b01}  Let $G = \set{\textmatrix 1b01  : b \in \F_q}$.
A quotient map is $Q_G(x)=x^q-x$. If $Q_G(v)=\tau \in \F_q$, then $v^q=v+\tau = \textmatrix 1\tau 01(v)$. Thus, $\inv(\tau)=\textmatrix 1\tau 01$.

More generally,
let $q=p^n$, where $p$ is a prime or a prime power, let $W\subset \F_q$ be a $d$-dimensional $\F_p$-vector subspace
of $\F_q$, and $G_W=\set{\textmatrix 1w01 : w \in W}$. As is well known, the quotient map is $Q_W(x)= \prod_{w\in W}(x-w)$, which is an $\F_p$-additive polynomial,
	\ie, $Q_W(\l x+y)=\lambda Q_W(x)+Q_W(y)$ for $\lambda \in \F_p$. $Q_W$ can be regarded as an $\F_p$-linear map from $\F_q$ to $\F_q$ with kernel $W$.
Let $Y=Q_W(\F_q)$. Then $Y$
is an $(n-d)$-dimensional $\F_p$-vector subspace of $\F_q$, and for $\tau \in \F_q$ we will prove in Section~\ref{sec:unipotent} 
that 
	$$\inv_{Q_W}(\tau) = \textmatrix {1\ } {Q_Y(\tau)} {0\ } 1 \in G_W.$$
In particular, $Q_Y(\tau) \in W$, which is not otherwise obvious.  This observation implies the
following symmetric relation between $W$ and~$Y$.  Though easy to prove, to our knowledge it was not previously noticed.

\medskip 

\noindent{\bf Proposition.}\ {\it Let $q=p^n$, let $W$ be an $\F_p$-vector subspace of $\F_q$, $Q_W(x) = \prod_{w \in W} (x-w)$, and $Y=Q_W(\F_q)$. Then
$W=Q_Y(\F_q)$, and there are short exact sequences
$$0 \to W \stackrel{inc.}{\longrightarrow} \F_q \stackrel{Q_W}{\longrightarrow} Y \to 0  $$
and 
$$0 \to Y \stackrel{inc.}{\longrightarrow} \F_q \stackrel{Q_Y}{\longrightarrow} W \to 0.  $$ }
\qed
\end{example}

These observations lead to a simple characterization of when an $\F_p$-additive polynomial splits in $\F_q$.
(See Proposition~\ref{prop:linPoly}).

\newpage

\noindent{\bf Splitting Criterion.}\ 
{\it
Let $q=p^n$ where $p$ is a prime or a prime power, and let $L(x) = x^{p^d} + \sum_{i=0}^{d-1} a_i\, x^{p^i}$ be an $\F_p$-additive polynomial, where $a_i \in \F_q$, $a_0 \ne 0$, and $d\ge1$. Then all the roots of $L$
are in $\F_q$ if and only if there is an $\F_p$-additive polynomial $M(x) = x^{p^{n-d}} + \sum_{i=0}^{n-d-1} b_i x^{p^i} \in \F_q[x]$ with $M \circ L(x) = x^q - x$.
In that case, it is also true that $L \circ M(x) = x^q - x$.}

\medskip

This splitting criterion is simpler than
characterizations that are currently in the literature.
The criterion in current use is as follows.  (See McGuire and Sheekey \cite{MS} and Csajb\'{o}k {\it et al} \cite{CMPZ19}.)
Let $L(x) = x^{p^d} + \sum_{i=0}^{d-1} a_i x^{p^i}\in \F_q[x]$, and define $d \x d $ matrices $C_L$ and $A_L$ by
$$C_L = \begin{pmatrix} 0 & 0 \cdots & 0 & -a_0 \\
1 & 0 \cdots & 0 & -a_1 \\
0 & 1 \cdots & 0 & -a_2 \\
\vdots & \vdots & \cdots & \vdots \\
0 & 0 \cdots & 1 & -a_{d-1}
\end{pmatrix} $$
$$A_L = C_L C_L^{(p)}  \cdots C_L^{(p^{n-1})}$$
where $C^{(p^i)}$ means raising every matrix entry to the power $p^i$.  Then $L$ has all its roots in $\F_q$ if and only if $A_L$ is equal to the identity matrix.

As an example, let $q=P^7$ and $$L(x) = x^{P^3} - b x^P - a x,$$ where $a,b \in \F_{q}$ and $a\ne0$. It was shown by Csajb\'{o}k {\it et al} \cite[Theorem 3.3]{CMPZ18} using combinatorial arguments that $L$ can have all its roots in $\F_{q}$ only if $q$ is even. A complete characterization of when $L$ has all its roots in $\F_{q}$ was found by G. McGuire and D. Mueller \cite{MM}. One can obtain this result more simply using the new splitting criterion. Namely, let $M(x) = x^{P^4}+u_3 x^{P^3}+u_2 x^{P^2}+u_1 x^P+u_0 x\in\F_{P^7}[x]$ and try to solve $M\circ L(x) = x^{P^7}-x$.  From the coefficients of $x^{P^6}$, $x^{P^5}$, $x^{P^4}$, and $x$ one finds that $u_3=0$, $u_2=b^{P^4}$, $u_1 = a^{P^4}$, and $u_0=1/a$. The equation $M\circ L(x)=x^{P^7}-x$ then simplifies to
$$0 = (1/a-b^{P^4+P^2})x^{P^3} - (a^{P^2}b^{P^4}+a^{P^4}b^P)x^{P^2}
-(a^{P^4+P} + b/a)x^P.$$
Then $a=b^{-P^4-P^2}$, and in particular $b\ne0$. From the coefficient of $x^{P^2}$, and using $b^{P^7}=b$, we have 
$0=a^{P^2} b^{P^4} + a^{P^4} b^P = b^{-P^6-P^4} b^{P^4} + b^{-P^8-P^6}b^P = 2b^{-P^6}$. Thus, $2=0$, showing $q$ is even. Finally,
$a^{P^4+P}+b/a=0$ yields $b=a^{P^4+P+1}=(b^{-P^4-P^2})^{P^4+P+1}$, which simplifies to ${\rm N}_{\F_{q}/\F_P}(b)=1$. The conclusion is that $L$ has all its roots in $\F_{q}$ iff $q$ is even,
${\rm N}_{\F_{q}/\F_P}(b)=1$, and $a=b^{-P^4-P^2}$.

\begin{example} \label{example:pgl2}  
	The case $G=\PGL_2(\F_q)$, studied in Section~\ref{sec:PGL2},  reveals nontrivial information about conjugacy classes of $\PGL_2(\F_q)$.   
If $K$ is any field and if $g=\textmatrix abcd \in \GL_2(K)$, let 
\begin{equation}
\iota(g) = \frac{(a+d)^2}{ad-bc}. \label{eq:iota2}
\end{equation}
Then $\iota(cg) = \iota(g)$ for $c\in K^\x$, so $\iota$ is well defined on $\PGL_2(K)$.  Also, $\iota(hgh^{-1})=\iota(g)$ for $h \in \GL_2(K)$,
	so $\iota$ is constant on conjugacy classes.  If $e_1,e_2$ are the roots of the characteristic polynomial of $g$, then $\iota(g)=e_1/e_2+e_2/e_1+2$.

	The following results concerning the map $\iota$ will be proved in Section~\ref{sec:PGL2}.

	\noindent{\bf Theorem.}\ {\it A quotient map for $\PGL_2(\F_q)$ is
	$Q(x) = (x^{q^2}-x)^{q+1}/(x^q-x)^{q^2+1}$, and $\tau \in \F_q$ is regular with respect to $Q$ iff $\tau \ne 0$.
	Let $\calK$ denote the set of conjugacy classes $\calC_\g$ such that $\circ(\g)\ge 3$. 
	Then $\iota$ induces a bijection from $\calK$ onto $\F_q^\x$, and the inverse bijection is $\inv_Q$.  
	\qed}

The bijections in the theorem are pictured here:
$$ \fbox{$\calC_\g$ s.t. $\circ(\g)\ge 3$}\qquad      {{\stackrel\iota\longrightarrow} \atop \stackrel{\inv_Q}\longleftarrow}\qquad   \fbox{$\F_q^\x$ } $$

	Let $H$ be any subgroup of $G=\PGL_2(\F_q)$, $Q_H$ a quotient map for $H$, and $Q_G(x)=(x^{q^2}-x)^{q+1}/(x^q-x)^{q^2+1}$.
	It is shown in Lemma~\ref{prop:H} that there is a unique rational function $h\in \F_q(x)$ such that $Q_G=h\circ Q_H$. We will prove:

	\bigskip\noindent{\bf Theorem.}\ {\it Let $H,Q_H,h$ be as above. Suppose $\tau \in \F_q$ is regular with respect to $Q_H$ and let $\inv_{Q_H}(\tau,q)=\calC_{\g,H}$. If $\g=1$ then
	$h(\tau)=\infty$. If $\g\ne1$ then $h(\tau)=\iota(\g)$. \qed}
\end{example}

\begin{example} \label{example:order3} Section~\ref{sec:three} considers 
	$G=\{I,\beta,\beta^2\}\subset \PGL_2(K)$, where $\beta = \textmatrix 1{-1}10$ and $K$ is any field. 
	A quotient map is $Q(x)=(x^3-3x+1)/(x(x-1))$, and $\tau\in\cj K$ is regular if and only if $\tau^2-3\tau+9\ne0$.  
	Let $\s \in \Aut(\cj K/K)$, $\tau\in \cj K$ 
	such that $\tau^2-3\tau + 9\ne 0$ and $\s(\tau)=\tau$. Let $v\in\cj K$ such that $Q(v)=\tau$. 
	Then there is a unique $\ell\in \Z/3\Z$ such that $\s(v)=\b^\ell(v)$. By definition, $\beta^\ell=\inv_Q(\tau,\s)$.  
	Section~\ref{sec:three} presents formulae for $\ell$ in terms of $\tau$, as follows. 

	\noindent{\bf Theorem.}\ {\it With notation as above, if $\inv_Q(\tau,\s)=\b^\ell$ then $\ell \pmod 3$ is determined by:
\begin{enumerate} 
	\item If char$(K)\ne 3$, let $\omega\in \cj K$ denote a primitive cube root of unity and let $\zeta\in\cj K$ satisfy
		$\zeta^3=(\tau+3\omega^2)/(\tau+3\omega)$.  Then $\s^2(\zeta)/\zeta = \omega^{\ell}$.
	\item If char$(K)=3$, let $\zeta\in\cj K$ satisfy $\zeta^3-\zeta=1/\tau$. Then $\ell=\s(\zeta)-\zeta$.
	\item In the special case where $K=\F_q$ and $\sigma$ is the Frobenius, $\s(x)=x^q$, then 
		$$\begin{cases} 
			\omega^\ell= \left(\frac{\tau+3\omega^2}{\tau+3\omega}\right)^{(q^2-1)/3} & \text{if $3\nmid q$,}\\
			\ell=\Tr_{\F_q/\F_3}(1/\tau) & \text{if $3|q$.}
		\end{cases}
		$$
\end{enumerate}
	}
\qed
\end{example}

This article has three parts. Part~{\ref{part:general} (Sections~\ref{sec:orbits}--\ref{sec:invariant}) presents the general theory.
Specifically, Section~\ref{sec:orbits} discusses $G$-orbits in $\cj K \cup\{\infty\}$,
Section~\ref{sec:Q} discusses existence and computation of quotient maps, 
and Section~\ref{sec:invariant} defines the invariant $\inv(\tau)$.
Part~\ref{part:K} (Sections~\ref{sec:known}--\ref{sec:six}) considers finite groups that are naturally defined over $\PGL_2(K)$ for any $K$. 
Section~\ref{sec:known} considers $G=\{\textmatrix 1001,\textmatrix 0110\}\subset \PGL_2(K)$ and shows how the Artin invariant for this group
in the case $K=\F_q$
is related to the well-known fact: $\F_q = \set{\z+1/\z : \z^{q-1}=1\ {\rm or}\ \z^{q+1}=1}$.
Section~\ref{sec:K} generalizes Examples~\ref{example:Kummer} and~\ref{example:Klein} to arbitrary fields~$K$.
Section~\ref{sec:three} considers the order-3 group given in Example~\ref{example:order3}. 
Section~\ref{sec:six} considers the dihedral group of order~six in $\PGL_2(K)$ generated by $\textmatrix1{-1}10$ and $\textmatrix 0110$.
Part~\ref{part:q} considers subgroups of $\PGL_2(\F_q)$, including Borel subgroups, unipotent subgroups, cyclic subgroups,
$\PGL_2(\F_q)$, and $\PSL_2(\F_q)$.  Some applications are given.

\part{General Theory} \label{part:general}

\section{Orbits} \label{sec:orbits}

Let $K$ be any field and $\cj K$ its algebraic closure. 
The projective linear group $\PGL_2(K)$ is defined as the group of invertible $2\x2$ matrices
with entries in $K$, modulo the scalar matrices, $\textmatrix c00c$, where $c \in K^\x$.
As is well known, if $K \subset L$ where $L$ is a field, then $\PGL_2(K)$ acts on $L\cup \{\infty\}$ via
$$\textmatrix abcd  (v)  = \frac{av+b}{cv+d}.$$
This equation is self-explanatory if $v \in L$ and $cv+d\ne 0$. If $v=\infty$, then 
$\textmatrix abcd (v) = a/c$, where we interpret $a/0 = \infty$.  Also, if $v\in L$ and $cv+d=0$, then $\textmatrix abcd (v) = \infty$.
The reader can verify that if $\gamma,\delta \in \PGL_2(K)$ then $\gamma\left(\delta (v)\right) = (\gamma \delta) v$.
$\PGL_2(K)$ acts triply transitively on $K\cup\{\infty\}$, \ie, for any distinct $a,b,c \in K \cup \{\infty\}$ there is $\gamma \in \PGL_2(K)$
taking $\infty$ to $a$, 0 to $b$, and 1 to $c$.  In fact, $\gamma$ is uniquely determined and it equals
$$\gamma = \begin{pmatrix} a(b-c) & b(c-a) \\ b-c & c-a \end{pmatrix}$$
if $a,b,c$ are all finite, or $\textmatrix {c-b} b01$ if $a=\infty$, $\textmatrix a{c-a}10$ if $b=\infty$, $\textmatrix a{-b}1{-1}$ if $c=\infty$.
Thus, $\PGL_2(K)$ is in one-to-one correspondence with the set of ordered triples $(a,b,c)$ of
distinct elements in $K \cup \{\infty\}$, and in particular
\begin{equation} |\PGL_2(\F_q)| = (q+1)q(q-1). \label{eq:cardPGL2} \end{equation}

Let $G$ be a finite subgroup of $\PGL_2(K)$, and let $|G|$ denote its cardinality. 
If $v\in L\cup \{\infty\}$, then the $G$-orbit containing $v$, or simply {\it orbit} if $G$ is clear from context, is
defined as
$$\calO_v = \set{\gamma(v) : \gamma \in G}.$$
Orbits partition $L \cup \{\infty\}$ into disjoint sets.
We will say an orbit is {\it short} if it has fewer than $|G|$ elements;
otherwise the orbit is {\it full-sized}. 

\begin{lemma} \label{lem:short} Let $G$ be a finite subgroup of $\PGL_2(K)$ and let $L/K$ be an extension of fields.
An element $v\in L\cup\{\infty\}$ belongs to a short orbit if and only if there is $\gamma \in G$, $\gamma \ne \textmatrix 1001$, 
such that $\gamma (v) = v$. Every short orbit is contained in $K\cup\{\infty\}$ or a quadratic extension of $K$.
The union of short orbits has at most $2(|G|-1)$ elements.
\end{lemma}

\begin{proof}  Let $S$ denote the union of short orbits in $L$.  Then
$$ v \in S \iff \text{$\calO_v$ is short} \iff \gamma_1(v) = \gamma_2(v)\ \text{for a pair of distinct elements $\gamma_1,\gamma_2 \in G$},$$
and in that case, $\gamma_1^{-1}\gamma_2$ fixes $v$. Thus,
$$S = \cup_{1\ne \gamma \in G} F_\gamma,$$
where $F_\gamma = \set{v \in L\cup\{\infty\} : \gamma(v)=v}$.   Since $|F_\g|\le 2$ and consists of rational elements in $K\cup\{\infty\}$ or a pair of conjugate elements, the result follows.
\end{proof}

\begin{lemma} \label{lem:orbitMultiplicity}  Let $G$ be a finite subgroup of $\PGL_2(K)$, and let $\calO \subset L \cup \{\infty\}$ be a $G$-orbit,
where $L/K$ is an extension field.
Then $|\calO|$ divides $|G|$, and each element of $\calO$ is fixed by exactly $|G|/|\calO|$ elements of $G$.  The integer
$\mult(\calO)=|G|/|\calO|$ is called the {\it multiplicity of $\calO$}. 
If $\calO \ne \calO_\infty$ and $v \in \calO$, then
$$\prod_{\gamma \in G} (x-\gamma(v)) = \left(\prod_{w \in \calO} (x-w)\right)^{\mult(\calO)}.$$
\end{lemma}
\begin{proof} This follows from standard facts about groups acting on sets, as can be found for example in \cite[Section 4.1, Prop.~2]{DF}. \end{proof}

If $v \in \cj K$,
define $\deg_K(v) = [K(v):K]$. Note that $\deg_K(v) = \deg_K(\gamma(v))$ for all
$\gamma \in \PGL_2(K)$, because $v$ and $\gamma(v)$ generate the same field over $K$. Consequently, 
$\deg_K(v)$ is constant on orbits. If $K=\F_q$, we write $\deg_q(v)$ instead of $\deg_K(v)$. Then $\F_q(v)=\F_{q^t}$, where $t = \deg_q(v)$.

\begin{lemma} \label{lem:degree} Let $\gamma \in \PGL_2(\F_q)$ and suppose $\circ(\gamma)=t>1$, where $\circ(\gamma)$ denotes the order of $\g$.
If $v \in \cj\F_q$ and
$v^q = \gamma(v)$ then $v^{q^i}=\g^i(v)$ for all $i\ge1$, and $\deg_q(v)$ divides $t$. If in addition $v,\g(v),\ldots,\g^{t-1}(v)$ are distinct, 
then $\deg_q(v) = t$. 
\end{lemma}

\begin{proof} Note that $v^{q^2} = (v^q)^q = (\gamma(v))^q$. Since $\gamma$ has entries in $\F_q$, this equals $\gamma(v^q) = \gamma(\gamma(v))=\gamma^2(v)$.
By induction one can show that $v^{q^i} = \gamma^i(v)$ for all $i\ge1$. Thus, the $G$-orbit $\calO_v$ is the set of $\F_q$-conjugates of $v$, where
$G$ is the cyclic group of order $t$ generated by~$\g$. Then $\deg_q(v)=|\calO_v|=t/\mult(\calO_v)$. This shows that $\deg_q(v)$
divides $t$, and $\deg_q(v)=t$ iff $\calO_v$ has full size, \ie, iff $\g^i(v)$ for $0\le i < t$ are distinct.
\end{proof}

The next proposition will be useful in determining how many field elements $\tau$ have the same invariant $\calC_\g$, assuming that $\circ(\gamma)\ge 3$.
See Section~\ref{subsec:vqgamma} and Proposition~\ref{prop:factoring} for 
further study of the equation $v^q=\gamma(v)$.

\begin{proposition} \label{prop:counting} Suppose that $\gamma \in \PGL_2(\F_q)$ has order $t \ge 3$. Let
$$A_{\g,q} = \set{v \in \cj\F_q : v^q = \gamma(v)},\qquad \calL_{\g,q} = A_{\g,q} \setminus \F_{q^2}.$$
Let $G \subset \PGL_2(\F_q)$ be a group that contains $\gamma$ and let $\calC_\g = \set{\a\g\a^{-1}:\a\in G}$. 
\begin{enumerate} 
\item[({\it i})] If $v \in A_{\g,q}$ then the $G$-orbit $\calO_v = \set{\b(v):\b \in G}$ has full size if and only if $v \in \calL_{\g,q}$.  
\item[({\it ii})] $|\calL_{\g,q}|=q+\kappa$ with $\kappa \in \{0,1,-1\}$, and $t$ divides $q+\kappa$. (Note that $\kappa$ is uniquely determined
from $t$, since $t\ge 3$ and $\kappa \equiv - q \pmod t$).
\item[({\it iii})] Let $\calL = \cup_{\b \in \calC_\g} \calL_{\b,q}$. 
Then $\calL$ decomposes into exactly $r$ $G$-orbits, all of full size, where $r = |\calC_\g|(q+\kappa)/|G|$. 
\end{enumerate}
\end{proposition}

\begin{proof} {\it (i)} Let $v \in A_{\g,q}$, so $v^q=\gamma(v)$. We will show that the $G$-orbit $\calO_v$ has full size if and only if $v\in \calL_{\g,q}$.
Since all short orbits are contained in $\F_{q^2} \cup \{\infty\}$, $v \in \calL_{\g,q}$ implies that $\calO_v$ has full size.
Conversely, if $\calO_v$ has full size then $\gamma^i(v)$ are distinct for $i=0,1,\ldots,t-1$, so $\deg_q(v)=t$
by Lemma~\ref{lem:degree}.  In particular, $v \not \in \F_{q^2}$ so $v \in \calL_{\g,q}$.

{\it (ii)}  Write $\gamma = \textmatrix abcd$. $A_{\gamma,q}$ is the set
of solutions in $\cj\F_q$ to $f(x) = 0$, where $f(x)=x^q(cx+d)-ax-b$. We claim the roots are distinct. 
For if $r$ is a repeated root, then $f(r)=f'(r)=0$, so $r^qc-a=0$.
Either $c=a=0$ (contradicting that $ad-bc \ne 0$) or $r=a/c$.  But $r=a/c$ implies 
$f(r)=(a/c)(a+d)-a^2/c-b=(ad-bc)/c\ne0$, contradicting that $r$ is a root of $f$.  This establishes that $f$ has no repeated roots, so it has $\deg(f)$
distinct roots in $\cj\F_q$.  So $|A_{\g,q}|=\deg(f)$, which equals $q+1$ if $c\ne 0$, or $q$ if $c=0$.

Let $X=A_{\g,q}\cap \F_{q^2}$.  Any $v \in X$ satisfies $v = v^{q^2} = \gamma^2(v)$, so it is a fixed
point of $\gamma^2$. There are at most two fixed points, so $|X|\le 2$.  Further, if $c=0$ then $\gamma^2$ fixes $\infty$,
so it can fix at most one other point, and it follows that $|X|\le 1$ when $c=0$. 

Since $A_{\g,q}$ is the disjoint union of $\calL_{\g,q}$ and $X$, $|\calL_{\g,q}|=|A_{\g,q}| - |X|$.   If $c=0$ then $|A_{\g,q}|=q$ and $|X|  \in \{0,1\}$,
and if $c\ne 0$ then $|A_{\g,q}| = q+1$ and $|X| \in \{0,1,2\}$. In either case, $|\calL_{\g,q}|\in \{q-1,q,q+1\}$, \ie, 
$|\calL_{\g,q}|=q+\kappa$ where $\kappa \in \{-1,0,1\}$.

To see that $t$ divides $|\calL_{\g,q}|$, observe that $\g$ permutes $\calL_{\g,q}$ and has no fixed points.  The permutation breaks into cycles, each of order $t$,
so the cardinality of $\calL_{\g,q}$ must be a multiple of~$t$.  

{\it (iii)}  If $v \in \calL$ and $\a \in G$ then we claim $\a(v) \in \calL$. Indeed, $v^q = \b (v)$ for some 
$\b=\ep\g\ep^{-1} \in \calC_\g$,
so if $w=\a(v)$ then $w^q=\a(v^q)=\a\ep\g\ep^{-1}(v)= \a\ep\g\ep^{-1}\a^{-1}(w)=(\a\ep)\g(\a\ep)^{-1}(w)$. Further, $\deg_q(w)=\deg_q(v)=t$, so $w\not\in \F_{q^2}$.
This shows $\a(v)\in \calL$, as claimed. Then $\calL$ splits into $G$-orbits. All have full length by {\it (i)}, so $|G|$ divides $|\calL|$, and the number of $G$-orbits is $|\calL|/|G|$.
Finally, $\calL$ is a disjoint union of the sets $\calL_{\b,q}$ with $\b \in \calC_\g$, because if $v^q=\b(v)$ and $v^q=\b'(v)$ then $\b(v)=\b'(v)$, $\b^{-1}\b'(v)=v$.
Since $\deg_q(v)>2$, this forces $\b=\b'$. Note that $\order(\ep\g\ep^{-1})=\circ(\g) = t$ for each $\b = \ep\g\ep^{-1} \in \calC_\g$, so each set\
$\calL_{\b,q}$ has the same cardinality, $q + \kappa$, where $\kappa\in \{-1,0,1\}$ and $\kappa \equiv -q \pmod t$.   We conclude that $|\calL|=(q+\kappa)|\calC_\g|$
and the number of $G$-orbits is $r=(q+\kappa)|\calC_\g|/|G|$.
\end{proof}

\section{Quotient maps} \label{sec:Q}

Let $G$ be a finite subgroup of $\PGL_2(K)$, where $K$ is any field.
A {\it quotient map} for $G$ is a rational function $Q(x)$ such that the extension field $K(x)/K(Q)$ has Galois group $G$. We further require that $Q(\infty)=\infty$.
This section gives proof of existence, properties, examples, and computational aspects of quotient maps.
We remark that existence of $Q$ is well known to algebraic geometers; see
Mumford \cite{Mumford}.

\subsection{Existence of quotient maps.}

The existence of a quotient map essentially follows from Galois theory of the field $K(x)$ (see Artin \cite{Artin}), together with some facts about subfields of $K(x)$
(see van der Waerden \cite{VdW}), where $K(x)/K$ is transcendental.

Every nonzero $f \in K(x)$ can be written uniquely as $p_1(x)/p_2(x)$, where
$p_1$ and $p_2$ are relatively prime polynomials and $p_2$ is monic.  Define $$\deg(f) = \max\{\deg(p_1),\deg(p_2)\}.$$
If $\gamma = \textmatrix abcd \in \PGL_2(K)$, then it may be viewed as an element of $K(x)$ of degree~1: $$\gamma(x) = (ax+b)/(cx+d),\qquad \deg(\g)=1.$$  

\begin{lemma} \label{lem:degf} {\it (i)}\ If $f\in K(x)$ is nonconstant then $[K(x):K(f)] = \deg(f)$. \\ \noindent{\it (ii)}\  If $f,g \in K(x)$
are nonconstant then $\deg(f\circ g)= \deg(f)\,\deg(g)$. 
\end{lemma}

\begin{proof}  {\it (i)}\ See \cite[Section 10.2]{VdW}. \\ {\it (ii)}\ Let $y=g(x)$. Since $[K(x):K(y)] = \deg(g)$, we have $[K(y):K] = \infty$ and so $y$ is transcendental.
By {\it (i)},  $\deg(f \circ g) = [K(x):K(f(g(x)))] = [K(x):K(y)] [K(y):K(f(y))] = \deg(g)\,\deg(f)$. \end{proof}

\begin{corollary} \label{cor:autKx}
If $\gamma \in \PGL_2(K)$ and $f \in K(x)$, define $A_\gamma(f) = f\circ \gamma^{-1}$. Then $\gamma \mapsto A_\gamma$ is an isomorphism
from $\PGL_2(K)$ onto $\Aut(K(x)/K)$.
\end{corollary}

\begin{proof} This is well known, but we give a proof for completeness. $\Aut(K(x)/K)$ is defined as the group of isomorphisms from $K(x)$ to $K(x)$ that
fix all elements of $K$.
First, $A_\gamma \in \Aut(K(x)/K)$ because $A_\g(f+g)=A_\g(f)+A_\g(g)$, $A_\g(f) A_\g(g) = A_\g(fg)$ when $f,g \in K(x)$, 
$A_\g(c) = c$ when $c \in K$, and $A_\g^{-1}=A_{\g^{-1}}$.  Clearly $A_\g=1 \iff \g^{-1}(x)=x \iff \g=1$.
Further, $A_{\gamma} (A_\delta(f)) = A_\gamma (f \circ \delta^{-1}) = f\circ \delta^{-1} \circ
\gamma^{-1} = f \circ (\gamma\delta)^{-1} = A_{\gamma\delta}(f)$.  So $\PGL_2(K)$ injects into $\Aut(K(x)/K)$, and we just need to show it is surjective.
Let $A$ be any automorphism of $K(x)/K$.
Since $x$ generates $K(x)$ over $K$, so does $A(x)$. Then $\deg(A(x))=1$ by Lemma~\ref{lem:degf}.
Write $A(x)= (ax+b)/(cx+d)$, where $ax+b$ and $cx+d$ have no common factor and
$a$ or $c$ is nonzero.  Then $ad-bc\ne 0$, so $\gamma=\textmatrix abcd$ is in $\PGL(2,K)$. Evidently $A(x)=A_{\gamma^{-1}}(x)$, and since an automorphism
of $K(x)/K$ is determined by the image of $x$, it follows that $A = A_{\gamma^{-1}}$. This proves surjectivity.
\end{proof}

Let $\Sigma$ be the fixed field of $G$:
\begin{equation} \Sigma = \set{f(x) \in K(x) : \text{$f\circ \gamma(x) = f(x)$ for all $\gamma \in G$} }. \label{eq:Sigma}
\end{equation}

\begin{proposition} \label{prop:Q0}
There is a function $Q(x) \in K(x)$ of degree $|G|$ such that $\Sigma = K(Q)$.  Moreover, $[K(x):\Sigma]=|G|$,  $K(x)/\Sigma$ is Galois,
and its Galois group is isomorphic to $G$. If $Q'(x) \in \Sigma$ and $\deg(Q')=|G|$ then there is $\alpha \in \PGL_2(K)$ such that $Q' = \alpha \circ Q$.
\end{proposition}

\begin{proof} Let $x,y$ be independent transcendentals and consider 
\begin{equation}  F(y) = \prod_{\gamma \in G} y - \gamma(x) \in K(x)[y]. \label{eq:Fy} \end{equation}
$F(y)$ has degree $|G|$ and its coefficients are in $\Sigma$. 
Since $F \in \Sigma[y]$ and $F(x)=0$, this shows $K(x)$ is an algebraic extension of $\Sigma$ and $[K(x):\Sigma] \le |G|$.
The group $G$ is contained in $\Aut(K(x)/K)$ by Corollary~\ref{cor:autKx}, and it fixes all elements of $\Sigma$, therefore $[K(x):\Sigma] \ge |G|$
by Galois theory (see the corollary to Theorem~13 in \cite{Artin}). Combining these inequalities gives $[K(x):\Sigma]=|G|$. Since the degree of the extension equals 
the order of the group of automorphisms of $K(x)$ that fix $\Sigma$, the extension is Galois.

L\"uroth's Theorem \cite[\S10.2, {\it p.}~218]{VdW} states that any field $E$ such that $K\subset E \subset K(x)$ and $[K(x):E]<\infty$ has the form $E=K(f)$,
where $f\in K(x)\setminus K$. Therefore, $\Sigma = K(Q)$ for some $Q \in K(x)$. By Lemma~\ref{lem:degf}{\it (i)},
$\deg(Q)=[K(x):K(Q)]$, which equals $|G|$. 

If $Q' \in \Sigma=K(Q)$ then $Q'=h(Q)$ for some $h \in K(x)$. By Lemma~\ref{lem:degf}{\it (ii)}, 
if $\deg(Q')=|G|$ then $\deg(h)=1$, so $h \in \PGL_2(K)$.
\end{proof}

\begin{proposition} \label{prop:QExists}
Let $K$ be any field and let $G$ be a finite subgroup of $\PGL_2(K)$.
There is a rational function $Q\in K(x)$ such that 
\begin{enumerate}
\item $Q(\gamma x) = Q(x)$ for all $\gamma \in G$; 
\item If $Q$ is written as a reduced fraction, \ie, $Q=f/g$ where $f,g \in K[x]$
and $\GCD(f,g)=1$, then $|G|=\deg(f) > \deg(g)$.
\end{enumerate}
Further, if $\widetilde Q$ is another function with these properties, then
$\widetilde Q(x) = a Q(x) + b$ for some $a\in K^\x$ and $b\in K$.
\end{proposition}
\begin{proof} Let $Q_0=f_0/g_0$ be a function as in Proposition~\ref{prop:Q0}, so $\deg(Q_0)=|G|$ and $Q_0\circ\gamma = Q_0$ for all $\gamma \in G$.  
Then $\alpha \circ Q_0$ satisfies
these conditions also, for any $\alpha \in \PGL_2(K)$. We claim that $\alpha$ can be chosen so that $\alpha\circ Q_0 = f/g$, where $\deg(f) = |G| > \deg(g)$.
If $\deg(g_0) < |G|$, take $\alpha=\textmatrix 1001$, the identity map. If $\deg(g_0)=|G|$ and $\deg(f_0)<|G|$, then take $\alpha = \textmatrix 0110$, the reciprocal map.
Finally, if $\deg(f_0) = \deg(g_0) = |G|$, then there is $c \in K$ such that $\deg(f_0+cg_0) < |G|$, and $\textmatrix 0110 \textmatrix 1c01 \circ Q =
g_0/(f_0+cg_0)$ has the desired form.  

For the last statement, $\widetilde Q = \alpha \circ Q$ for $\alpha \in \PGL_2(K)$ by Proposition~\ref{prop:Q0}. The condition on the degrees of the denominators
forces $\alpha$ to have the form $\textmatrix a b 0 1$.
\end{proof}

Because the functions in Proposition~\ref{prop:QExists} are so central to this article, we give them a name.

\begin{definition} \label{def:quotient} A function $Q\in K(x)$ that satisfies the two conditions of Proposition~\ref{prop:QExists} is called
a {\it quotient map} for $G$.
\end{definition}

\begin{proposition} \label{prop:KQSigma}  If $Q(x)$ is a quotient map for $G$ then $K(Q(x))=\Sigma$, where $\Sigma$ is defined in~(\ref{eq:Sigma}). 
In particular, $K(x)/K(Q)$ is Galois, and
its automorphism group is isomorphic to~$G$.
\end{proposition}
\begin{proof} $Q \in \Sigma$ by the first part of the definition, so $K(Q)\subset \Sigma$.
Also, $\deg(Q)=|G|$ by the second part of the definition, so $[K(x):K(Q)]=\deg(Q)=|G|=[K(x):\Sigma]$. We conclude that $K(Q)=\Sigma$. Then $K(x)/K(Q)=K(x)/\Sigma$ is Galois, and its Galois group is isomorphic to $G$ by
Proposition~\ref{prop:Q0}.
\end{proof}

\subsection{Properties of quotient maps.} \label{sec:Qproperties}

\begin{proposition} \label{prop:H}
	If $H \subset G \subset \PGL_2(K)$ are finite subgroups and $Q_H$, $Q_G$ are quotient maps for these groups then $Q_G=h(Q_H)$ for a unique $h\in K(x)$,
	and $\deg(h)=|G|/|H|$.
\end{proposition}
\begin{proof} Let $\Sigma_H = \set{ u \in K(x) : \text{$u \circ \gamma = u$ for all $\g \in H$}}$ and define $\Sigma_G$ similarly.
By Proposition~\ref{prop:KQSigma},  $\Sigma_G=K(Q_G)$ and $\Sigma_H = K(Q_H)$.
Since $Q_G \in \Sigma_G \subset \Sigma_H=K(Q_H)$, $Q_G=h(Q_H)$ for some rational function~$h$. Since $Q_H$ is transcendental over $K$, $h$ is unique.
	$|G|=\deg(Q_G)=\deg(h(Q_H)) = \deg(h)\,\deg(Q_H) = \deg(h)|H|$ by Definition~\ref{def:quotient} and Lemma~\ref{lem:degf}. Thus, $\deg(h)=|G|/|H|$.
\end{proof}

The next proposition (especially {\it (i)}) illustrates that quotient maps have very strong arithmetic properties. 

\begin{proposition} \label{prop:Qproperties}
	Let $Q(x)\in K(x)$ be a quotient map for $G$, where $G \subset \PGL_2(K)$. Write $Q(x)=f(x)/g(x)$ where $f,g$ are relatively prime polynomials and $f$ is monic. 
Let $L/K$ be an extension field and let $x,y$ be independent transcendentals over $L$. Then 
\begin{enumerate} 
\item[{\it (i)}] {$f(y)-Q(x)g(y) = \prod_{\gamma \in G} \left(y - \gamma(x)\right)$.} \\
\item[{\it (ii)}] If $v_1,v_2 \in L$ and $Q(v_2)\ne \infty$ then $Q(v_1)=Q(v_2)$ if and only if $v_2 = \gamma(v_1)$ for some $\gamma \in G$. Consequently, if $w \in L$ then $Q^{-1}(w)$
is a $G$-orbit in $\cj L$. \\
\item[{\it (iii)}] If $w \in L$ and $\calO=Q^{-1}(w)$ is the corresponding orbit in $\cj L$, then $f(x)-w g(x) = \left(\prod_{v \in \calO} x-v\right)^{\mult(\calO)}$. \\
\item[{\it (iv)}] 
	$g(x) = a\prod_{v \in \calO_\infty, v \ne \infty} (x-v)^{\mult(\calO_\infty)}$, where $a \in K^\x$. Here
		$\calO_\infty=\set{\gamma(\infty) : \gamma \in G}$ and $\mult(\calO_\infty)=|H|$, 
		where $H=\{\g \in G : \g(\infty)=\infty\} = \{\textmatrix abcd \in G : c=0\}$.
\end{enumerate}
\end{proposition}

\begin{proof} {\it (i)}\ By Definition~\ref{def:quotient}, $\deg(f)=|G|>\deg(g)$. The left and right sides of the equation in~{\it (i)}, when regarded as polynomials in~$y$, 
	are both monic polynomials 
	of degree $|G|$  with coefficients in $\Sigma$, where $\Sigma$ is defined in (\ref{eq:Sigma}). Also, both vanish at $y=x$. 
	Since $[K(x):\Sigma]=|G|$ by Proposition~\ref{prop:Q0}, both are minimal polynomials for $x$ over $\Sigma$. Then each divides the other, so they are equal. \\
	{\it (ii)}\ If $v_2=\g(v_1)$ for some $\g\in G$, then $Q(v_2)=Q\circ \g(v_1)=Q(v_1)$, since $Q\circ \g=Q$. Now assume
	$Q(v_1)=Q(v_2)$. By hypothesis, this is finite, so $g(v_2)\ne 0$. By part {\it (i)},
$$f(v_2)-Q(v_1) g(v_2) = \prod_{\gamma \in G} v_2 - \gamma(v_1).$$
The left side vanishes since $Q(v_1)=Q(v_2)$. Thus, $v_2 = \gamma(v_1)$ for some $\gamma \in G$. \\
{\it (iii)}\ Let $v\in Q^{-1}(w)$. Set $x=v$ in the identity of part {\it (i)}
	to obtain that $f(y)-wg(y) = \prod_{\gamma \in G} y - \gamma(v)$, then
	apply 
Lemma~\ref{lem:orbitMultiplicity}. \\
	{\it (iv)}\   
	Let $F(x,y)=g(x)f(y)-f(x)g(y) \in K[x,y]$. Since $\deg(f)=|G|>\deg(g)$, this polynomial has degree~$|G|$ in each variable, and by~{\it (i)}, 
	$$F(x,y) =g(x)\prod_{\g\in G} \left(y-\g(x)\right).$$  
	Let
	$$ u(x) = \prod_{\textmatrix abcd \in G} (cx+d).$$
	(Since $G$ is projective, $u(x)$ is well-defined only up to a nonzero scalar multiple in $K^\x$.)
	Let $H=\{\textmatrix abcd \in G : c = 0\}$, so $|H|=\mult(\calO_\infty)$.
	Then $\deg(u)=|G|-|H|$, and
	$$F(x,y) = \frac{g(x)}{u(x)}  \prod_{\textmatrix abcd\in G} \left((cx+d)y-(ax+b)\right).$$
	Note that $(cx+d)y-(ax+b)$ has degree~1 in $y$ because $c$ or $d$ is nonzero; also it has degree~1 in $x$ because $c$ or $a$ is nonzero. Let
	$$P(x,y) = \prod_{\textmatrix abcd\in G} \left((cx+d)y-(ax+b)\right).$$
	This has degree $|G|$ in $x$ and in~$y$.
	Also, $P(x,y)$ is not divisible by any linear factor $rx+s \in K[x]$ with $r\ne0$, because $(cx+d)y+(ax+b)$ can be divisible by $rx+s$ only if
	it vanishes at $x=-s/r$, in which case $\textmatrix abcd \choose{-s/r}1 = \choose 00$, contradicting that $\textmatrix abcd$ is invertible.
	In particular, $P(x,y)$ is not divisible by any nonconstant factor $cx+d$ of $u(x)$, and so it is relatively prime to~$u(x)$.  
	Since $u(x)F(x,y)=g(x)P(x,y)$, $u(x)$ is relatively prime to $P(x,y)$,  $u(x)$ must divide~$g(x)$.
	$F(x,y)$ and $P(x,y)$ both have degree~$|G|$ in~$x$, therefore $\deg_x(g/u)=0$, \ie,  it is constant.
	To complete the proof, it remains only to prove that $u(x)$ is a constant multiple of $\prod_{v \in \calO_\infty, v \ne \infty} (x-v)^{|H|}$.

	Since $\textmatrix abcd^{-1} = \textmatrix d{-b}{-c}a$ in $\PGL_2(K)$ and $\g \to \g^{-1}$ is a bijection of~$G$, 
	$$ u(x) \equiv^*  \prod_{\textmatrix abcd \in G} (-cx + a),$$
	where the symbol $\equiv^*$ indicates ``up to a constant multiple in~$K^\x$''. Now
	$-cx+a\equiv^* 1$ if $c=0$, $-cx+a\equiv^* x-a/c$ if $c\ne 0$, and $a/c = \textmatrix abcd(\infty)\in \calO_\infty$.  
	Thus,
	$$u(x) \equiv^*  \prod_{\g \in G\setminus H} (x-\g(\infty)).$$
	Let $R$ be a complete set of coset representatives for $G/H$, excluding the identity coset, so $G \setminus H$ is the disjoint union of $rH$ for $r \in R$.  
	Writing $\gamma = r h$ with $r \in R$, $h\in H$ we have $\gamma(\infty)=rh(\infty)=r(\infty)$, so $\calO_\infty  \setminus \{\infty\} = \set{r(\infty):r \in R}$.
	Also, $\{r(\infty):r\in R\}$ are distinct, for $r(\infty)=r'(\infty)$ would imply $r^{-1}r'\in H$ and consequently $r'H=rH$.
	Thus, 
	\begin{equation*} u(x) \equiv^*  
		\prod_{r \in R} \prod_{h \in H} (x-rh(\infty)) 
		=  \prod_{r \in R} (x-r(\infty))^{|H|} 
		=  \prod_{v\in \calO_\infty,\ v\ne \infty}(x-v)^{|H|}.
	\end{equation*}
	Since $g$ is a constant multiple of $u$ and both have coefficients in~$K$, this proves the result.
\end{proof}

Parts {\it (ii)} and {\it (iv)} of Proposition~\ref{prop:Qproperties}
together imply:
\begin{equation}  \text{If $v_1,v_2 \in L \cup \{\infty\}$ then $Q(v_1)=Q(v_2)$ if and only if
$v_2 = \gamma(v_1)$ for some $\gamma \in G$.} \label{orbitStatement} 
\end{equation}

\begin{proposition} \label{prop:Qorbit} Let $G$ be a finite subgroup of $\PGL_2(K)$ and $Q$ a quotient map for~$G$. If $w \in \cj K \cup \{\infty\}$ then 
$Q^{-1}(w)$ is a $G$-orbit in $\cj K \cup \{\infty\}$. \end{proposition}

\begin{proof} First we show that $Q^{-1}(w)$ is nonempty. If $w = \infty$ then $\infty \in Q^{-1}(w)$. Now assume $w$ is finite, and let $v\in \cj K$
be a root of $g(x)w-f(x)$, where $Q=f(x)/g(x)$ and $f,g$ are relatively prime. If $g(v)=0$ then the equation $g(v)w-f(v)=0$ forces $f(v)=0$, contradicting
that $f,g$ are relatively prime. We conclude that $g(v) \ne 0$, so $w = f(v)/g(v) = Q(v)$ and $v \in Q^{-1}(w)$. The fact that $Q^{-1}(w)$ is a
$G$-orbit follows from Proposition~\ref{prop:Qproperties}{\it (ii)} if $w\in\cj K$ or Proposition~\ref{prop:Qproperties}{\it (iv)} if $w=\infty$.
\end{proof}

\subsection{Computation of quotient maps.} \label{sec:computeQ}

From the perspective of Galois theory,
quotient maps arise from invariant theory. We show in this section that they may also be computed 
by considering their zeros and poles.

\begin{theorem} \label{thm:computeQ}
Let $G$ be a finite subgroup of $\PGL_2(K)$. Let $\calO\subset \cj K$ be a $G$-orbit that does not contain $\infty$ and let $\mult(\calO)=|G|/|\calO|$ be its multiplicity. 
Let 
$$f_\calO(x)= \left(\prod_{v \in \calO} (x-v)\right)^{\mult(\calO)}\quad {\rm  and}\quad 
	g(x) = \left(\prod_{v \in \calO_\infty, v \ne \infty}(x-v)\right)^{\mult(\calO_\infty)}.$$ 
Then there is $w \in \cj K$ such that $f_\calO(x)/g(x) + w $ is a
quotient map for $G$.
\end{theorem}

\begin{proof}
Let $Q$ be a quotient map for $G$. By Proposition~\ref{prop:Qproperties}{\it (iv)}, $Q=f/g$ for some $f \in K[x]$. On replacing $Q$ by a constant multiple, we can assume that $f$ is monic.
Let $v \in \calO$ and $w=Q(v)\in \cj K$. 
By Proposition~\ref{prop:Qproperties}{\it (iii)}, $f(x)-wg(x)=f_\calO(x)$. Thus, $Q(x)=f(x)/g(x) = f_\calO(x)/g(x) + w$.
\end{proof}

\begin{example} \label{example:G3Q}
Let $G= \set{1,\b,\b^2}$ where $\b = \textmatrix 1{-1}10$. Then 
	$$\calO_\infty = \set{\infty,\beta(\infty),\beta^2(\infty)} = \set{\infty,1,0}.$$ 
	The denominator of $Q$ is therefore $x(x-1)$. To compute the numerator, select any $v \in \cj K \setminus \{0,1\}$ and compute its orbit $\calO$. 
If the characteristic is not~2, taking $v=-1$ gives
	$\calO=\{-1,2,1/2\}$, and $f_\calO(x)=(x+1)(x-2)(x-1/2)$.
The formula for $Q$ will be prettier if we add 3/2, so we take
	$$Q(x) = \frac{f_\calO(x)}{g(x)} + \frac 32 = \frac{(x+1)(x-2)(x-1/2)}{x(x-1)} + \frac 32 = \frac{x^3-3x+1}{x(x-1)}.$$
It turns out that this formula works for characteristic~2 as well. To see this, suppose that $K$ has characteristic~2 and let $\omega$ be a primitive
cube root of 1 in $\cj K$. Then $\{\omega\}$ is a $G$-orbit of multiplicity~3, and a quotient map is
$$\frac{(x-\omega)^3}{x(x-1)} + \omega = \frac{x^3+\omega x^2 + \omega^2 x + 1}{x(x-1)} + \omega = \frac{x^3+x+1}{x(x-1)}.$$
This equals $(x^3-3x+1)/(x(x-1))$ since $-3=1$ in characteristic~2. Thus, the formula $Q(x)=(x^3-3x+1)/(x(x-1))$ works for all fields.
\qed
\end{example}

\begin{example} \label{example:PGL2Q}
Consider $G=\PGL_2(\F_q)$.
Then $\calO_\infty = \F_q \cup \{\infty\}$, the multiplicity of this orbit is $(q^3-q)/(q+1)=q^2-q$, and
$\prod_{v\in \calO_\infty, v\ne \infty}(x-v)=\prod_{v\in\F_q}(x-v)=x^q-x$, so $g(x)=(x^q-x)^{q^2-q}$. 
To compute the numerator, select $v \in \F_{q^3} \setminus \F_q$.
Every element $\gamma(v)$ has the same degree as $v$, so the
orbit is contained in $\F_{q^3}\setminus \F_q$. 
Further, $\calO_v$ has full size, because all short orbits are contained
in $\F_{q^2} \cup \{\infty\}$ by Lemma~\ref{lem:short}.  Since $|\PGL_2(\F_q)|=q^3-q=|\F_{q^3}\setminus \F_q|$, it
follows that $\calO_v=\F_{q^3}\setminus \F_q$
and we may take the numerator to be
$$f(x) = \prod_{v \in \F_{q^3}}(x-v)/\prod_{v\in \F_q}(x-v) = (x^{q^3}-x)/(x^q-x).$$ 
We obtain $f(x)/g(x) = (x^{q^3}-x)/(x^q-x)^{q^2-q+1}$.  After working out formulas for $\inv(\tau)$, we decided to alter the definition to
$Q_G(x)=(x^{q^3}-x)/(x^q-x)^{q^2-q+1} + 1$ because that made the statement of our main theorem for $\PGL_2(\F_q)$ more aesthetic.

If instead we had selected $v \in \F_{q^2} \setminus \F_q$, the orbit would be $\calO=\F_{q^2} \setminus \F_q$, with multiplicity $|G|/|\calO|
	=(q^3-q)/(q^2-q)=q+1$.  Then $f_\calO(x)=\left((x^{q^2}-x)/(x^q-x)\right)^{q+1}$ and 
	a quotient map is $f_\calO(x)/g(x)=(x^{q^2}-x)^{q+1}/(x^q-x)^{q^2+1}$.
	It turns out that $Q_G$ and $f_\calO(x)/g(x)$ are equal. In fact, since $(x^{q^2}-x)^{q+1}=(x^{q^2}-x)^q(x^{q^2}-x)$,
\begin{eqnarray*} \frac{x^{q^3}-x}{(x^q-x)^{q^2-q+1}}-\frac{(x^{q^2}-x)^{q+1}}{(x^q-x)^{q^2+1}} &=& \frac{(x^{q^3}-x)(x^q-x)^q-(x^{q^3}-x^q)(x^{q^2}-x)}
{(x^{q^2}-x)^{q^2+1}} \\ 
&=& \frac{-x^{q^3+q}+x^{q^3+1}+x^{q^2+q}-x^{q^2+1}}{(x^q-x)^{q^2+1}} \\ 
&=& \frac{-x^{q^3}(x^q-x)+x^{q^2}(x^q-x)}{(x^q-x)^{q^2+1}} \\ 
&=& \frac{-x^{q^3}+x^{q^2}}{(x^q-x)^{q^2} } = -1.
\end{eqnarray*}
Then
\begin{equation}  
Q_G(x) = \frac{x^{q^3}-x}{(x^q-x)^{q^2-q+1}} + 1 = \frac{(x^{q^2}-x)^{q+1}}{(x^q-x)^{q^2+1}}.\label{eq:QGidentity}
\end{equation}
\qed
\end{example}

The above method to compute $f$ and $g$ by creating orbits finds quotient maps for most groups we considered. However, it did not work well 
when attempting to find a quotient map
for a cyclic group of $\PGL_2(\F_q)$ of order $\ell$ when $\ell \ge 3$ and $\ell| q+1$, as it is difficult to find an expression for $\prod_{v' \in \calO_v} (x-v')$.
Instead, we took advantage that such $G$ is conjugate over a quadratic extension to a diagonal subgroup, and for this group it is easy to find
$Q$. By composing $Q$ with the element of $\PGL_2(\F_{q^2})$ that establishes the conjugacy, one obtains an invariant function $Q_0$, but it is not rational
and the denominator has degree $|G|$. By applying an appropriate linear fractional transformation to $Q_0$, one can regain rationality and the property that the
degree of the denominator is smaller than the degree of the numerator.  
This lengthy computation was done in an earlier draft of the article. Fortunately, Xander Faber found a much easier method to compute this quotient map,
given in Proposition~\ref{prop:Qell}.

The following lemma describes how quotient maps are related when $G_1$ and $G_2$ are conjugate subgroups of
$\PGL_2(K)$, where $K$ is any field. 

\begin{lemma} \label{lem:conjugateGroups} Let $G_1$, $G_2$ be finite subgroups of $\PGL_2(K)$ that are conjugate to one another; \ie, there is $\a\in \PGL_2(K)$
such that $G_2=\set{\alpha\gamma \alpha^{-1}  : \gamma \in G_1 }$.
If $Q_1$ is a quotient map for $G_1$ then let $Q'=Q_1 \circ\alpha^{-1}$ and $k=Q'(\infty)$.
Let $\beta$ be any element of $\PGL_2(K)$ such that $\beta(k)=\infty$. Then $Q_2= \beta \circ Q_1 \circ \alpha^{-1}$ is a quotient map for $G_2$.
\end{lemma}

\begin{proof}  $Q_2 \circ \delta(x) = Q_2(x)$ for all $\delta \in G_2$ because  for $\gamma \in G_1$,
$$Q_2(\alpha\gamma \alpha^{-1} (x)) = \beta \circ Q_1(\gamma\alpha^{-1} (x)) = \beta \circ Q_1(\alpha^{-1} (x)) = Q_2(x).$$
Further, $\deg(Q_2)=\deg(Q_1)=|G|$ since linear fractional transformations do not affect the degree.
Finally, $Q_2(\infty) = \beta \circ Q'(\infty) = \beta(k)=\infty$, therefore the degree of the numerator of $Q_2$ exceeds the degree of the denominator. 
Thus, $Q_2(x)$ is a quotient map for $G_2$.
\end{proof}

\section{Artin invariant} \label{sec:invariant}

Let $K$ be a field and let $G$ be a finite subgroup of $\PGL_2(K)$.
In the previous section we defined a quotient map for $G$ to be a $G$-invariant function $Q(x)=f(x)/g(x)\in K(x)$ such that $\deg(f)=|G|>\deg(g)$, and we proved
existence and some properties.    In particular, if $\tau \in \cj K\cup\{\infty\}$ then $Q^{-1}(\tau)$ is a $G$-orbit in $\cj K\cup\{\infty\}$. 

\begin{definition} Let $\tau \in \cj K \cup \{\infty\}$.  
	If the $G$-orbit $Q^{-1}(\tau)$ has full size, \ie, $|Q^{-1}(\tau)|=|G|$, then we say that $\tau$ is {\it regular (with respect to~$Q$)}; otherwise it is irregular. 
\end{definition}

\begin{proposition}[Definition of Artin invariant] 
Let $\tau \in \cj K$ and $\s \in \Aut(\cj K/K)$ such that $\s(\tau)=\tau$. 
Let $v \in \cj K$ such that $Q(v)=\tau$. Then there is $\gamma \in G$ such that $\sigma(v)=\gamma(v)$.
If $\tau$ is regular, then the conjugacy class $\calC_\g =\{\d \g \d^{-1} : \d \in G\}$ is independent of the 
choice of $v\in Q^{-1}(\tau)$.
In that case, we write $\inv_Q(\tau,\sigma)=\calC_\g$, and we call $\inv_Q(\tau,\sigma)$ the {\it Artin invariant} of $\tau$ with respect to $Q$ and~$\sigma$.
If $\infty$ is regular, \ie, $Q^{-1}(\infty)$ has full size, then we define $\inv_Q(\infty,\s)=\calC_\textmatrix1001$.
\end{proposition}

\begin{proof} $Q^{-1}(\tau)$ is a $G$-orbit by Proposition~\ref{prop:Qorbit}. Since 
$Q\left(\sigma(v)\right)=\sigma\left(Q(v)\right)=\sigma(\tau)=\tau$, $v$ and $\sigma(v)$
are both in $Q^{-1}(\tau)$, therefore there is $\gamma \in G$ such that $\sigma(v)=\gamma(v)$. 
Now suppose that $\tau$ is regular. Then $|Q^{-1}(\tau)|=|G|$, so $\gamma$ is uniquely determined from $\sigma$ and~$v$. We claim that $\calC_\g$
depends only on $\sigma$ and $\tau$, and not on the choice of $v\in Q^{-1}(\tau)$. Indeed, suppose that $w\in Q^{-1}(\tau)$, and we will show that
$\s(w)=\a(w)$ where $\a \in \calC_\g$. 
There is $\d \in G$ such that $w=\d(v)$. Since the entries of $\delta$ are rational, $\s(w)=\s(\d(v)) = \d(\s(v))=\d(\g(v))=\d\g\d^{-1}(w)$. Here
$\d\g\d^{-1}\in\calC_\g$, as required.
\end{proof}

When $\tau=\infty$, then $Q^{-1}(\tau)=\{\g(\infty) : \g \in G\} \subset \{\infty\} \cup K$. If one defines $\s(\infty)=\infty$ for all $\s \in \Aut(\cj K/K)$,
then $\s(v)=v$ for all $v\in Q^{-1}(\tau)$. This is why it makes sense to define $\inv(\infty,\s)=\calC_{\textmatrix 1001}$ when $\infty$ is regular.
A benefit of this defintion is that it makes certain statements cleaner, for example Proposition~\ref{prop:cyclicCount}.

If $K=\F_q$ and $\sigma={\rm Frob}_q \in \Aut(\cj\F_q/\F_q)$
is the $q$-power Frobenius, then we will write $\inv_Q(\tau,q)$ instead of $\inv_Q(\tau,{\rm Frob}_q)$, or simply
$\inv(\tau)$ if $Q$ and $q$ are understood from the context.

For future reference, the definition of Artin invariant may be briefly summarized as follows when $K=\F_q$ and $\tau \in \F_q$ is regular:
\begin{equation}\text{If $Q(v)=\tau$, then $v^q=\gamma(v)$ for some $\g \in G$, and $\inv_Q(\tau,q)=\calC_\g.$} \label{eq:invariantFq}
\end{equation}
For arbitrary $K$, when $\tau \in \cj K$ is regular, $\sigma \in \Aut(\cj K/K)$, and $\s(\tau)=\tau$, the criterion is:
\begin{equation} \text{If $\tau=Q(v)$, then $\sigma(v)=\g(v)$ for some $\g \in G$, and $\inv_Q(\tau,\s)=\calC_\g$.}
\label{eq:invariantK}
\end{equation}

If $G$ is abelian, then $\calC_\g = \{\g\}$. In that case, we sometimes write $\inv_Q(\tau,\sigma)=\g$, instead of $\inv_Q(\tau,\sigma)=\calC_\g=\{\g\}$. 
\medskip

For every subgroup of $\PGL_2(\F_q)$ that we have investigated,
$\inv(\tau)$ can be described directly in terms of $\tau$, {\it e.g.}, involving Legendre symbols or other numerical invariants, without reference to $v$.
As a matter of notation, we often use a symbol
$[\tau/q]$ to denote these values that are directly computed from $\tau$. For instance, in Example~\ref{example:Klein},  we can define $[\tau/q] = 
\left(\jacobi\tau q,\jacobi{\tau-1}q\right)$ for $\tau \in \F_q \setminus \{0,1\}$, and then (\ref{eq:Klein}) describes $\inv(\tau,q)$ directly in terms of $[\tau/q]$.

\begin{proposition} \label{prop:main} Let $q$ be a prime power, $G$ a subgroup of $\PGL_2(\F_q)$, and $Q\in \F_q(x)$ 
a quotient map for $G$, as in Definition~\ref{def:quotient}. 
\begin{enumerate}
\item[{\it (i)}] For $\gamma\in G$, let 
\begin{equation} \text{$V_{\gamma,q} = \set{v \in \cj\F_q\setminus\calO_\infty : v^q = \gamma (v)}$ and
$V_{G,q} = \cup_{\b \in G} V_{\b,q}$.} \label{eq:VgammaDef}
\end{equation}
Then $V_{G,q}$ decomposes into exactly $q$ $G$-orbits, and $Q$ induces a bijection between these orbits and $\F_q $.  
\item[{\it (ii)}] If $v \in V_{G,q}$ and $|\calO_v|=|G|$, then there is a unique $\gamma\in G$ 
such that $v\in V_{\gamma,q}$.  
\item[{\it (iii)}] For each $\tau \in\F_q$ there is a conjugacy class $\calC \subset G$ 
such that $Q^{-1}(\tau) \subset \cup_{\gamma \in \calC} V_{\gamma,q}$.
If $\tau$ is regular (\ie, $|Q^{-1}(\tau)|=|G|$), then $\calC$ is uniquely determined and $\calC=\inv_Q(\tau,q)$.
\item[{\it (iv)}] If $G$ is abelian, then for each $\tau \in \F_q$ there is $\gamma \in G$ such that $Q^{-1}(\tau) \subset V_{\gamma,q}$. If in addition
$\tau$ is regular,  then $\g = \inv_Q(\tau,q)$.
\item[{\it (v)}] Suppose $\g\in G$ has order $t\ge 3$, and let $\calC_\g = \{\a\g\a^{-1} : \a\in G\}$. Then $t|(q+\kappa)$ for a unique $\kappa\in\{0,1,-1\}$,
	and the number of regular elements $\tau\in \F_q$ with $\inv_Q(\tau,q)=\calC_\g$ is
		exactly $|\calC_\g|(q+\kappa)/|G|$. In particular, if $G$ is abelian then there are exactly
		$(q+\kappa)/|G|$ regular elements $\tau\in\F_q$ with $\inv_Q(\tau,q)=\g$.
\end{enumerate}
\end{proposition}

\begin{proof}  {\it (i)}\ Let $v \in \cj\F_q \setminus \calO_\infty $ and $\tau = Q(v)\in\cj\F_q$. 
Since each preimage set $Q^{-1}(\tau)$ is a $G$-orbit by Proposition~\ref{prop:Qorbit}, 
\begin{eqnarray*} \tau\in\F_q &\iff& \tau = \tau^q \iff Q(v)=Q(v^q) \\
&\iff& \text{$v^q = \gamma (v)$ for some $\gamma \in G$} \iff v \in V_{G,q}.
\end{eqnarray*}
This shows 
\begin{equation} Q^{-1}(\F_q) = \cup_{\gamma \in G} V_{\gamma,q} = V_{G,q}.\label{eq:VGQ} \end{equation}
Since each preimage set $Q^{-1}(\tau)$ is a $G$-orbit 
and $\F_q$ has $q$ elements, $V_{G,q}$ partitions into exactly $q$ orbits.

{\it (ii)}\ If $v \in V_{G,q}$, then $v \in V_{\gamma,q}$ for some $\gamma \in G$, so $v^q=\gamma(v)$. If in addition the orbit of $v$ has full size,
then the elements $\gamma(v)$ for $\gamma \in G$ are distinct, so that $\gamma$ is uniquely determined from $v$ and $q$.

{\it (iii)}\ and {\it (iv)}\  Suppose $\tau \in \F_q$ and $v \in Q^{-1}(\tau)$.
By part {\it (i)}, which we have already proved, 
there is $\gamma \in G$ such that $v^q = \gamma(v)$. 
Let $\calC = \set{\alpha \gamma \alpha^{-1} : \alpha \in G}$, the conjugacy class of $\gamma$. We claim that
$\calO_v \subset \cup_{\beta \in \calC} V_{\beta,q}$.  To see this, let $w = \alpha(v) \in \calO_v$, where $\a \in G$. 
Since the entries of $\a$ are in $\F_q$, 
$$w^q = \left(\a(v)\right)^q = \a(v^q) = \a \gamma(v) = \a \gamma \a^{-1} (w),$$
therefore $w \in V_{\beta,q}$ where $\beta = \a \gamma \a^{-1} \in \calC$. This proves the claim. Now suppose $\tau$ is regular.
	For any $v\in Q^{-1}(\tau)$ there is $\b\in\calC$ such that $v^q=\beta(v)$. Then
	$\inv_Q(\tau,q)=\calC_\beta$ by~(\ref{eq:invariantFq}). Since $\beta\in\calC$, $\calC=\calC_\b$.

{\it (v)}\ By~{\it (iii)}, the number of regular $\tau\in\F_q$ with $\inv_Q(\tau,q)=\calC_\g$ is the number of full-sized $G$-orbits in
	$\cup_{\b\in\calC_\g} V_{\b,q}$. By Proposition~\ref{prop:counting}, this number is $|\calC_\g|(q+\kappa)/|G|$.
\end{proof}

Recall that if $G \subset \PGL_2(K)$, there was some choice in the definition of $Q$, as one could change it to $aQ+b$, where $a,b \in K$ and $a\ne0$.
Since $\tau = Q(v) \iff a\tau + b = (aQ+b)(v)$, the set of preimages of $\tau$ under $Q$ is the same as the set of preimages of $a\tau + b$ under $aQ+b$.
In particular, $\tau$ is regular with respect to $Q$ iff $a\tau+b$ is regular with respect to $aQ+b$, and in that case
\begin{equation} \inv_{a Q+b}(a\tau+b) = \inv_Q(\tau). \label{eq:aQb} \end{equation}
We select $a,b$ so that the invariants, when expressed in terms of $\tau$, have simple and natural expressions.

If two finite subgroups of $\PGL_2(K)$ are conjugate to one another by a rational linear fractional transformation, then their Artin invariants are 
essentially equivalent, as shown below.
Thus, we are free to normalize groups via rational conjugation when possible.  If $\s\in\Aut(\cj K/K)$ then we define $\s(\infty)=\infty$.

\begin{lemma}  \label{lem:conjugateInv}
Suppose that $G_1,G_2$ are finite subgroups of $\PGL_2(K)$ that are conjugate to one another; that is, there is $\alpha \in \PGL_2(K)$ such that
$G_2=\set{\alpha\gamma \alpha^{-1}  : \gamma \in G_1 }$. If $Q_1$ is a quotient map for $G_1$ then let $Q_2 = \b \circ Q_1 \circ \a^{-1}$
	be a quotient map for $G_2$, as in Lemma~\ref{lem:conjugateGroups}. Then $\tau \in \cj K\cup\{\infty\}$ is regular with respect to $Q_1$ iff $\beta(\tau)$ is regular
	with respect to $G_2$. Further, if $\tau\in \cj K\cup\{\infty\}$ is regular, $\s \in \Aut(\cj K/K)$, and $\s(\tau)=\tau$ then
$$\inv_{Q_2}(\beta(\tau),\s) = \a\, \inv_{Q_1}(\tau,\s)\, \a^{-1}.$$
\end{lemma}

\begin{proof} $Q_1(v) = \tau$ iff $\beta\circ Q_1 (v) = \beta(\tau)$
iff $Q_2(\alpha(v))=\beta(\tau)$.  Therefore, $Q_2^{-1}(\beta(\tau))=\a (Q_1^{-1}(\tau))$.
It follows that $Q_1^{-1}(\tau)$ has full size iff $Q_2^{-1}(\beta(\tau))$ has full size, so $\tau $ is regular wrt $Q_1$ iff $\beta(\tau)$ is regular
	wrt $Q_2$. In that case, $\g\in \inv_{Q_1}(\tau,\s)$ iff there is $v\in Q_1^{-1}(\tau)$ with $\s(v) = \gamma(v)$ iff there is $w=\a(v) \in Q_2^{-1}(\b(\tau))$ with
	$\sigma(w)=\a(\sigma(v))=\a\g(v)=\a\g\a^{-1}(w)$ iff $\a\g\a^{-1} \in \inv_{Q_2}(\beta(\tau))$. 
	Note that this proof is valid even when $\tau=\infty$ or $\b(\tau)=\infty$.
\end{proof}

The next lemma shows that if $H\subset G$ are finite subgroups
of $\PGL_2(K)$ then their Artin invariants
are closely related. If $\d\in H$, let $\calC_{\d,H} = \{\g\d\g^{-1} : \g \in H\}$.
Then $\calC_{\d,H} \subset \calC_{\d,G}$.
Let $Q_G$, $Q_H$ be quotient maps for $G$ and $H$.
By Proposition~\ref{prop:H}, there is a unique rational function 
$h\in K(x)$ of degree~$|G|/|H|$ such that $Q_G = h\circ Q_H$.

\begin{lemma} \label{lem:Hinv} Let $G,H,Q_G=h(Q_H)$ be as above.
	Suppose $h(\tau)$ is regular with respect to~$G$, where $\tau \in \cj K\cup\{\infty\}$.
Then $\tau$ is regular with respect to~$H$, and for any $\s\in\Aut(\cj K/K)$
	such that $\s(\tau)=\tau$, there is $\d\in H$ such that
	$$ \inv_{Q_H}(\tau,\s)= \calC_{\d,H}, \qquad \inv_{Q_G}(h(\tau),\s) = \calC_{\d,G}.$$
\end{lemma}
\begin{proof} Let $V=Q_H^{-1}(\tau)$; this is an $H$-orbit by Proposition~\ref{prop:Qorbit}.
	Let $v\in V$. Then $\s(v)\in V$, and there is $\d \in H$ such that $\s(v)=\d(v)$.
	Since $Q_G(v)=h(Q_H(v))=h(\tau)$, $v$ is in the $G$-orbit $Q_G^{-1}(h(\tau))$.
	By hypothesis, $h(\tau)$ is regular with respect to $Q_G$, therefore 
	$\g(v)$ for $\g \in G$ are distinct.
	Then $V=\{\g(v) : \g \in H\}$ has $|H|$ distinct elements, so $\tau$ is regular with respect to $Q_H$.
	Since $Q_H(v)=\tau$, $Q_G(v)=h(\tau)$, and $\s(v)=\d(v)$, (\ref{eq:invariantK}) implies
	$\inv_{Q_H}(\tau,\s) = \calC_{\d,H}$ and $\inv_{Q_G}(h(\tau),\s) = \calC_{\d,G}$.
\end{proof}

\part{Small subgroups of $\PGL_2(K)$} \label{part:K}

This part considers finite subgroups that are contained in $\PGL_2(K)$ for any $K$. Sections~\ref{sec:known} and~\ref{sec:K}
contain known examples, and Sections~\ref{sec:three} and \ref{sec:six} contain new ones.

\section{An example related to Dickson polynomials} \label{sec:known}

Anyone who has studied Dickson polynomials is probably familiar with the lemma that $v \mapsto v + 1/v$
gives a surjective map from $\mu_{q-1} \cup \mu_{q+1}$ onto $\F_q$.
This lemma appears in a 1961 article by Brewer \cite{Brewer}, and was probably known earlier.
This section examines that lemma from the Artin invariant perspective.  

Let $K$ be any field, $c \in K^\x$,  and
$G = \set{\textmatrix 1001,\textmatrix 0c10 }\subset\PGL_2(K)$.
The $G$-orbits are $\calO_v = \set{v,c/v}$ for $v\in\cj K \cup\{\infty\}$, and the short orbits are $\{\sqrt c\}$ and $\{-\sqrt c\}$.
A quotient map is $Q(x)=x+c/x$.
The short orbits $\{\sqrt c\}$ and $\{-\sqrt c\}$ are sent under $Q$ to $2\sqrt c $ and $-2\sqrt c $,  respectively; 
these are the irregular elements of $\cj K \cup\{\infty\}$.

\begin{proposition} \label{prop:0c10}  Let $G$ and $Q$ be as above. 
	Let  $\tau \in \cj K$ be regular (equivalently, $\tau^2-4c\ne0$) and let $\s \in \Aut(\cj K/K)$
	such that $\s(\tau)=\tau$. If char$(K)\ne 2$, then 
	$\inv(\tau,\sigma)= \textmatrix 0c10^{(1-A)/2}$, where $A=\sigma\left(\sqrt{\tau^2-4c\,}\,\right)/\sqrt{\tau^2-4c\,}\in \{1,-1\}$. 
	If char$(K)=2$, let $X$ be a root of $X^2+X+c/\tau^2=0$. Then $\s(X)=X+j$, where $j\in \F_2$, and $\inv(\tau,\s)=\textmatrix 0c10^j$.
\end{proposition}
\begin{proof} The solutions to $Q(v)=\tau$ are the roots of $x^2-\tau x +c$. Denote these roots by $v$ and $v'$. Let $\g=\inv(\tau,\s)\in G$, so $\s(v)=\g(v)$. 
	If char$(K)\ne 2$ then $\{v,v'\}=\{(\tau\pm \sqrt{\tau^2-4c})/2\}$.
	If $\s(v)=v$ then $\s(\sqrt{\tau^2-4c})=\sqrt{\tau^2-4c}$ and $\g=1$; otherwise $\s(v)=v'$, $\s(\sqrt{\tau^2-4c})=-\sqrt{\tau^2-4c}$, $\g=\textmatrix 0c10$.
	This proves the proposition when the characteristic is not~2.
	If char$(K)=2$, then $v/\tau$, $v'/\tau$ are the two solutions to $X^2+X=c/\tau^2$, and $v'/\tau = (v/\tau)+1$. 
	If $\s(v)=v$ then $\s(X)=X$ and $\g=\textmatrix 1001$. If $\s(v)\ne v$ then $\s(v)=v'$, $\s(X)=X+1$ and $\g=\textmatrix 0c10$. The result follows.
\end{proof}

Consider the case $K=\F_q$ and $\s(x)=x^q$. If $\tau\in \cj \F_q$ then $Q^{-1}(\tau)$ is a $G$-orbit, say $\{v,c/v\}$, and $\tau^q=\tau$ iff $v^q \in \{v,c/v\}$.
When $c=1$, then $v^q\in \{v,1/v\}$ iff $v\in \mu_{q-1}\cup \mu_{q+1}$, and one obtains the result mentioned above
that $v\mapsto v+1/v$ gives a surjective map of $\mu_{q-1}\cup \mu_{q+1}$ onto $\F_q$.
Proposition~\ref{prop:0c10} in this case is due to Brewer \cite{Brewer} for $q$ odd, and Dillon-Dobbertin \cite{DD} when $q$ is even.
In the case where $q$ is even, the element $j$ in Proposition~\ref{prop:0c10} is equal to $\Tr_{\F_q/\F_2}(c/\tau^2)$.

\begin{proposition} \label{prop:BDD}
{\bf (Brewer \cite{Brewer} for $q$ odd; Dillon and Dobbertin \cite{DD} for $q$ even).}
Let $q$ be a prime power and $Q(x) = x + c/x$, where $c \in \F_q^\x$.  Every $\tau \in \F_q$ may be written as $\tau = Q(v)$, where $v\in\cj\F_q^\x$
and where
	$v^q = v$ or $v^q = c/v$. Conversely, if $v^q=v$ or $v^q=c/v$ then $Q(v) \in \F_q$.
Let $\tau=Q(v)\in\F_q$.
If $\tau^2=4c$ then $v=\pm \sqrt c$ and $v^q=v=c/v$.  If $\tau^2 \ne 4c$, then
\begin{enumerate}
\item for $q$ odd: $v^q=v \iff \jacobi{\tau^2-4c}q = 1$;\qquad $v^q=c/v \iff \jacobi{\tau^2-4c}q = -1$;
\item for $q$ even: $v^q=v \iff \Tr_{\F_q/\F_2}(c/\tau^2)=0$;\qquad $v^q=c/v \iff \Tr_{\F_q/\F_2}(c/\tau^2)=1$.
\end{enumerate}
\end{proposition}

Brewer proved Proposition~\ref{prop:BDD}(1) and used it to compute the number of $\F_p$-rational points
of a curve $y^2 = D_n(x)$ over a prime field $\F_p$, where $D_n(x)$ is a Dickson polynomial, determined by the property that $D_n(x+1/x)=x^n+1/x^n$.
These point-counting formulas were applied to determine some character sums.  
Dillon and Dobbertin proved Proposition~\ref{prop:BDD}(2) and used it to show that for $q$ even, $D_n(S_j) \subset S_j \cup \{0\}$, where 
$S_j= \set{\tau\in\F_q^\x : \Tr_{F_q/\F_2}(1/\tau)=j} $ for $j \in \{0,1\}$.
Analogues of the result of Dillon and Dobbertin for the case of odd characteristic are formulated in \cite{Permutation}.

\medskip
\noindent{\bf Remark.}\ The large difference in behavior between odd or even characteristic in Proposition~\ref{prop:BDD}
is attributable to the fact that the transformation $x\mapsto c/x$ has a unique fixed point in the algebraic closure in characteristic~2, but two fixed points in odd characteristic.
Consequently, $G$ is conjugate to the unipotent group $\{\textmatrix 1001,\textmatrix 1101\}$ in characteristic~2.
Unipotent groups are related to additive or linear
maps (see Section~\ref{sec:unipotent}), consistent with the appearance of a trace map in Proposition~\ref{prop:BDD}(2).
This phenomenon will revisit us in Section~\ref{sec:three}, when we study a group of order~3.  
\qed
\medskip

As shown in \cite[Section 9]{Wilson-like}, Proposition~\ref{prop:BDD}
can be used to obtain a quick proof of Legendre symbol formulas $\jacobi 2q$, $\jacobi 3q$, and $\jacobi 5q$. 

\section{Kummer and Klein examples} \label{sec:K}

This section shows that two groups that were considered in the introduction for finite fields generalize to an arbitrary field $K$.
\bigskip

\noindent {\it Kummer extensions.}\  (See Example~\ref{example:Kummer}.)
Suppose $K$ contains the primitive $n$th roots of unity, so in particular $p\nmid n$ if $K$ has finite characteristic~$p$.
Let $G= \set{\textmatrix a001 : a^n=1}\subset \PGL_2(K)$. Then $Q(x)=x^n$ is a quotient map. If $\tau \in K^\x$ then $\tau$ is regular, $Q^{-1}(\tau)=
\set{a \tau^{1/n} : a^n=1}$, and if $\s \in \Aut(K(\tau^{1/n})/K)$ then $\inv_Q(\tau,\sigma)=\textmatrix a001$, where $a=\sigma(\tau^{1/n})/\tau^{1/n}\in \mu_n$.

\bigskip

\noindent{\it Klein group.}\ (See Example~\ref{example:Klein} and \cite{Structure}.) 
Let $G= \set{\textmatrix 1001,\textmatrix {-1}001,\textmatrix 0b10, \textmatrix 0{-b}10}\subset \PGL_2(K)$, where $b$ is a fixed element of $K^\x$.  Assume the
characteristic of $K$ is not~2. The short orbits are
$\calO_\infty=\{\infty,0\}$, $\{\pm\sqrt b\}$, and $\{\pm\sqrt{-b}\}$. A quotient map is $Q(x)=(x+b/x)^2/4$, and $b,0,\infty$ are irregular.
Let $\tau \in \cj K$. Then $Q^{-1}(\tau) = \set{\pm \sqrt\tau \pm \sqrt{\tau-b}}$, and if $v=\sqrt\tau+\sqrt{\tau-b}$
(for some choices of square root), then $b/v=\sqrt\tau-\sqrt{\tau-b}$. Suppose that $\tau$ is regular ({\it i.e.}, $\tau \not \in \{0,b\}$), 
$\s \in \Aut(\cj K/K)$, and $\s(\tau)=\tau$.
Let $A=\sigma(\sqrt\tau)/\sqrt\tau$
and $B=\sigma(\sqrt{\tau-b})/\sqrt{\tau-b}$. Then $\inv(\tau,\sigma) = \textmatrix A001 \textmatrix 0b10^{(1-AB)/2}$.

\section{Artin invariant for a transformation group of order~3} \label{sec:three}

Consider 
$G_3= \set{I,\beta,\beta^2}$ where 
$$\beta = \textmatrix{1\,}{-1}{1\,}{\,\,0}\in\PGL_2(K)$$
and $K$ is any field.
The $G_3$-orbits of $\cj K \cup \{\infty\}$ are $\calO_\infty = \{\infty,1,0\}$ and
$$\calO_v = \set{v, 1-1/v, 1/(1-v)}, \quad
\text{$v\in \cj K \setminus \{0,1\}$.}$$

$G_3$ is normalized by the map $\rho = \textmatrix 0110$ that takes $v$ to $1/v$:
$\rho\beta\rho=\beta^2=\beta^{-1}$.  This feature will appear in our
analysis. (See Lemma~\ref{lem:G3invProperty}.)

The characteristic-3 case turns out to be very different from other characteristics.  
In fact, our result for characteristic~3 is strikingly similar to the theorem 
of Dillon and Dobbertin (see Prop.~\ref{prop:BDD}(2)).
This phenomenon will be explained at the end of this section.
\newpage

\subsection{Short orbits, quotient map, irregular elements.} \label{sec:G3Q}

\begin{lemma} \label{lem:G3short} The short orbits of $G_3$ (that is, the orbits with fewer than three
elements) are $\set{-1}$ in characteristic~3, or $\set{-\omega}$
and $\set{-\omega^2}$ in characteristic different from~3, where $\omega$
is a primitive cube root of unity in $\cj K$.
\end{lemma}

\begin{proof}
By Lemma~\ref{lem:short}, $v$ is in a short orbit if 
and only if $v=\beta(v)$ or $v=\beta^{-1}(v)$, or equivalently,
$v^2-v+1=0$. In characteristic~3, this factors as $(v+1)^2=0$, so
$\{-1\}$ is the only short orbit.  If the characteristic is not three,
then $v^2-v+1=0 \iff v \in \set{-\omega,-\omega^2}$.
\end{proof}

\begin{lemma} \label{lem:G3Q}
$Q_3(x) = (x^3-3x+1)/(x(x-1))$ is a quotient map for $G_3$ over any field. 
The set of irregular elements of $\cj K \cup \{\infty\}$ is $\{0\}$ if char$(K)=3$ and
$\{-3\omega ,-3\omega^2\} $ if char$(K)\ne3$, where $\omega$ is a 
primitive cube root of unity in~$\cj K$. That is, the irregular elements are the roots
of $\tau^2-3\tau+9=0$.
\end{lemma}
\begin{proof} The formula for the quotient map was computed in Example~\ref{example:G3Q}.
The irregular points are the images under $Q_3$ of the short orbits. These are
$Q_3(-1)=0$ in char.~3, and $Q_3(-\omega)=-3\omega$,
$Q_3(-\omega^2)=-3\omega^2$ in char.~$\ne 3$. 
\end{proof}

\begin{lemma} \label{lem:reciprocal} $Q_3(1/x) = 3-Q_3(x)$.
\end{lemma}

\begin{proof} This is a simple computation. \end{proof}

\begin{lemma} \label{lem:G3invProperty} If $\tau \in \cj K$ is regular then so is $3-\tau$, and for any $\s \in \Aut(\cj K/K)$ 
such that $\s(\tau)=\tau$ we have
$\inv_{Q_3}(3-\tau,\sigma)=\inv_{Q_3}(\tau,\sigma)^{-1}$. 
\end{lemma}
\begin{proof} By Lemma~\ref{lem:reciprocal}, $\tau=Q_3(v) \iff 3-\tau=Q_3(1/v)$. The short elements of $\cj K$ are closed under reciprocal,
therefore $\tau$ is regular iff $3-\tau$ is regular. Let $\rho = \textmatrix 0110$. Now 
\begin{eqnarray*}
\inv(\tau,\s)=\beta^j &\iff& \sigma(v)=\beta^j(v) \\
&\iff& \sigma(1/v)=1/\sigma(v)=\rho\sigma(v)=\rho\beta^j(v)=\beta^{-j}\rho(v)=\beta^{-j}(1/v).
\end{eqnarray*}
Since $Q_3(1/v)=3-\tau$ and $\s(1/v)=\b^{-j}(1/v)$, this shows that $\inv(3-\tau,\sigma)=\b^{-j}=\inv(\tau,\s)^{-1}$.
\end{proof}

If $\tau\in \F_q$ is regular then $\inv(\tau)\in G_3$
is the unique element $\gamma \in G_3$ such that $v^q=\gamma(v)$ for any (hence every) $v \in Q_3^{-1}(\tau)$.   Recalling that $\beta = \textmatrix 1{-1}10$,
and noting that $v \not \in \calO_\infty = \{\infty,0,1\}$, 
\begin{equation}
\inv(\tau) = \begin{cases} \beta^0 & \text{iff $v^{^q-1}=1$,} \\ 
\beta  & \text{iff $v^{q+1}-v+1=0$,} \\ 
\beta^{-1} & \text{iff $v^{q+1}-v^q+1=0$.} \end{cases} \label{eq:type}
\end{equation}

\subsection{Explicit description of $\inv(\tau,\sigma)$ when char$(K)\ne 3$.}

We wish to express $\inv(\tau,\sigma)$ 
purely in terms of $\tau$, without reference to~$v$.
Our method is to solve for $v$ in terms of $\tau$
and then determine how the Galois group permutes the solutions.
The equation relating $v$ and $\tau$ is cubic in $v$,
namely, 
\begin{equation} v^3-3v+1-\tau v(v-1)=0. \label{eq:cubic} \end{equation}
There are well-documented ways to explicitly solve a cubic dating back to the 1500's, however
these fail in characteristic~3.  This section assumes char$(K)\ne3$,
and the characteristic-3 case will be considered
in Section~\ref{sec:G3char3}.  
Let $\omega$ denote a {\it fixed} primitive cube root of~1 in $\cj K$. 
Then $\omega^2+\omega+1=0$.

The first step to solve a cubic $x^3+Ax^2+Bx+C$ is to make a change of
variables $y=x+A/3$ so as to eliminate the $x^2$ term.  Writing
$y^3 + Dy + E =0$, substitute $y=z+k/z$ to obtain  
$$z^3 + k^3/z^3 + (3k+D)z + (3k^2+Dk) z^{-1} + E = 0.$$
By setting $k=-D/3$, the $z$ and $z^{-1}$ terms drop out:
$$z^6 + E z^3 - D^3/27 = 0.$$
Then $z^3 = \left(-E \pm \sqrt{E^2+4D^3/27}\right)/2$.
The right side has two possible values, so $z$ has six possible values,
however $y=z+k/z$ turns out to have only three possible values.

To solve the cubic (\ref{eq:cubic}), 
it is convenient to let 
$$\hat\tau = \tau/3,\qquad R = \hat\tau^2-\hat\tau+1=(\hat\tau+\omega)(\hat\tau+\omega^2).$$
Substitute $y=v-\hat \tau$ to obtain
$$y^3 -3Ry - (2\hat\tau-1)R = 0.$$
Next, substitute $y = z + R/z$ to obtain
$$z^3 + R^3/z^3 - (2\hat\tau-1)R = 0.$$
Then $z^6-(2\hat\tau-1)Rz^3+R^3=0$, so $(z^3/R)^2-(2\hat\tau-1)(z^3/R)+R=0$ and
$$z^3/R = (1/2)\left( 2\hat\tau-1 \pm \sqrt{(2\hat\tau-1)^2-4R}\right).$$
Note that $(2\hat\tau-1)^2-4R=-3$ and $\omega = (-1+\sqrt{-3})/2$ for
an appropriate choice of $\sqrt{-3}$.  Thus, one choice of $z$ satisfies
$$z^3 = R\cdot(\hat\tau + \omega).$$
Fix $\lambda,\mu \in
\cj K$ such that 
\begin{equation} \text{$\lambda^3 = \tau/3+\omega$ and $\mu^3 = \tau/3+\omega^2$.} 
\label{eq:lambdamu} \end{equation}
Then $\lambda\mu$ is a cube root of $R$, so one solution for $z$ is
$z = \lambda^2 \mu$.  Then $v=y+\hat\tau=\hat\tau + z+R/z
=\hat\tau + \lambda^2\mu + \lambda^3 \mu^3/(\lambda^2\mu)
= \hat\tau + \lambda^2 \mu + \lambda \mu^2$.
If we set $\lambda'=\omega^j \lambda$ for $j\in\Z/3\Z$ and use $\lambda'$ in the above
construction instead of $\lambda$, then we arrive at a solution
\begin{equation} v_j = \tau/3 + \omega^{-j} \lambda^2 \mu + \omega^{j} \lambda \mu^2.
\label{eq:vj} \end{equation}

\begin{proposition} \label{prop:vjRel} Suppose char$(K)\ne 3$, let $\tau\in \cj K$
	and $\l,\mu,v_j$ as in (\ref{eq:lambdamu}) and~(\ref{eq:vj}).
Then $Q_3(v_j) = \tau$, where $Q_3(x)=(x^3-3x+1)/(x^2-x)$.  Also,
\begin{equation}
v_1 = 1-1/v_0, \quad v_2 = 1-1/v_1.  \label{eq:vjRel}
\end{equation} 
\end{proposition}

\begin{proof} That $Q_3(v_j)=\tau$ was proved above. By Proposition~\ref{prop:Qorbit}, $v_0$, $v_1$, and $v_2$ 
belong to the same $G_3$-orbit. Either $v_1 = 1-1/v_0$, in which case $v_0 v_1 = v_0-1$,
or $v_1=1/(1-v_0)$, in which case $v_0v_1 = v_1-1$.
To see which of these holds, we compute $v_0 v_1$. As before, let $\hat\tau = \tau/3$.  Then
\begin{eqnarray*} v_0 v_1 &=& (\hat\tau + \lambda^2 \mu + \lambda \mu^2)(\hat\tau + \omega^2 \lambda^2\mu + \omega \lambda \mu^2) \\
&=& (\hat\tau + \lambda \mu (\lambda + \mu))(\hat\tau + \lambda\mu (\omega^2 \lambda + \omega \mu)) \\
&=& \hat\tau^2 + \hat\tau \lambda\mu (\lambda + \mu + \omega^2 \lambda+\omega \mu )
     + \lambda^2 \mu^2 ( \omega^2 \lambda^2 + \omega^2\lambda\mu + \omega\lambda\mu + \omega\mu^2) \\
&=& \hat\tau^2 + \hat\tau \lambda\mu (-\omega\lambda - \omega^2 \mu )
     + \lambda^2 \mu^2 ( \omega^2 \lambda^2 - \lambda\mu + \omega\mu^2) \\
&=& \hat\tau^2 - \hat\tau (\omega \lambda^2\mu + \omega^2 \l\mu^2 )
     + \lambda^3 \omega^2\lambda\mu^2 -\lambda^3\mu^3 + \mu^3\omega\lambda^2\mu \\
&=& \hat\tau^2 - \hat\tau (\omega \lambda^2\mu + \omega^2 \l\mu^2 )
     + (\hat\tau+\omega) \omega^2\lambda\mu^2 -(\hat\tau+\omega)(\hat\tau+\omega^2) + (\hat\tau+\omega^2)\omega\lambda^2\mu  \\
&=& \hat\tau^2 +\lambda\mu^2 - (\hat\tau^2-\hat\tau+1) + \lambda^2\mu \\
&=& \hat\tau + \lambda^2\mu + \lambda\mu^2 - 1 = v_0-1.
\end{eqnarray*}
This computation shows that $v_1=1-1/v_0$, \ie, $v_1=\beta(v_0)$, where $\beta=\textmatrix1{-1}10$.   
Since $\set{v_0,v_1,v_2}$ is an orbit, it must be that $v_2 = \beta(v_1)$.
\end{proof}

It is interesting to note what happens when $\tau$ is irregular. 
Then, as shown in Section~\ref{sec:G3Q}, $\tau \in \{-3\omega,-3\omega^2\}$, therefore
$\lambda\mu=0$ and $v_j = \tau/3$ for all $j$.

The labeling of $v$'s depends on the choices for $\omega$, $\lambda$ and $\mu$.
However, for any such choice, $\beta(v_j)=1-1/v_j = v_{j+1}$ for each $j$. In other words,
\begin{equation} \beta^k(v_j) = v_{j+k}\quad\text{for $j,k\in\Z/3\Z$.} \label{eq:betavj} \end{equation}
The next theorem describes $\inv_{Q_3}(\tau,\sigma)$ directly in terms of $\sigma$ and $\tau$.

\begin{theorem} \label{thm:G3charneq3} Suppose char$(K)\ne 3$, and let $\s \in \Aut(\cj K/K)$.  
If $\tau \in \cj K \setminus\{-3\omega,-3\omega^2\}$ and $\s(\tau)=\tau$, then
$\inv_{Q_3}(\tau,\sigma)=\b^\ell$ where $\ell$ is determined from $\tau$ as follows.\\
{\it (i)}\ Let $\z \in \cj K$ such that 
	$$\z^3=\frac{\tau+3\omega^2}{\tau+3\omega}.$$
Then $\s^2(\z)/\z = \omega^\ell$.  \\
{\it (ii)}\ 
	If $K=\F_q$ and $\s(x)=x^q$, where $3\nmid q$, then 
$$\left(\frac{\tau+3\omega^2}{\tau+3\omega}\right)^{(q^2-1)/3} = \omega^\ell.$$
\end{theorem}

\begin{proof} 
	Since $Q(v_j)=\tau$ and $\inv(\tau,\s)=\beta^\ell$, (\ref{eq:invariantK}) and (\ref{eq:betavj}) imply 
	$\sigma(v_j)=\beta^{\ell}(v_j)=v_{j+\ell}$, hence
	$$\sigma^2(v_j)=\beta^{2\ell}(v_j)=v_{j+2\ell}.$$

	Let $\l,\mu$ be as in (\ref{eq:lambdamu}) and $\z_0=\mu/\l$. 
	Then $\z_0^3=(\tau+3\omega^2)/(\tau+3\omega)=\z^3$, therefore $\z=\omega^i\z_0$
	for some $i\in \Z/3\Z$.  Since $\omega$ is defined over at most a quadratic
	extension of $K$, $\sigma^2(\omega)=1$ and hence $\s^2(\z)/\z=\s^2(\z_0)/\z_0$. 
	Now $\sigma^2(\z_0)/\z_0$
	is a cube root of unity, because its cube is equal to $\sigma^2(\z_0^3)/\z_0^3=1$.
	Let $\omega^k = \sigma^2(\z_0)/\z_0$. 

	By~(\ref{eq:vj}),
	$$v_j/\l^3= \tau/(3\l^3) + \omega^{-j} \mu/\l + \omega^j(\mu/\l)^2 =
	\frac\tau{\tau+3\omega} + \omega^{-j} \z_0 + \omega^j \z_0^2.$$
	Then
$$ \sigma^2(v_j/\l^3) =
	\frac\tau{\tau+3\omega} + \omega^{-j} \omega^k\z_0 + \omega^j \omega^{2k}\z_0^2
	= v_{j-k}/\l^3,$$
	which implies $\sigma^2(v_j)=v_{j-k} = v_{j+2k}$. 
	We have shown $\s^2(v_j)=v_{j+2k}=v_{j+2\ell}$. Since $v_0,v_1,v_2$
	are distinct when $\tau$ is regular, it follows that $k\equiv \ell \pmod3$.
This proves {\it (i)}.

	If $K=\F_q$ and $\s(x)=x^q$, then
	$$ \omega^\ell = \s^2(\z)/\z =  \z^{q^2-1} = (\z^3)^{(q^2-1)/3} 
	= \left( \frac{\tau+3\omega^2}{\tau+3\omega}\right)^{(q^2-1)/3},$$
	proving {\it (ii)}.
\end{proof}

If the choice of $\omega$ is changed to $\tilde\omega=\omega^2$, then $\ell$ does not change:
$\l$ and $\mu$ are exchanged, and so the new $\z_0$ value is $\widetilde \z_0=1/\z_0$ 
and $\s^2(\widetilde\z)/\widetilde\zeta = \s^2(\z^{-1})/\z^{-1}=\omega^{-\ell}=\widetilde \omega^\ell$. 
This is to be expected; for example when $K=\F_q$, $\b^\ell$ reflects the form of
equation satisfied by $v\in Q_3^{-1}(\tau)$ (see (\ref{eq:type})),
and this is certainly independent of the choice of $\omega$.

\subsection{Explicit description of $\inv(\tau,\sigma)$ in characteristic 3. } \label{sec:G3char3}

As promised, we return to the case of char.~3.   
Then the only short $G_3$-orbit is $\{-1\}$, and the only irregular element is $Q_3(-1)=0$. Let $\tau \in \cj K^\x$ and write $\tau = Q_3(v)=(v+1)^3/(v(v-1))$.
Note that $v \not \in \F_3$ since $\tau \not \in \{\infty,0\}$.
Suppose $\s \in \Aut(\cj K/K)$ and $\s(\tau)=\tau$.
Then $\inv_{Q_3}(\tau,\sigma)=\beta^\ell$ is determined from
$\sigma(v)=\beta^\ell(v)$.  We wish to describe $\inv_{Q_3}(\tau,\sigma)$ 
purely in terms of $\tau$ and $\sigma$, without reference to $v$.  The approach of solving for $v$ in terms of $\tau$
no longer works in characteristic~3, so we must try something different.

\begin{proposition} \label{prop:G3char3} Let $K$ be a field of characteristic~3 and let $\tau \in \cj K$ be regular (so $\tau \ne 0$). 
Let $\sigma \in \Aut(\cj K/K)$
and assume $\s(\tau)=\tau$.  Let $\z\in \cj K$ satisfy $\z^3-\z=1/\tau$. Then there is $\ell \in \F_3$ such that $\sigma(\z)=\z+\ell$, 
and we have $\inv_{Q_3}(\tau,\sigma)=\beta^\ell$. 
\end{proposition}

\begin{proof}
Let $v \in Q^{-1}(\tau)$.  Since $\tau = (v^3+1)/(v^2-v)$, $v^3 - \tau v^2 + \tau v + 1 = 0$. Since $\tau\not\in\{ 0,\infty\}$, 
$v\not\in\F_3$.  The substitution $y=v+1$ eliminates the linear term:
$$y^3 - \tau y^2 + \tau = 0.$$
Let $\z = -1/y = -1/(v+1)$.  (Here note that $y \ne 0$ since $v \not\in\F_3$.)  Then
$$\z^3-\z = 1/\tau.$$
	Let $\ell=\sigma(\z)-\z$. Then $\ell\in\F_3$, because
	$$\ell^3=\s(\z^3)-\z^3=\s(\z+1/\tau)-(\z+1/\tau)=\s(\z)-\z=\ell.$$
	The other roots of $x^3-x=1/\tau$ are $\z+i$, $i \in \F_3$, and $\s(\z+i)-(\z+i)=\ell$ for all three roots.
	Since $\z=-1/(v+1)= \textmatrix0{-1}11(v)$, 
\begin{eqnarray*} \s(v) &=&  \s\left(\textmatrix {1}1{-1}0 (\z)\right) = \textmatrix {1}1{-1}0 \textmatrix 1\ell01 (\z) \\
&=& \textmatrix {1}1{-1}0 \textmatrix 1\ell01 \textmatrix 0{-1}11 (v) =\begin{pmatrix} \ell+1&\ell \\ -\ell & 1-\ell \end{pmatrix}(v) = \beta^\ell(v).
\end{eqnarray*}
Thus, $\inv_{Q_3}(\tau,\sigma)=\beta^\ell$, where $\ell = \sigma(\z)-\z$.
\end{proof}

\begin{corollary} \label{cor:G3char3Fq}
If $\tau \in \F_q^\x$, where $q=3^n$, then $\inv_{Q_3}(\tau,q)=\b^\ell$, where $\ell=\Tr_{\F_q/\F_3}(1/\tau)$.
\end{corollary}

\begin{proof}   Let $\z \in \cj\F_q$ satisfy $\z^3-\z=1/\tau$. By Proposition~\ref{prop:G3char3}, $\z^q=\z+\ell$ for some $\ell \in \F_3$, and $\inv_{Q_3}(\tau,q)=\b^\ell$.
We claim that $\ell=\Tr_{\F_q/\F_3}(1/\tau)$. Indeed, 
\begin{equation*} \Tr_{\F_q/\F_3}(1/\tau) = \sum_{i=0}^{n-1} (1/\tau)^{3^i} = \sum_{i=0}^{n-1}(\z^3-\z)^{3^i}  = \sum_{i=0}^{n-1} \left( \z^{3^{i+1}} - \z^{3^i}\right)
= \z^q-\z = \ell. \end{equation*}
\end{proof}

\subsection{A symbol with values in $\Z/3\Z$.} \label{sec:symbol}
Let $K=\F_q$, where $q$ is any prime power. Then
$\inv(\tau,q)\in G_3 = \set{I,\beta,\beta^{-1}}$ for regular $\tau \in \F_q$. $G_3$ is isomorphic to $\Z/3\Z$ abstractly, but making this explicit
requires selecting a preferred generator, which seemingly could equally well be $\beta$ or $\beta^{-1}$.  On the other hand, Corollary~\ref{cor:G3char3Fq} relates 
$\inv(\tau,q)$ with the absolute trace map in char.~3, which 
is genuinely a map to $\Z/3\Z$. Specifically, if $q=3^n$ and $\tau = Q_3(v) \in \F_q^\x$, then 
$$\text{$v^q = \beta^j(v) \iff \Tr_{\F_q/\F_3}(1/\tau)=j$, \ie, 
$\log_\beta(\inv(\tau,q)) = \Tr_{\F_q/\F_3}(1/\tau)$.}$$
The trace map 
determines the preferred generator $\beta\in G_3$, or equivalently a preferred isomorphism $\log_\beta : G_3 \to \Z/3$.
Then, for all characteristics, $\log_\beta (\inv(\tau))$ takes values in $\Z/3\Z$.  We summarize this in the theorem below.

\begin{theorem}  \label{thm:G3symbol} {\bf (A tripartite symbol)} 
Let $q$ be any prime power.  Let $\tau \in \F_q$, and assume $\tau^2-3\tau+9 \ne 0$.
Define a symbol $[\tau/q] \in  \Z/3\Z$ by
$$ \left[\tau/ q\right] = \Tr_{\F_q/\F_3}(1/\tau), \qquad \text{if $3|q$;} $$
$$  \left(\frac{\tau + 3\omega^2}{\tau+3\omega}\right)^{(q^2-1)/3} = \omega^{[ \tau/ q]},\qquad \text{if $3\nmid q$}, $$
where $\omega$ is any primitive cube root of unity in $\cj\F_q$ when $3\nmid q$.  
Let $\beta = \textmatrix1{-1}10 \in \PGL_2(\F_q)$.
Then $\beta$ has order~3, so $\beta^j$ is well defined for $j \in \Z/3\Z$.
If $v \in \cj\F_q$ such that $v^3-3v+1 = v(v-1)\tau$,  then
\begin{equation} v^q=\beta^{[\tau/ q]}(v). \label{eq:vqbeta} \end{equation}
Moreover, (\ref{eq:vqbeta}) determines $[\tau/ q]\in\Z/3\Z$ uniquely and could serve as an alternative definition for the symbol.  
\end{theorem}
\begin{proof} This is immediate from Theorem~\ref{thm:G3charneq3} and 
	Corollary~\ref{cor:G3char3Fq}.
\end{proof}

Theorem~\ref{thm:G3symbol} weaves together fields of different characteristic in a remarkable way.  
The group $G_3$ is defined over the integers, so $G_3\subset \PGL_2(K)$ for any field $K$. Taking $K=\F_q$,
this is consistent with the fact that 3 divides $|\PGL_2(\F_q)|=q(q-1)(q+1)$; in particular, 3 could divide $q$, $q-1$, or $q+1$.
Generally speaking, these three cases behave very differently, \eg,
when $3|q$, $\beta$ has a single eigenvalue and it is conjugate to $\textmatrix 1101$. If $3|q-1$ then $\beta$ has two rational eigenvalues and it is
conjugate to a diagonal matrix of order~3. If $3|q+1$, then $\beta$ has a pair of quadratic irrational eigenvalues. The equation~(\ref{eq:vqbeta})
applies to all these cases simultaneously.

We point out the similarity between Corollary~\ref{cor:G3char3Fq} 
(when $q=3^n$) and Proposition~\ref{prop:BDD}(2) (when $q=2^n$).
(See the remark following Proposition~\ref{prop:BDD}.)
In Corollary~\ref{cor:G3char3Fq}, the order-3 element $\beta$ is conjugate to $\textmatrix 1101$ in characteristic~3, and 
$$\inv(\tau,q) = \beta^{\Tr_{\F_q/\F_3}(1/\tau)}.$$
In Proposition~\ref{prop:BDD}(2) with $G=\{\textmatrix 1001,\textmatrix 0110\}$, the order-2 element $\gamma=\textmatrix 0110$ is conjugate to 
$\textmatrix 1101$ in characteristic~2, and 
$$\inv(\tau,q) = \gamma^{\Tr_{\F_q/\F_2}(1/\tau)}.$$
In both cases, the appearance of the trace map can be explained by Lemma~\ref{lem:conjugateInv}
and Proposition~\ref{prop:dim1codim1}{\it (i)}.

\section{The dihedral group of order~6} \label{sec:six}

Consider the subgroup $G_6 \subset \PGL_2(K)$ that is generated by the transformations 
$\beta(x)=1-1/x$ and $\rho(x)=1/x$.  This group is dihedral  of order~6 since $\beta^3=1$, $\rho^2=1$,
and $\rho\beta^i\rho = \beta^{-i}$.

The orbit of $v\in\cj K \cup \{\infty\}$ is
$$\calO_v = \set{v,1-1/v,1/(1-v),1/v,1-v,v/(v-1)}.$$
Note that $\calO_\infty= \{\infty,1,0\}$ is short.
To find the other short orbits, we find all $v \in \cj K \setminus \{0,1\}$ such that $v=\gamma(v)$ for some $1\ne \gamma \in G$.
If char$(K) \not \in \{2,3\}$, then
\begin{itemize}
\item $v = 1-1/v$ iff $v \in \{-\omega,-\omega^2\}$, where $\omega^2+\omega+1=0$;
\item $v = 1/(1-v)$ iff $v \in \{-\omega,-\omega^2\}$;
\item $v = 1/v$ iff $v  = -1$.  (Here $v=1$ is excluded since it belongs to $\calO_\infty$.)
\item $v=v/(v-1)$ iff $v = 2$. (Here $v=0$ is excluded since it belongs to $\calO_\infty$.)
\item $v = 1-v$ iff $v = 1/2$.
\end{itemize}
The short orbits are $\{-\omega,-\omega^2\}$, $\{-1,2,1/2\}$, and $\calO_\infty$ in that case.
If char$(K)=2$ then the latter three equations have no solutions in $\cj K  \setminus \F_2$, and the only short orbits in $\cj K \cup \{\infty\}$
are $\{-\omega,-\omega^2\} = \F_4\setminus \F_2$ and $\calO_\infty = \F_2 \cup \{\infty\}$. Finally, if char$(K)=3$ then
the short orbits are $\{-1\}$ and $\calO_\infty$.

\begin{lemma} \label{lem:G6Q} The function
\begin{equation} \label{def:Q6}
Q_6(x) = \frac{(x^3-3x+1)(x^3-3x^2+1)}{x^2(x-1)^2}. 
\end{equation}
is a quotient map for $G_6$ over any field $K$.
	The set of irregular elements of $\cj K\cup \{\infty\}$ is $S\cup\{\infty\}$, where
\begin{equation} \label{def:S}
	S = \begin{cases} \{-9,-9/4\} & \text{if char$(K) \not \in \{2,3\}$} \\ \{1\} & \text{if char$(K)=2$} \\ \{0\} & \text{if char$(K)=3$.} \end{cases} 
\end{equation}
\end{lemma}
\begin{proof} Note that $Q_6(x)=-Q_3(x)Q_3(1/x)$, where $Q_3(x)=(x^3-3x+1)/(x(x-1))$ is the quotient map for $G_3$
given in Lemma~\ref{lem:G3Q}. 
It is clear that $Q_6(x)=Q_6(1/x)$. Also, 
$Q_6(\beta(x)) = -Q_3(\beta(x)) Q_3(\rho(\beta(x))) = -Q_3(x) Q_3(\beta^2\rho(x)) = Q_6(x)$.  This proves $G$-invariance.  The degree of
the numerator is $6=|G_6|$ and the degree of the denominator is $<|G_6|$.  Thus, $Q_6$ satisfies all required properties to be a quotient map for $G_6$.

	The irregular elements are $Q_6(v)$, where $v$ is in a short $G_6$-orbit. The short $G_6$-orbits were computed in the paragraph preceding the statement of the lemma.
	Computing their images under $Q_6$ demonstrates that the irregular elements are $S\cup\{\infty\}$.
\end{proof}

Amusingly, Artin \cite[\S II.G]{Artin} considers the 
particular example of finding the fixed field in $K(x)$ to the set of automorphisms $f(x) \mapsto f(\gamma(x))$ for $\gamma \in G_6$.
He finds that the fixed field is $K(I)$ where $I(x)=(x^2-x+1)^3/\left(x^2(x-1)^2\right)$.
Note that $I(x)=Q_6(x) + 9$. In the notation of Theorem~\ref{thm:computeQ}, $I(x) = f_\calO(x)/g(x)$, where $\calO=\{-\omega,-\omega^2\}$.

Recall that $\inv_{Q_6}(\tau,\sigma)$ is a certain conjugacy class in $G_6$. The conjugacy classes are
$$\calC_\beta = \{\beta,\beta^2\}, \quad \calC_\rho = \set{\rho,\rho\beta,\rho\beta^2},\quad \calC_1=\{1\}.$$
The next proposition computes $\inv_{Q_6}(\tau,\s)$.

\begin{proposition} Let $\tau \in \cj K$ be regular with respect to $Q_6$
and let $\sigma \in \Aut(\cj K/K)$ such that $\s(\tau)=\tau$.  
Let $z \in \cj K$ satisfy $z^2-3z=\tau$. Then \\
{\it (i)}\ $z$ is regular with respect to $Q_3$. \\
{\it (ii)}\ If $\sigma(z)=z$ then $\gamma = \inv_{Q_3}(z,\sigma) \in G_3$ is defined, and  $\inv_{Q_6}(\tau,\sigma)=\calC_\g$. \\
{\it (iii)}\ If $\sigma(z)\ne z$, then $\inv_{Q_6}(\tau,\sigma) = \calC_\rho$. \\
{\it (iv)}\  If {\rm char}$(K)\ne 2$, then $\sigma(z)=z$ iff $\s\!\left(\sqrt{9+4\tau\,}\right)=\sqrt{9+4\tau\,}$. If $K=\F_{2^n}$ and $\sigma(x)=x^{2^n}$, then
$\sigma(z)=z$ iff $\Tr_{\F_{2^n}/\F_2}(\tau)=0$.
\end{proposition}

\begin{proof} By Lemma~\ref{lem:reciprocal}, $Q_6(x)=-Q_3(x)Q_3(1/x) = -Q_3(x)(3-Q_3(x))$. Thus, $Q_6(x)=h\circ Q_3(x)$
	and $\tau = h(z)$, where $h(x)=x^2-3x$. The statements {\it (i)} and {\it (ii)} are true by Lemma~\ref{lem:Hinv}, and {\it (iv)} is well known. 
	
	It remains only to prove~{\it (iii)}. Let $v\in Q_3^{-1}(z)$.
	Then $Q_6(v)=h(Q_3(v))=h(z)=\tau$, so that by Proposition~\ref{prop:Qorbit}, $Q_6^{-1}(\tau)=\{\g(v) : \g \in G\}$.
	Since $Q_6(\s(v))=\s(\tau)=\tau$, $\s(v)=\g(v)$ for some $\g \in G_6$, and by (\ref{eq:invariantK}), $\inv_{Q_6}(\tau,\s)=\calC_\g$.
	If $\g \in H$, then
	$$\sigma(z)=\sigma(Q_3(v))=Q_3(\sigma(v))=Q_3(\g(v)) = Q_3\circ\g(v)=Q_3(v)=z,$$
	contrary to the hypothesis that $\s(z)\ne z$. Then $\g\in G \setminus H
	=\calC_\rho$, so $\calC_\rho=\calC_\g=\inv_{Q_6}(\tau,\sigma)$.
\end{proof}

\part{Subgroups of $\PGL_2(\F_q)$} \label{part:q}

This part considers many different subgroups of $\PGL_2(\F_q)$, including Borel subgroups, unipotent subgroups, cyclic subgroups, 
$\PGL_2(\F_q)$, and $\PSL_2(\F_q)$. As indicated in the introduction (Examples~\ref{example:1b01} and~\ref{example:pgl2}), explicit computation of Artin invariants
reveals arithmetic information about additive polynomials and conjugacy classes of $\PGL_2(\F_q)$.

\section{Borel subgroup of $\PGL_2(\F_q)$}

The Borel subgroup $B_q\subset \PGL_2(\F_q)$ is defined as 
$$B_q=\set{\g \in \PGL_2(\F_q) : \g(\infty)=\infty} = \set{\textmatrix ab01  \in \PGL_2(\F_q) : a \in \F_q^\x,\ b \in \F_q}.$$
The cardinality is $q(q-1)$, and the short orbits are $\{\infty\}$ and $\F_q$.  The orbit $\F_q$ has multiplicity $|B_q|/q=q-1$. By
Theorem~\ref{thm:computeQ}, a quotient map is given by
$$Q(x) = (x^q-x)^{q-1}.$$
The irregular elements are the images of the short orbits under $Q$, namely 0 and $\infty$. 

The conjugacy classes of $B_q$ are $\calC_{\textmatrix a001}$ for $a \in \F_q^\x$ and $\calC_{\textmatrix 1101}$.
No element $\tau \in \F_q^\x$ has $\inv(\tau) = \calC_\textmatrix 1001$, because $v^q = v$ implies that $v$ belongs to a short $B_q$-orbit.

\begin{proposition} \label{prop:Borel} If $\tau \in \F_q^\x$ then 
$$ \inv(\tau) = \begin{cases} \calC_{\textmatrix 1101} & \text{if $\tau = 1$} \\
\calC_{\textmatrix \tau 001} & \text{otherwise.} 
\end{cases}
$$
\end{proposition}
\begin{proof} If $v \in \cj\F_q\setminus \F_q$ and $v^q = \textmatrix 1101 (v) = v+1$ then $Q(v) = (v^q-v)^{q-1} = 1$. Thus, $\inv(1)=\calC_{\textmatrix 1101}$.
	If $\tau\ne1$  and $v^q= \textmatrix \tau 001(v) = \tau  v$ then
	$Q(v) = (v^q-v)^{q-1}=v^{q-1}(v^{q-1}-1)^{q-1}=\tau (\tau -1)^{q-1}=\tau $. Thus, $\inv(\tau )=\calC_{\textmatrix \tau 001}$ when $\tau\ne1$.
\end{proof}

Now we generalize. Suppose that $q=P^e$, where $P$ is a prime power, and let 
$$H = B_P  = \set{ \textmatrix ab01 : a \in \F_P^\x,\ b \in \F_P}.$$
Then $Q_H(x) = (x^P-x)^{P-1}$ is a quotient map, $\{\infty\}$ and $\F_P$ are the only short $B_P$-orbits, and 0 and $\infty$ are irregular. 
The conjugacy classes of $H$ are $\calC_{\textmatrix a001,H}$ with $a \in \F_P^\x$ and $\calC_{\textmatrix 1101,H}$, where $\calC_{\g,H}$ denotes $\{\a\g\a^{-1} : \a\in H\}$.

Let N and Tr denote the polynomials in $\F_q[x]$:
$$N(x) = \prod_{i=0}^{e-1} x^{P^i} = x^{1+P+P^2+\cdots + P^{e-1}} = x^{(q-1)/(P-1)},\qquad
\Tr(x) = \sum_{i=0}^{e-1} x^{P^i}.$$
If $\tau \in \F_q$ then $N(\tau)=N_{\F_q/\F_P}(\tau)$ and $\Tr(\tau)=\Tr_{\F_q/\F_P}(\tau)$.
For each $a\in\F_P$, there are exactly $q/P$ elements in $\F_q$ with trace~$a$,
and $q/P$ is the degree of Tr$(x)$, thus 
$$\text{If $\tau \in \cj\F_q$ and Tr$(\tau)=a\in\F_P$ then $\tau \in \F_q$ and $\Tr_{\F_q/\F_P}(\tau)=a$.}$$
Likewise, if $\tau\in\cj\F_q$ and $N(\tau)\in\F_P^\x$ then $\tau\in\F_q^\x$, because $\tau^{q-1}=N(\tau)^{P-1}=1$.
Note also that if $\tau=s^{P-1}$ and $N(\tau)=1$ then $s\in\F_q^\x$, because $s^{q-1}=N(s^{P-1}) = N(\tau)$.

\begin{proposition} \label{prop:BorelGen} 
With respect to $H$ and $Q_H$ given above, for $\tau \in \F_q^\x$: \\
	{\it (i)}\ If N$(\tau)\ne 1$ then $ \inv_{Q_H}(\tau,q) = \ \calC_{\textmatrix {N(\tau)} 001,H}$. \\
{\it (ii)}\ If N$(\tau) = 1$ then we may write $\tau=s^{P-1}$ with $s \in \F_q^\x$, and 
	$$ \inv_{Q_H}(\tau,q) = \begin{cases} \calC_{\textmatrix 1001,H} & \text{if $\Tr(s)=0$} \\
\calC_{\textmatrix 1 101,H} & \text{if $\Tr(s)\ne 0$.}
\end{cases}
$$
\end{proposition}
\begin{proof} We will apply Lemma~\ref{lem:Hinv}, taking $G=B_q$ and $H=B_P$.  First, we compute a function $h$ such that $Q_G=h\circ Q_H$:
	\begin{eqnarray*}  Q_G &=& (x^q-x)^{q-1} = \left( \sum_{i=0}^{e-1} (x^P-x)^{P^i}\right)^{q-1} = (x^P-x)^{q-1}\left( \sum_{i=0}^{e-1}(x^P-x)^{P^i-1} \right)^{q-1} 
		\\ &=& N(Q_H)f(Q_H)^{q-1}, \quad
		\text{where $f(x)=\sum_{i=0}^{e-1} x^{(P^i-1)/(P-1)}$.}
	\end{eqnarray*}
	Thus, $h(x)=N(x)f(x)^{q-1}$.
	Let $\tau \in \F_q^\x$ and $v \in Q_H^{-1}(\tau)$, so $\tau = (v^P-v)^{P-1}$. Then $h(\tau)=h\circ Q_H(v)=Q_G(v)=(v^q-v)^{q-1}$. Since 
	$f(\tau)\in\F_q$, $h(\tau)=N(\tau)f(\tau)^{q-1} \in \{ N(\tau),0\}$. On the other hand, $h(\tau)=(v^q-v)^{q-1}$ vanishes if and only if $v\in \F_q$. Thus,
	$$h(\tau) = \begin{cases} N(\tau) & \text{if $v \not \in \F_q$} \\ 0 & \text{if $v\in \F_q$.} \end{cases} $$

		Since $(v^P-v)^{P-1}=\tau\ne 0$, $v \not \in \F_P$. 
		Let $s\in \cj\F_q$ such that $\tau=s^{P-1}$. Then 
		$(s/(v^P-v))^{P-1}=1$, so $s=c(v^P-v)$ with $c\in\F_P^\x$, and $\Tr(s)=c\Tr(v^P-v)=c(v^q-v)$. In particular, $\Tr(s)=0$ iff $v\in \F_q\setminus\F_P$. Thus, if $\Tr(s)=0$ then
	        $N(\tau)=(v^P-v)^{q-1}=1$,
		and the formulas $\tau=Q_H(v)$, $v^q=v$ imply $\inv_{Q_H}(\tau)=\calC_{\textmatrix1001}$.

	Now suppose $\Tr(s)\ne 0$, so $v\not\in\F_q$.
	Then $h(\tau)=N(\tau)\ne0$, so it is regular with respect to $Q_G$. By Lemma~\ref{lem:Hinv},
	there is $\d \in H$ such that $\inv_{Q_H}(\tau,q)=\calC_{\d,H}$ and $\inv_{Q_G}(N(\tau),q)=\calC_{\d,G}$. On the other hand, Proposition~\ref{prop:Borel} implies
	$\inv_{Q_G}(N(\tau),q) = \calC_{\g,G}$, where $\g=\textmatrix{N(\tau)}001$ if $N(\tau)\ne 1$ and $\g=\textmatrix 1101$ if $N(\tau)=1$.  In both cases,
	$\calC_{\g,G}=\left\{\textmatrix {N(\tau)} b01 : b \in \F_q\right\}$ and $\calC_{\g,G}\cap H=\calC_{\g,H}$. Since $\d \in \calC_{\d,G}\cap H = \calC_{\g,G}\cap H=\calC_{\g,H}$, 
	it follows that $\inv_{Q_H}(\tau,q)=\calC_{\d,H}=\calC_{\g,H}$.   
\end{proof}

\section{Unipotent subgroups of $\PGL_2(\F_q)$} \label{sec:unipotent}

Let $q = p^n$ where $p$ is prime. 
A {\it unipotent subgroup} of $\PGL_2(\F_q)$ is a group $G_W=\set{\textmatrix 1 w 0 1 : w \in W}$, where $W$ is an $\F_p$-vector subspace of $\F_q$.
Note that $\{\infty\}$ is the only short orbit.  All other orbits are cosets of $W$ in $\cj\F_q$, and each has cardinality $|G_W|=p^{\dim(W)}$.
Every element of $\cj\F_q$ is regular.
Thus, for each $\tau \in \F_q$, there is a unique $\g = \textmatrix 1j01 \in G_W$ such that
$v^q=\gamma(v)$ for all $v \in Q_W^{-1}(\tau)$. Otherwise put, if $Q_W(v)=\tau \in \F_q$, then $j := v^q-v \in W$, and $j$ depends only on $\tau$.

A simple example is $W = \F_p\subset \F_q$, and $G_W\subset \PGL_2(\F_q)$ is the subgroup of order $p$ generated by $\textmatrix 1101$. Then $Q(x)=x^p-x$
is a quotient map\footnote{The polynomial $x^p-x-\tau$ is called an Artin-Schreier polynomial. 
If $\Tr_{\F_q/\F_p}(\tau)\ne 0$ then its splitting field has degree~$p$
and Galois group $\Z/p\Z$. See \cite[VI, \S6, Th.~6.4]{Lang}.}.
Let $\tau = Q(v) \in \F_q$. Then
$v^q= \textmatrix 1j01 v = v+j$ for some $j \in \F_p$, and $\inv(\tau)=\textmatrix 1j01$.
To relate $\inv(\tau)$ to a quantity that is directly computable from $\tau$, we note that
\begin{eqnarray*}
\Tr_{\F_q/\F_p}(\tau) &=& \tau + \tau^p + \tau^{p^2} + \cdots + \tau^{q/p} \\
&=& (v^p-v)+(v^{p^2}-v^p)+(v^{p^3}-v^{p^2})+ \cdots + (v^q-v^{q/p})
 = v^q-v=j.
\end{eqnarray*}
Thus, the invariant coming from this group is essentially the absolute trace of $\tau$.

Now let $W$ be an arbitrary $\F_p$-vector subspace of $\F_q$. Then
$$Q_W(x) = \prod_{w \in W} (x-w)$$
is easily seen to be $G_W$-invariant, so it is a quotient map for $G_W$. 

\begin{lemma} $Q_W(x)$ is an additive polynomial, \ie, $Q_W(x+y)=Q_W(x)+Q_W(y)$.  \end{lemma}
\begin{proof} See Goss \cite{Goss}, Theorem 1.2.1. Alternatively, observe that $Q_W(x)$ and $Q_W(x+y)$ are both quotient maps for $G_W$ 
over the field $\F_q(y)$, as both are invariant under $G_W$ and of the right form.  By Proposition~\ref{prop:QExists},
there are $a,b\in \F_q(y)$ with $a\ne0$ such that $Q_W(x+y)=aQ_W(x)+b$.  Since $Q_W(x+y)$ and $Q_W(x)$ are both monic as polynomials in~$x$,
$a$ must be~1.  The value for $b$ can be found by setting $x=0$: $Q_W(y)=Q_W(0)+b=b$. Hence, $Q_W(x+y)=Q_W(x)+Q_W(y)$.
\end{proof}

Suppose now that $P$ is a power of $p$ and $\F_p\subset \F_P \subset \F_q$, and that $W$ is an $\F_P$-vector subspace of $\F_q$. In that case, $Q_W$ is $\F_P$-additive, meaning that it is additive and it satisfies the additional property that $Q_W(ax)=aQ_W(x)$ for all $a\in \F_P$.

\begin{lemma} A monic polynomial $L(x)\in\F_q[x]$ is $\F_P$-additive if and only if it has the form $x^{P^d}  + \sum_{i=0}^{d-1} a_i x^{P^i}$,
where $a_i \in \F_q$.
\end{lemma}
\begin{proof} See Goss \cite{Goss}, Proposition 1.1.5. \end{proof}

\begin{proposition} \label{prop:YW} Let $q=P^e$.
Let $W\subset\F_q$ be a $d$-dimensional $\F_P$-vector subspace and $Y = Q_W(\F_q)$. Then $Y$ is an $(e-d)$-dimensional
$\F_P$-vector subspace of $\F_q$, and $Q_Y \circ Q_W(x) = x^q-x$.
\end{proposition}

\begin{proof} Since $Q_W$ is $\F_P$-additive, it may be viewed as an $\F_P$-linear map from $\F_q$ to $\F_q$. Its image $Y$ is then an $\F_P$-vector subspace of $\F_q$.
Note that $W$ is the kernel. Because of the exact sequence $0 \to W \to \F_q \to Y \to 0$, $\dim_{\F_P}(W) + \dim_{\F_P}(Y) = \dim_{\F_P}(\F_q)=e$. 
Let $Z$ be a complementary subspace to $W$, that is, $\dim_{\F_P}(Z)=e-d$ and $Z+W=\F_q$. Then $Q_W$ maps $Z$ isomorphically onto $Y$, and
\begin{eqnarray*}
x^q-x &=& \prod_{a \in \F_q} (x-a) = \prod_{w \in W, z \in Z}(x-w-z) \\
&=& \prod_{z\in Z}Q_W(x-z) =  \prod_{z \in Z} (Q_W(x)-Q_W(z)) \\
&=& \prod_{y\in Y} (Q_W(x)-y) = Q_Y\circ Q_W(x).
\end{eqnarray*}
\end{proof}

\begin{proposition} \label{prop:YWinv}
Let $W$ and $Y$ be as in Proposition~\ref{prop:YW}. 
If $\tau\in \F_q$ then $Q_W^{-1}(\tau)$ is a $G_W$-orbit (\ie, a coset $v+W\subset \cj\F_q$), and there is a unique $\gamma = \textmatrix 1w01 \in G_W$
such that $v^q = \gamma(v) = v + w$ for all $v \in Q_W^{-1}(\tau)$.  Moreover, $w = Q_Y(\tau)\in W$.
\end{proposition}

\begin{proof} Writing $Q_W(v)= \tau$, we have $0=\tau^q-\tau=Q_W(v^q)-Q_W(v)=Q_W(v^q-v)$, so $v^q-v \in W$. Set $w=v^q-v$.
By Proposition~\ref{prop:YW}, $w=v^q-v=Q_Y(Q_W(v))=Q_Y(\tau)$.
\end{proof}

In the above proposition, $\gamma = \inv_{Q_W}(\tau)$ is expressed directly in terms of $\tau$:
$$\inv_{Q_W}(\tau) = \textmatrix 1 {Q_Y(\tau)} 0 1 \in G_W.$$
This is surprising, as it is not even obvious that $Q_Y(\tau)$ belongs to $W$.
The following corollary may be of independent interest.

\begin{corollary} \label{cor:duality} {\bf (Reciprocity)} If $q=P^e$, $W$ is a $d$-dimensional $\F_P$-vector subspace of $\F_q$, and $Y=Q_W(\F_q)$, then there are short exact
	sequences of $\F_P$-vector spaces:
$$0 \to W \stackrel{inc.}{\longrightarrow} \F_q \stackrel{Q_W}{\longrightarrow} Y \to 0  $$
and 
$$0 \to Y \stackrel{inc.}{\longrightarrow} \F_q \stackrel{Q_Y}{\longrightarrow} W \to 0.  $$ 
Moreover, $Q_Y \circ Q_W(x) = Q_W \circ Q_Y(x) = x^q-x$.
\end{corollary}

\begin{proof} $Q_W$ is an $\F_P$-linear map from $\F_q$ to $\F_q$, has kernel~$W$, and has image~$Y$, thus the first short exact sequence holds, and $Y$ is an $\F_P$-subspace of $\F_q$ of dimension
	$e-d$. Likewise, there is a short exact sequence of $\F_P$-vector spaces
$$0 \to Y \stackrel{inc.}{\longrightarrow} \F_q \stackrel{Q_Y}{\longrightarrow} V \to 0,  $$ 
	where $\dim_{\F_P}(V)=e-(e-d)=d$. 
By Proposition~\ref{prop:YWinv}, if $\tau \in \F_q$ then
$Q_Y(\tau) \in W$, thus 
$V=Q_Y(\F_q) \subset W$. Since $V$ and $W$ have the same dimension, they are equal. Proposition~\ref{prop:YW}
shows $Q_Y \circ Q_W(x) = x^q-x$.  Applying Proposition~\ref{prop:YW} again, but with the roles of $Y$ and $W$ reversed, 
and using that $W = Q_Y(\F_q)$ (which we have already proved), we deduce $Q_W \circ Q_Y = x^q-x$ also.
\end{proof}

\begin{proposition} \label{prop:linPoly} Let $q=P^e$ and let $L(x) = x^{P^d} + \sum_{i=0}^{d-1} a_i x^{P^i}$ be an $\F_P$-additive polynomial, 
	where $a_i \in \F_q$ and $a_0 \ne 0$. Then all the roots of $L$
	are in $\F_q$ if and only if there is an additive polynomial $M(x) = x^{P^{e-d}} + \sum b_i x^{P^i}\in\F_q[x]$ with $M \circ L(x) = x^q - x$.
	In that case, it is also true that $L \circ M(x) = x^q - x$, all the roots of $M$ are in $\F_q$, and $M=Q_Y$ where $Y=L(\F_q)$.
\end{proposition}

\begin{proof} The roots of $L$ in $\cj\F_q$ comprise a $d$-dimensional $\F_P$-vector subspace $W\subset \cj\F_q$. If $W\subset \F_q$, then $L=Q_W$ and the result follows
from Corollary~\ref{cor:duality}. Conversely, if there is an $\F_P$-additive polynomial $M(x)$ satisfying $M\circ L(x) = x^q-x$, then for any root $w\in\cj\F_q$ of $L$ we have
$0 = M\circ L(w) = w^q-w$, showing that $w \in \F_q$.  Thus, $L=Q_W$ where $W\subset \F_q$. By Corollary~\ref{cor:duality}, $Q_Y \circ L=x^q-x$, where $Y=Q_W(\F_q)$.
Then $M \circ L(x) = x^q-x = Q_Y \circ L(x)$. Set $z=L(x)$, which is transcendental. Since $M(z)=Q_Y(z)$, $M=Q_Y$.
By Corollary~\ref{cor:duality},  $L \circ M = Q_W \circ Q_Y = x^q-x$.
\end{proof}

\begin{proposition} \label{prop:dim1codim1} Suppose $q=P^e$, where $P$ is a prime power. \\
{\it (i)}\ Let $W$ be a one-dimensional $\F_P$-subspace of $\F_q$, so $W=c\,\F_P$ where $c \in \F_q^\x$. A quotient map is $Q_W(x)=x^P-c^{P-1}x$.
If $\tau \in \F_q$ then $\inv_{Q_W}(\tau) = \textmatrix 1w01$, where $w=c \Tr_{\F_q/\F_P}(\tau/c^P)$. That is, $Q_W(v)=\tau\in\F_q$ implies $v^q-v = c \Tr_{\F_q/\F_P}(\tau/c^P)$.\\
{\it (ii)}\ If $Y$ is an $(e-1)$-dimensional $\F_P$-vector subspace of $\F_q$, then for $\tau \in \F_q$, 
$$\inv_{Q_Y}(\tau) = \begin{pmatrix} 1 & {\tau^P-\tau/a_0} \\ 0& 1 \end{pmatrix},\qquad \text{where $a_0 = \prod_{0\ne y \in Y} y$.} $$
That is, $Q_Y(v)=\tau\in\F_q$ implies $v^q-v=\tau^P-\tau/a_0$.
	Also, $a_0=1/c^{P-1}$ for some $c \in \F_q^\x$, and $Q_Y(\F_q)=c \F_p$.
\end{proposition}

\begin{proof} Corollary~\ref{cor:duality} implies $W\longleftrightarrow Y$ gives a bijection between the set of all 1-dimensional $\F_P$-vector spaces $W=c\,\F_P$ and the set of all
	$(e-1)$-dimensional vector spaces $Y$. In this correspondence, $Y=Q_W(\F_q)$, $W=Q_Y(\F_q)$, and $Q_W\circ Q_Y=Q_Y\circ Q_W=x^q-x$. 
	Assume $W$ and $Y$ are so related.

First, $Q_W(x)= \prod_{j\in\F_P} (x-jc) = c^P \prod_{j \in \F_P} ((x/c)-j) = c^P ( (x/c)^P - (x/c) ) = x^P - c^{P-1}x$.
	If Tr is the polynomial $\sum_{i=0}^{e-1} x^{P^i}$ and $L(x)=c\Tr(x/c^P)$ then 
	$$L\circ Q_W(x)=c\Tr(x^P/c^P-x/c)=c( (x/c)^q-x/c)=x^q-x,$$
	therefore $L=Q_Y$. Since $\inv_{Q_W}(\tau)=\textmatrix 1 {Q_Y(\tau)} 01 $, {\it (i)} follows.

	Next, $\inv_{Q_Y}(\tau)=\textmatrix 1{Q_W(\tau)}01$,  and $Q_W(\tau)=\tau^P-c^{P-1}\tau$. Let $a_0 = \prod_{0\ne y\in Y} y$. Since $Q_Y(x) = \prod_{y\in Y}(x-y) =
	\prod_{y\in Y}(x+y)$, $a_0$ is
	the coefficient of $x$ in $Q_Y(x)$. Since $Q_Y=L$, this coefficient is $1/c^{P-1}$. Thus, $Q_W(\tau)=\tau^P-\tau/a_0$. This proves {\it (ii)}.
\end{proof}

\section{Cyclic subgroups of $\PGL_2(\F_q)$} \label{sec:cyclic}
In this section, we find quotient maps and Artin invariants of cyclic subgroups $G\subset \PGL_2(\F_q)$, and prove that the Artin invariant is equidistributed -- every $\g\in G$
has the same number of regular elements $\tau \in \F_q\cup\{\infty\}$ such that $\inv(\tau)=\g$.  We also study 
the equation $v^q=\g(v)$ when $\g\in \PGL_2(\F_q)$ and $v \in \cj\F_q\cup\{\infty\}$.

\subsection{Dickson's analysis.}
Cyclic subgroups of $\PSL_2(\F_q)$ were analyzed
by Dickson \cite{Dickson}. We modify his analysis to obtain the cyclic subgroups of $\PGL_2(\F_q)$. 

It is useful to introduce the following matrices.
For $\lambda\in \F_{q^2} \setminus \F_q$, define $C_\lambda \in \GL_2(\F_{q^2})$ by
\begin{equation} C_\lambda  = \begin{pmatrix} \lambda & -\lambda^q \\ 1 & -1 \end{pmatrix},\qquad {\rm so}\quad
C_\l^{-1} = (\l-\l^q)^{-1} \begin{pmatrix} 1&-\l^q \\ 1 & -\l \end{pmatrix}.
\label{eq:Clambda}
\end{equation}
Then
\begin{equation} C_{\l^q} = \begin{pmatrix} \l^q & -\l \\ 1&-1 \end{pmatrix} 
= C_\l \begin{pmatrix} 0&-1 \\-1&0 \end{pmatrix},
 \label{eq:Clambdaq} \end{equation}
so $C_{\l^q}(x) = C_\l(1/x)$. Note that $C_\l(\infty)=\l$, $C_\l(0)=\l^q$, $C_\l(1)=\infty$.

If $M \in \cj\F_q^\x\GL_2(\F_q)$ then its order as an element of $\PGL_2(\F_q)$ is
the least $\ell\ge1$ such that $M^\ell$ is a scalar matrix.

\begin{proposition} \label{prop:Dickson}
Suppose that $M \in \GL_2(\F_q)$ and that $M$ has order $\ell > 1$ as an element of $\PGL_2(\F_q)$. 
Then exactly one of the following holds: (a)\  $M$ has a unique fixed point $z \in \cj\F_q \cup \{\infty\}$, which is rational over $\F_q$; 
(b)\ $M$ has two distinct fixed points $z_1,z_2 \in \cj\F_q \cup \{\infty\}$, which are both rational; 
or (c)\ $M$ has two fixed points in $\cj \F_q \cup \{\infty\}$, which
are a conjugate pair $\{\l,\l^q\}$, with $\l \in \F_{q^2} \setminus \F_q$.

If (a) holds, then $\ell=p$, the prime that divides $q$.  If $z = \infty$ then
$kM = \textmatrix 1b01$ for some $b,k\in \F_q^\x$. If $z \in \F_q$ then $kM = \a \textmatrix 1 b01 \a^{-1}$ for some $b,k \in \F_q^\x$, 
where $\a=\textmatrix z011 \in \PGL_2(\F_q)$.

If (b) holds then $\ell$ divides $q-1$.
Let $\a \in \GL_2(\F_q)$ such that $\a(z_1) = \infty$ and $\a(z_2) = 0$, for example $\a=\textmatrix 1{-z_2}1{-z_1}$. 
Then $\a M \a^{-1} = \textmatrix a00d$, where $a,d \in \F_q^\x$
and $a/d$ has order $\ell$.

If (c) holds then $\ell$ divides $q+1$, and there is 
$\d\in \F_{q^2}\setminus \F_q$ such that $M = D_{\d,\l}$, where
\begin{equation} D_{\d,\l} = C_\l \textmatrix {\d^q} 0 0 {\delta} C_\l^{-1}. \label{def:Ddel} \end{equation}
Further, $\d^{q-1}$ is a primitive $\ell$th root of unity.
Conversely, if $\d,\l \in \F_{q^2} \setminus \F_q$ and $M=D_{\d,\l}$, then $M \in \GL_2(\F_q)$, the fixed points of $M$ in $\cj\F_q\cup\{\infty\}$ are
$\{\l,\l^q\}$, and the order of $M$ as an element of $\PGL_2(\F_q)$ equals $\circ(\d^{q-1})$, the multiplicative order of $\d^{q-1}$.
\end{proposition}

\begin{proof}  Write $M = \textmatrix abcd$. The fixed points of $M$ in $\cj\F_q$ are $z$ such that 
$az+b=z(cz+d)$. Note that $\infty$ is a fixed point if and only if $c=0$, and it is the unique fixed point if and only if $c=0$, $a=d$, and $b\ne 0$.
The quadratic either has a single repeated root in $\F_q \cup \{\infty\}$, 
two distinct roots in $\F_q \cup \{\infty\}$, or a pair of conjugate roots
$\l,\l^q$ where $\l \in \F_{q^2}\setminus \F_q$.  This gives rise to the mutually exclusive cases (a), (b), and (c). 

(a)\ Suppose the quadratic equation has a single repeated root $z\in\F_q\cup \{\infty\}$. Let $\a \in \GL_2(\F_q)$ such that $\a(z)=\infty$. 
Then $\a M\a^{-1}$ fixes $\infty$ only, so it
is a scalar multiple of $\textmatrix 1b01$, where $b \in \F_q$ and $b\ne 0$. The order of $M$ as an element of $\PGL_2(\F_q)$ is 
equal to $p$, the characteristic of $\F_q$.

(b)\ Suppose the quadratic equation has two distinct roots $z_1,z_2 \in \F_q \cup \{\infty\}$. Let $\a\in \GL_2(\F_q)$ such that $\a(z_1)=\infty$ and $\a(z_2)=0$.
Then $\a M\a^{-1}$ fixes 0 and $\infty$, so it has the form $\textmatrix a00d$.  

(c)\ 
Finally, suppose that the quadratic equation has no rational roots. Then $c\ne0$. The roots of the quadratic are a pair $\l,\l^q \in \F_{q^2}$.
$C_\l^{-1} M C_\l$ fixes $\infty$ and $0$, so it has the form $\textmatrix \a 0 0 \d$, where
$\a + \d = \Tr(M)$ and $\a\d = \det(M)$.  Thus, $M=C_\l \textmatrix \a 00\d C_\l^{-1}$.  Note that $(x-\a)(x-\d) = x^2 -\Tr(M) x + \det(M)$, so $\a$
and $\d$ are either rational, or they form a conjugate pair in $\F_{q^2}$.  If they are rational, then by applying the Frobenius to all coefficients we find:
$$M=C_{\l^q} \textmatrix \a 00 \d C_{\l^q}^{-1}  = C_\l \textmatrix 0{-1}{-1}0 \textmatrix \a 0 0 \d \textmatrix 0{-1}{-1}0 C_\l^{-1} = C_\l \textmatrix \d 0 0 \a C_\l^{-1},$$
which would imply that $\a = \d$, so that $M$ is a scalar matrix. However, we assumed that $M$ has order $\ell>1$ as an element of $\PGL_2(\F_q)$, so 
we obtain a contradiction.  We conclude that $\a,\d$ are a conjugate pair in $\F_{q^2}$, \ie, $\d \in \F_{q^2}\setminus \F_q$ and  $\a=\d^q$. Then
$$M = C_\l \textmatrix {\d^q} 00 \d  C_\l^{-1},\qquad\text{where $\d \in \F_q^{2} \setminus \F_q$.} $$
As an element of $\PGL_2(\cj\F_q)$, $M$ is equivalent to $\d^{-1}M = C_\l \textmatrix \z001 C_\l^{-1}$, where $\z = \delta^{q-1}$.
If $E$ denotes the matrix on the right, then it is clear that $E^i$ is scalar if and only if $\z^i=1$, therefore the order of $M$ in $\PGL_2(\F_q)$
is $\circ(\zeta)$.  Since the order of $\d$ divides $q^2-1$, $\ell=\circ(\z)$ divides $q+1$.

Now we prove the final statement of (c). Let $\l,\d \in \F_{q^2}\setminus \F_q$ and let $M=D_{\d,\l}$. Applying the Frobenius, we find that
$$ M^{(q)} = C_{\l^q} \textmatrix \d 00{\d^q} C_{\l^q}^{-1} = C_\l \textmatrix 0{-1}{-1}0 \textmatrix \d 00{\d^q} \textmatrix 0{-1}{-1}0 C_\l^{-1} = M.$$
Thus, $M$ is rational. 
\end{proof}  

If $\z^{q+1} = 1$, define
\begin{equation}  E_{\z,\l} = C_\l \textmatrix \z001 C_\l^{-1}. \label{eq:Ezeta} \end{equation}
This is not rational as a matrix, but it is rational as an element of $\PGL_2(\F_q)$, and in fact if $\delta^{q-1}=\z$ (so $\delta^{q^2-1}=1$, \ie, 
$\delta \in \F_{q^2}^\x$) then
$$E _{\z,\l} = \delta^{-1} D_{\d,\l} \in \cj\F_q^\x \GL_2(\F_q).$$
In particular, $E_{\z,\l}=D_{\d,\l}$ as elements of $\PGL_2(\F_q)$.

\subsection{Quotient map and Artin invariant of a cyclic group.}

Let $G$ be a cyclic subgroup of $\PGL_2(\F_q)$ of order $\ell>1$ generated by
$M\in\PGL_2(\F_q)$.
Because of Lemma~\ref{lem:conjugateInv}, to understand the Artin invariant of $G$, we may first conjugate by any $\a\in\PGL_2(\F_q)$. By Proposition~\ref{prop:Dickson},
there are three cases:

(a)\quad $M=\textmatrix 1b01$ where $b \in \F_q^\x$ and $\ell=p$. 
By Proposition~\ref{prop:dim1codim1}{\it (i)}, a quotient map is $Q_G(x)=x^p-b^{p-1}x$, every $\tau \in \F_q$ is regular, and $\inv(\tau)= \textmatrix 1w01$,
where $w=b\Tr_{\F_q/\F_p}(\tau/b^p)$. 

(b)\quad $M=\textmatrix a001$ and $\ell=\circ(a)$. 
This is the Kummer case (Example~\ref{example:Kummer}).  
A quotient map is $Q(x)=x^\ell$. The irregular elements are 0 and $\infty$,
and for $\tau \in \F_q^\x$, $\inv(\tau) = \textmatrix {\tau^{(q-1)/\ell}}001$.

(c)\quad $M=D_{\d,\l}=E_{\z_\ell,\l}$, where $\l,\mu\in\F_{q^2}\setminus\F_q$, $\z_\ell=\d^{q-1}$, $\ell=\circ(\z_\ell)$, and $\ell|(q+1)$.  Then the group generated by $M$ is
$$ G_\ell = \{\,E_{\z_\ell^i,\l} : 0\le i<\ell\} =  \set{E_{\z,\l} : \z^\ell = 1 }.$$
We will compute a quotient map and Artin invariant for this group.
We begin by finding the short orbits.

\begin{lemma} \label{lem:AppBfixed} If $\zeta^{q+1}=1$ and $\zeta\ne1$ then for $v \in \cj\F_q \cup \{\infty\}$,
$E_{\zeta,\lambda}(v)=v$ if and only if $v \in \{\lambda,\lambda^q\}$. \end{lemma}
\begin{proof} Since $E_{\zeta,\lambda}=C_\lambda \textmatrix{\zeta} 001 C_\lambda^{-1}$ with $C_\l = \textmatrix \l{-\l^q}1{-1}$, 
\begin{eqnarray*} E_{\zeta,\lambda}(v)=v &\iff& \zeta C_\lambda^{-1}(v)=C_\lambda^{-1}(v)
\iff C_\lambda^{-1}(v)\in \{0,\infty\} \\
&\iff& v \in \{\lambda,\lambda^q\}.
\end{eqnarray*}
\end{proof}

\begin{lemma} \label{lem:AppBshort} The short $G_\ell$-orbits in $\cj\F_q \cup \{\infty\}$ are $\{\lambda\}$ and $\{\lambda^q\}$.  
	$\calO_\infty$ consists of $\infty$ together with $\ell-1$ elements of $\F_q$.
\end{lemma}
\begin{proof}
First, $\lambda$ and $\lambda^q$ are each fixed by every element of $G_\ell$ by Lemma~\ref{lem:AppBfixed}, so they form singleton orbits.
The same lemma shows that $\lambda$ and $\lambda^q$ are the only elements of $\cj\F_q$ that are fixed by a nontrivial element of $G_\ell$.
By Lemma~\ref{lem:short} it follows that $\{\lambda\}$ and $\{\lambda^q\}$ are the only short orbits.

For the last statement, note that $E_{\z,\l}(\infty) = D_{\d,\l}(\infty) \in \F_q$ where $\d^{q-1}=\z$, and no element of $\F_q$ is in a short orbit.
Therefore $\calO_\infty \subset \F_q \cup \{\infty\}$ and $|\calO_\infty|=\ell$. 
\end{proof}

\begin{proposition} \label{prop:Qell}  Let $\ell$ divide $q+1$ and $G_\ell = \set{ E_{\z,\l} : \z^\ell = 1}$.
A quotient map for $G_\ell$ is given by  
\begin{equation}
	Q_\ell(x)=C_\l \circ [\ell] \circ C_\l^{-1}(x) = \frac{\l(x-\l^q)^\ell - \l^q(x-\l)^\ell}{(x-\l^q)^\ell - (x-\l)^\ell}, 
\label{eq:Qell}
\end{equation}
where $[\ell]$ denotes the $\ell$th power map: $[\ell](x) = x^\ell$.
\end{proposition}
\begin{proof} Let $Q_\ell=C_\l\circ [\ell]\circ C_\l^{-1}$. Since $[\ell]\circ C_\l^{-1}(x)=\left(\frac{x-\l^q}{x-\l}\right)^\ell$, (\ref{eq:Qell}) holds.  
We need to prove that $Q_\ell$ is $G$-invariant, has degree $\ell$, is $\F_q$-rational, and carries $\infty$ to $\infty$.

Note that $[\ell]\circ\textmatrix \z001(x) = (\z x)^\ell=x^\ell = [\ell](x)$ when $\z\in \mu_\ell$. Therefore,
$$
Q_\ell \circ E_{\z,\l} = (C_\l \circ [\ell] \circ C_\l^{-1}) \circ (C_\l \circ \textmatrix \z001 \circ C_\l^{-1}) =Q_\ell.$$
By Lemma~\ref{lem:degf}, $\deg(Q_\ell) = \deg([\ell])=\ell$.
	By (\ref{eq:Qell}), the numerator of $Q_\ell$ 
has degree~$\ell$ and the denominator has degree~$<\ell$. (Alternatively, the degree of the denominator is less than the degree of the numerator
iff $Q_\ell(\infty)=\infty$. We have $Q_\ell(\infty)=C_\l[\ell]C_\l^{-1}(\infty)=C_\l[\ell](1)=C_\l(1)=\infty$.)

To complete the proof that $Q_\ell$ is a quotient map, it remains only to prove rationality.
When Frobenius is applied to the coefficients of the rational function $Q_\ell$, $\l$ and $\l^q$ are exchanged. 
Both the numerator and denominator of (\ref{eq:Qell}) are negated, so $Q_\ell$ remains invariant. Alternatively, 
since $C_{\l^q}(x)=C_\l(1/x)=C_\l \circ [-1]$, the conjugate of $Q_\ell$ is
$$C_{\l^q} \circ [\ell] \circ C_{\l^q}^{-1} = (C_\l \circ [-1]) \circ [\ell] \circ (C_\l \circ [-1])^{-1} = C_\l \circ [-1\cdot \ell \cdot -1] \circ C_\l^{-1}=Q_\ell.$$
\end{proof}

\begin{lemma} \label{lem:irregQell}
	The only irregular elements for $Q_\ell$ are $\l$ and $\l^q$.  
	In particular, every element of $\F_q\cup\{\infty\}$ is regular with respect to $Q_\ell$. 
\end{lemma}
\begin{proof} Recall $\tau$ is irregular iff $Q_\ell^{-1}(\tau)$ is a short orbit for $G_\ell$, and the only short orbits are $\{\l\}$ and
$\{\l^q\}$. Thus, the only irregular elements are $Q_\ell(\l)$ and $Q_\ell(\l^q)$. 
By~(\ref{eq:Qell}), $Q_\ell(\l)=\l$ and $Q_\ell(\l^q)=\l^q$.
\end{proof}

\begin{lemma} \label{lem:uv} If $v = C_\l(u)$ then $v^q = C_\l(u^{-q})$.
\end{lemma}
\begin{proof} By (\ref{eq:Clambdaq}), $v^q=C_{\l^q} (u^q) = C_\l\textmatrix 0{-1}{-1}0 (u^q) = C_\l(u^{-q})$.  
\end{proof}

\begin{theorem} \label{thm:invQell} If $\tau \in \F_q$, then 
	$\inv_{Q_\ell}(\tau) = E_{\z,\l}$ where 
	$\z = \left(\frac{\tau-\l}{\tau-\l^q}\right)^{(q+1)/\ell}$. 
\end{theorem}

\begin{proof} 
We have $\inv(\tau) = E_{\z,\l}$ for some $\z\in\mu_\ell$, and we must show $\zeta = \left(\frac{\tau-\l}{\tau-\l^q}\right)^{(q+1)/\ell}$. 

	Let $v\in Q_\ell^{-1}(\tau)$ and $u = C_\l^{-1}(v)$, so that $v^q = E_{\z,\l}(v) = C_\l(\z u)$. Since $v^q = C_\l(u^{-q})$ by Lemma~\ref{lem:uv}, 
$C_\l(\z u) = C_\l(u^{-q})$, and so $\z = u^{-(q+1)}$.
Since $\tau = Q_\ell(v) = C_\l\left([\ell](u)\right) = C_\l(u^\ell)$, it follows that $u^\ell = C_\l^{-1}(\tau)$. Consequently,
$$\z = u^{-(q+1)} =(u^\ell)^{-(q+1)/\ell} = \left(C_\l^{-1}(\tau)\right)^{-(q+1)/\ell} = \left(\frac{\tau-\l}{\tau-\l^q}\right)^{(q+1)/\ell}.$$
\end{proof}

The next lemmas find a symmetry of $\inv_{Q_\ell}$.  

\begin{lemma} \label{lem:QellR} Let
$R = C_\l \textmatrix 0110 C_\l^{-1}$.
Then
\begin{equation} Q_\ell\circ R(x) = 
\frac{\l^q (x-\l^q)^\ell -\l(x-\l)^{\ell}}{(x-\l^q)^\ell - (x-\l)^{\ell}},
\label{eq:QellR} \end{equation}
and 
\begin{equation} Q_\ell(x) + Q_\ell\circ R(x) = \l + \l^q. \label{eq:QellPlusQellR} \end{equation}
\end{lemma}
\begin{proof}
$R=C_\l \circ \textmatrix 0110 \circ C_\l^{-1} = C_\l \circ [-1] \circ C_\l^{-1}$, therefore
\begin{eqnarray*}
Q_\ell\circ R(x) &=& C_\l\circ[-\ell]\circ C_\l^{-1}(x)= 
	C_\l\left(\left(\frac{x-\l}{x-\l^q}\right)^\ell\right) \\
&=& \frac{\l^q (x-\l^q)^\ell -\l(x-\l)^{\ell}}{(x-\l^q)^\ell - (x-\l)^{\ell}}.
\end{eqnarray*}
	A simple computation using (\ref{eq:Qell}) then shows $ Q_\ell(x) + Q_\ell\circ R(x) = \l+\l^q$.
\end{proof}

\begin{lemma} \label{lem:QellSymmetry}
Let $\tau \in \F_q$. If $\inv_{Q_\ell}(\tau) = \g $ then $\inv_{Q_\ell}(\l+\l^q-\tau)=\g^{-1}$.  \end{lemma}
\begin{proof} Write $\tau = Q_\ell(v)$. Then $Q_\ell\left(R(v)\right)=\l+\l^q-Q_\ell(v)=\l+\l^q-\tau$. If $\inv_{Q_\ell}(\tau)=\g$ then $v^q = \g(v)$,
so $(R(v))^q = R(v^q) = R \left(\g(v)\right) = \g^{-1}\circ R(v)$. Let $w=R(v)$. Since $Q_\ell(w) = \l+\l^q-\tau$ and $w^q=\g^{-1}(w)$, it 
follows that $\inv(\l+\l^q-\tau)= \g^{-1}$.
\end{proof}

\subsection{Equidistribution of the Artin invariant when $G$ is cyclic.}
We show that when $G$ is cyclic, the number of regular $\tau\in\F_q\cup\{\infty\}$ with $\inv(\tau)=\g$ is the same for all $\g \in G$.

If $G$ is a subgroup of $\PGL_2(\F_q)$ and $\g \in G$, let 
$$N_{\g,G} = \#\{ \tau\in \F_q\cup\{\infty\}: \text{$\tau$ is regular and $\inv(\tau,q)=\calC_{\g,G}$}\},$$ 
{\it i.e.}, the number of $\tau\in \F_q\cup\{\infty\}$ such that $|Q^{-1}(\tau)|=|G|$ and $v^q=\g(v)$ for some $v\in Q^{-1}(\tau)$.
By (\ref{eq:aQb}), $N_{\g,G}$ does not depend on the particular choice of quotient map. 
Since $G$ maps $\F_q \cup \{\infty\}$ to itself, it is a union of $G$-orbits, and 
$N_{\textmatrix 1001,G}$ is the number of full-sized $G$-orbits in $\F_q\cup\{\infty\}$. 
If $G$ is abelian then $\sum_{\g\in G} N_{\g,G}$ is equal to the total number of regular elements in $\F_q\cup\{\infty\}$ wrt any quotient map.

\begin{proposition} \label{prop:cyclicCount}
	If $G$ is a cyclic subgroup of $\PGL_2(\F_q)$ of order $\ell\ge2$ and $M$ is a generator of $G$, then $N_{\g,G}=(q+\kappa)/\ell$ for all $\g\in G$, where $\kappa\in\{1,0,-1\}$ and 
	$1-\kappa$ is the number of $\tau\in\F_q\cup\{\infty\}$ that are fixed by~$M$.
\end{proposition}
\begin{proof} By Lemma~\ref{lem:conjugateInv}, if $G'=\a G \a^{-1}$ where $\a \in \PGL_2(\F_q)$ then $N_{\g,G} = N_{\a \g \a^{-1},G'}$ for all $\g\in G$.
	Also, the value $\kappa$ is the same for both $M$ and $\a M \a^{-1}$. Thus, the assertion holds for $G$ iff it holds
	for $\a G \a^{-1}$.  Then Proposition~\ref{prop:Dickson} reduces to the three cases (a) $M=\textmatrix 1b01$;
	(b) $M=\textmatrix a001$; or (c) $M=E_{\z,\l}$.  These three cases correspond to $\kappa=0$, $-1$, and 1, respectively.  The short orbits are
	$\{\infty\}$ in case~(a); $\{\infty\}$ and $\{0\}$ in case~(b); and $\{\l\}$ and $\{\l^q\}$ in case~(c).  In particular, $\F_q\cup\{\infty\}$ contains $1-\kappa$ short elements,
	so the number of full-sized orbits in $\F_q\cup\{\infty\}$ is $((q+1)-(1-\kappa))/\ell=(q+\kappa)/\ell$. Thus, $N_{\textmatrix 1001,G} = (q+\kappa)/\ell$.
	In particular, $\ell$ divides $q+\kappa$.
	Proposition~\ref{prop:main}{\it (v)} implies that $N_\g=(q+\kappa)/\ell$ when $o(\g)\ge 3$.  
	This completes the proof when $\ell$ is odd.  
	If $\ell$ is even, then $G$ contains exactly one element $\g_2$ of order~2.  The irregular elements of $\cj \F_q\cup \{\infty\}$ (with respect to the quotient maps described in
	this section) are $\{\infty\}$ in case~(a), $\{0,\infty\}$
	in case~(b), and $\{\l,\l^q\}$ in case~(c), so the total number of regular elements in $\F_q\cup \{\infty\}$ is $q+\kappa$.
	Then $\sum_{\g\in G} N_\g = q+\kappa$. We have already shown $N_\g=(q+\kappa)/\ell$ when $\g=1$ or $\circ(\g)\ge 3$,
	\ie, for all $\g\in G$, $\g\ne \g_2$. Then $\N_\g = (q+\kappa)/\ell$ for $\g=\g_2$ also.
\end{proof}

\subsection{The equation $v^q=\g(v)$.} \label{subsec:vqgamma}
Let $\g=\textmatrix abcd\in\PGL_2(\F_q)$ have order $t$, and consider the equation $v^q=\g(v)$, where $v \in \cj\F_q \cup \{\infty\}$. By convention, we interpret
$\infty^q=\infty$.  Denote the solution set by $S_{\g,q}$.  Then $\infty\in S_{\g,q}$ iff $c=0$,
and the remaining solutions are roots of the polynomial $f(x)=x^q(cx+d)-(ax+b)$.  As shown in the proof of Proposition~\ref{prop:counting}{\it (ii)}, this equation has no repeated roots. 
Thus, $|S_{\g,q}|=q+1$, where in the case $c=0$ the
``$+1$'' accounts for $v=\infty$.  Note that $S_{\g,q}\setminus \{\infty\}$ is the same as the set $A_{\g,q}$ that was studied in Section~\ref{sec:orbits}. If $d\ge 1$, let
$$S_{\g,q}^{(d)} = \{v \in S_{\g,q} : \deg_q(v)=d \},$$
where $\deg_q(v)=[\F_q(v):\F_q]$ when $v\in \cj\F_q$ and $\deg_q(\infty)=1$.  

\begin{lemma} \label{lem:SgammaConj} If $\a,\g\in \PGL_2(\F_q)$ and $d\ge1$, then $S_{\a \g \a^{-1},q}^{(d)} = \{\a(v) : v \in S_{\g,q}^{(d)} \}$.
\end{lemma}
\begin{proof} If $v\in\cj\F_q\cup\{\infty\}$ and $w=\a(v)$, then $\deg_q(w)=\deg_q(v)$, and
	$$v\in S_{\g,q} \iff v^q=\g(v) \iff (\a v)^q=\a\g(v) \iff w^q= \a\g\a^{-1}(w) \iff w \in S_{\a\g\a^{-1},q}.$$
\end{proof}

If $1\ne \g \in \PGL_2(\F_q)$, let
$$Z_\g=\{\a\in\PGL_2(\F_q) : \a \g = \g \a \}.$$ 
Lemma~\ref{lem:SgammaConj} immediately implies that $Z_\g$ acts on $S_{\g,q}^{(d)}$.  The next lemma is well known.

\begin{lemma} \label{lem:Zgamma}
	If $1\ne \g \in \PGL_2(\F_q)$, let $1-\kappa$ denote the number of fixed points of $\g$ in $\F_q\cup\{\infty\}$. Then $\kappa \in \{0,1,-1\}$, and 
	$|Z_\g|=q+\kappa$.
\end{lemma}
\begin{proof} Write $\g\sim\b$ to denote that $\g$ is conjugate to $\b$. By Proposition~\ref{prop:Dickson}, there are three cases: (a)\ $\g \sim \textmatrix 1b01$ for
	some $b \in \F_q^\x$ and $\kappa=0$; (b)\ $\g \sim \textmatrix a001$ for some $1\ne a\in\F_q^\x$ and $\kappa=-1$; or (c)\ $\g \sim D_{\d,\l}$ for some $\d,\l \in\F_{q^2}\setminus \F_q$ and
	$\kappa=1$.
	In case~(a), it is easy to see that $\textmatrix rstu$ commutes with $\textmatrix 1b01$ (in $\PGL_2$, \ie, up to a scalar multiple)  iff $t=0$ and $r=u$, therefore
	$Z_\g \cong Z_{\textmatrix 1b01}=\{\textmatrix 1 e01 : e \in \F_q\}$, and $|Z_\g|=q$. In case~(b), $\textmatrix rstu$ commutes with $\textmatrix a001$ in $\PGL_2$
	iff $s=t=0$, so
	$Z_\g\cong Z_{\textmatrix a001} = \{\textmatrix e001 : e \in \F_q^\x\}$,
	and $|Z_\g|=q-1$. In case~(c), $C_\l \textmatrix rstu C_\l^{-1}$ commutes with $D_{\d,\l}$ in $\PGL_2$ iff $\textmatrix rstu$ 
	commutes with $\textmatrix {\d^q}00\d$ in $\PGL_2$ iff
	$s=t=0$, and $C_\l \textmatrix r00u C_\l^{-1}$ is rational in $\PGL_2$ iff $C_{\l^q} \textmatrix {r^q}00{u^q}C_{\l^q}^{-1} = C_\l \textmatrix {kr} 00{ku} C_\l^{-1}$ with $k\ne0$
	iff $\textmatrix {u^q}00{r^q} = \textmatrix {kr}00{ku}$
	iff $(u/r)^q=r/u $ iff $u^{-1} C_\l \textmatrix r00u C_\l^{-1} = E_{\z,\l}$ with $\z=r/u \in \mu_{q+1}$. Thus, $Z_\g 
	= \{E_{\z,\l} : \z^{q+1}=1\}$ and $|Z_\g|=q+1$.  In each case, $|Z_\g|=q+\kappa$.
\end{proof}

\begin{proposition} \label{prop:vqgammav} Let $1\ne \g\in\PGL_2(\F_q)$ have order~$t$, and let $\kappa$ be as in Lemma~\ref{lem:Zgamma}. Then 
	$S_{\g,q}=S_{\g,q}^{(1)} \cup S_{\g,q}^{(t)}$, and $|S_{\g,q}^{(t)}|=q+\kappa$. If $v$ is any element of $S_{\g,q}^{(t)}$, then
	$S_{\g,q}^{(t)}=\{z(v) : z \in Z_\g\}$. 
\end{proposition}
\begin{proof}
	By Lemma~\ref{lem:SgammaConj}, if the proposition holds for $\g$ then it also holds for $\a\g\a^{-1}$, where $\a\in\PGL_2(\F_q)$.  Using Proposition~\ref{prop:Dickson},
	we may therefore assume that one of the three cases holds: (a)\ $\g=\textmatrix 1b01$, $b\in\F_q^\x$, $\kappa=0$; (b)\ $\g=\textmatrix a001$, $1\ne a \in \F_q^\x$, $\kappa=-1$;
	or (c)\ $\g=D_{\d,\l}$ and $\kappa=1$, where $\d,\l\in\F_{q^2}\setminus \F_q$. 

	In case~(a), $S_{\g,q}$ contains $\{\infty\}$, together with all $v\in\cj\F_q$ such that
	$v^q=v+b$. 
	Since $v^{q^i}=v+ib$, $v$ has exactly $p$ conjugates, where $p$ is the prime dividing~$q$, and so $\deg_q(v)=p=\circ(\g)$.
	If $v_0\in\cj\F_q$ is one solution to $v^q=v+b$, then the others are $v_0+c$ for $c\in\F_q$. Since $Z_\g=\{\textmatrix 1c01 : c \in \F_q\}$, the proposition holds in this case.

	In case~(b), $S_{\g,q}$ consists of $0,\infty$, and the nonzero solutions to $v^q=av$.  If $v\ne0$, then the distinct conjugates of $v$ are
	$v^{q^i}=a^iv$ for $0\le i<\circ(a)$, so $\deg_q(v)=\circ(a)=\circ(\g)$.
	If $v_0$ is one nonzero solution, so $v_0^{q-1}=a$, then the others
	are $cv_0$ for $c\in\F_q^\x$. Since $Z_\g=\{\textmatrix c001:c\in\F_q^\x\}$, the nonzero solutions form a single $Z_\g$-orbit. This analysis shows that $S_{\g,q}$
	contains two elements of $S_{\g,q}^{(1)}$ and the remaining $q-1$ elements comprise a single $Z_\g$ orbit of size $|Z_\g|=q-1$.

	In case~(c), we may write $\g=E_{\z_0,\l}$, and $Z_\g = \{E_{\z,\l} : \z^{q+1}=1\}$. As shown in Lemma~\ref{lem:AppBfixed}, if $1\ne \a\in Z_{\g}$ then $\l$ and $\l^q$
	are its only fixed points in $\cj \F_q \cup \{\infty\}$. In particular, $\a$ has no fixed points in $\F_q$, so $v^q=v=\a(v)$ has no solutions. Taking $\a=\g$, this implies that
	$S_{\g,q}^{(1)}=\emptyset$. If $v\in\{\l,\l^q\}$ then $\g(v)=v\ne v^q$, so $v \not \in S_{\g,q}$.  
	Now let $v$ be any element of $S_{\g,q}$. We have shown that $v\not\in\F_q\cup\{\infty\}\cup\{\l,\l^q\}$. Thus, $\deg_q(v)>1$ and $E_{\z,\l}(v)\ne v$ for every $E_{\z,\l}\in Z_\g$.
	Then $\{\a(v) : \a \in Z_\g\}$ are distinct. Calling this set $S$, we have $|S|=|Z_\g|=q+1$. Also, $S\subset S_{\g,q}$ by Lemma~\ref{lem:SgammaConj}. Since both have cardinality
	$q+\kappa$, $S=S_{\g,q}$.
	Since $\a(v)$ are distinct for $\a\in Z_\g$, $\g^i(v)$ are distinct for $0\le i<\circ(\g)$. Then $\deg_q(v)=\circ(\g)$ by Lemma~\ref{lem:degree}.
\end{proof}

\begin{corollary} \label{cor:factoring1} Let $1\ne \g=\textmatrix abcd \in \PGL_2(\F_q)$, let $t=\circ(\g)$, and let $1-\kappa$ be the number of fixed points of $\g$ in $\F_q\cup\{\infty\}$. 
	Then $\kappa\in \{0,1,-1\}$, and the polynomial $x^q(cx+d)-(ax+b)\in\F_q[x]$
	factors into exactly $(q+\kappa)/t$ irreducible polynomials of degree~$t$. The remaining factors are linear. 
	If $r$ is one irrational root of $f$, then the others are $\{z(r): z\in Z_\g\}$.
\end{corollary}
\begin{proof} The roots of $f$ are the finite elements of $S_\g$.  Each degree-$t$ factor of $f$ corresponds to $t$ conjugate roots in $S_{\g,q}^{(t)}$. The result now follows from
	Proposition~\ref{prop:vqgammav}.
\end{proof}

\section{$G=\PGL_2(\F_q)$.} \label{sec:PGL2}

This section considers the case $G=\PGL_2(\F_q)$, and we prove the statements from the introduction (Example~\ref{example:pgl2}). As shown in (\ref{eq:QGidentity}),
a quotient map is
\begin{equation*} Q(x) =  1 + \frac {x^{q^3} - x } {\left(x^q-x\right)^{q^2-q+1}} 
= \frac{(x^{q^2}-x)^{q+1}}{(x^q-x)^{q^2+1}}.\end{equation*}
As usual, we begin by considering short orbits, \ie, orbits of size less than $|G|$. Recall from Section~\ref{sec:orbits} that $|G|=q^3-q$.

\begin{lemma} \label{lem:PGL2Short} 
$v \in \cj\F_q$ belongs to a short orbit of $\PGL_2(\F_q)$ if and only if $v \in \F_{q^2}$.  The orbit of $\infty$ is $\F_q \cup \{\infty\}$.
\end{lemma}

\begin{proof} $\PGL_2(\F_q)$ maps $\F_{q^2}\cup\{\infty\} $ to itself. Since $q^2+1 < q^2+q \le q(q-1)(q+1)$, each element of $\F_{q^2}$ belongs to a short orbit.
Conversely, all elements of short orbits are in $\F_{q^2} \cup \{\infty\}$ by Lemma~\ref{lem:short}.
For the last statement, we know $\PGL_2(\F_q)$ preserves $\F_q \cup \{\infty\}$. It is a single orbit, because if $a \in \F_q$ then
$\textmatrix a 1 1 0(\infty) = a$.
\end{proof}

\begin{lemma}  \label{lem:pgl2Q} 
	The only irregular elements in $\cj\F_q\cup\{\infty\}$ with respect to $Q$ are~0 and~$\infty$. Their preimages are the short $G$-orbits
	\begin{equation} Q^{-1}(\infty)=\F_q\cup\{\infty\},\qquad Q^{-1}(0) = \F_{q^2}\setminus \F_q. \label{eq:deg2Orbit} \end{equation}
\end{lemma}

\begin{proof} By Lemma~\ref{lem:PGL2Short}, the union of the short orbits is $\F_{q^2} \cup \{\infty\}$, and $\calO_\infty=
	\F_q \cup \{\infty\}$. The images of the short orbits under $Q$ are the irregular elements. If $v \in\calO_\infty$ then $Q(v)=\infty$,
	and if $v \in \F_{q^2}\setminus \F_q$ then $Q(v)=0$.
	Then $\infty$ and 0 are the only irregular elements of $\cj\F_q\cup\{\infty\}$, and $\F_{q^2}\cup\{\infty\}$
	is the union of exactly two short orbits:
	$Q^{-1}(\infty)=\F_q\cup\{\infty\}$ and $Q^{-1}(0)=\F_{q^2}\setminus \F_q$. 
\end{proof}

\begin{lemma}  \label{lem:ge3}
Let $\tau\in\F_q^\x$.  Then $\tau$ is regular, and $\inv(\tau) = \calC_\g$ has the property that $\circ(\gamma)\ge3$. 
\end{lemma}
\begin{proof} By Lemma~\ref{lem:pgl2Q},  $\tau$ is regular. 
Then $Q^{-1}(\tau)$ is a full-sized orbit of $\PGL_2(\F_q)$, and $\inv(\tau)$ is defined as the unique conjugacy
class $\calC \subset \PGL_2(\F_q)$ such that $v^q = \gamma(v)$ with $\gamma \in \calC$ whenever $v \in Q^{-1}(\tau)$.
Since all elements of $\F_{q^2}$ are in short orbits, $Q^{-1}(\tau)$ misses $\F_{q^2}$, so $\deg_q(v)\ge 3$ for each $v \in Q^{-1}(\tau)$.
By Lemma~\ref{lem:degree}, $\deg_q(v) = \circ(\gamma)$.
Thus, $\inv(\tau) = \calC_\g$ always has the property that $\circ(\g)\ge 3$. 
\end{proof}

\begin{lemma} \label{lem:xyz}
If $\beta = \textmatrix abcd \in \PGL_2(K)$ then for any $x,y,z$ in a field containing $K$,
\begin{equation} \beta(x)-\beta(z) = \frac{(ad-bc)(x-z)}{(cx+d)(cz+d)} \label{eq:xz} \end{equation}
\begin{equation} \frac{\beta(x)-\beta(z)}{\beta(y)-\beta(z)} = \frac{(x-z)(cy+d)}{(y-z)(cx+d)} \label{eq:xyz} \end{equation}
\end{lemma}
\begin{proof} This is a straightforward computation. \end{proof}

Recall in (\ref{eq:iota2}) we defined
$$\iota\left(\textmatrix abcd\right) = \frac{e_1}{e_2} + \frac{e_2}{e_1} + 2 = \frac{(a+d)^2}{ad-bc}$$
where $e_1,e_2$ are the roots of the characteristic equation of $\textmatrix abcd$. We noted that $\iota$ is well defined on $\PGL_2(\F_q)$
and is constant on conjugacy classes.  Also, recall from Proposition~\ref{prop:BDD} (with $c=1$) that
\begin{equation} \F_q = \{ \z + 1/\z : \z \in \mu_{q-1}\cup \mu_{q+1} \}. \label{eq:zeta} \end{equation}

\begin{theorem} \label{thm:pgl2bijection}
The map $\iota$ induces a bijection between conjugacy classes $\calC_\g$ in $\PGL_2(\F_q)$ such that $\circ(\g)\ge 3$ and 
$\F_q^\x$.
	Further, $\iota(\gamma)=\tau$ iff $\inv_Q(\tau) = \calC_\g$. That is, $\inv_Q$ is the inverse bijection to $\iota$. 
\end{theorem}

\begin{proof}
Let $\gamma\in \PGL_2(\F_q)$ have order $\ell \ge 3$. Then case~(a), (b), or~(c) of Proposition~\ref{prop:Dickson} holds. 

In case~(a), 
	$\calC_\g = \calC_{\textmatrix 1b01}$ where $b\in\F_q^\x$, $\ell =p\ge 3$, and $\iota(\gamma) = 4$. 
Let $v \in \cj\F_q$ such that $v^q = \textmatrix 1b01 v = v +b$. Then
	$$Q(v)=\frac{(v^{q^2}-v)^{q+1}}{(v^q-v)^{q^2+1}} = \frac{((v+2b)-v )^{q+1}}{((v+b)-v)^{q^2+1}} = \frac{(2b)^2}{b^2}=4= \iota(\g).$$

	In case~(b), $\calC_\g = \calC_{\textmatrix a001}$ where $a \in \F_q^\x$ and $\circ(a)\ge 3$.
Let $v \in \cj\F_q^\x$ such that $v^q = \textmatrix a001 v = av $. Then
	\begin{eqnarray*} Q(v) &=& \frac{ (v^{q^2}-v)^{q+1}}{(v^q-v)^{q^2+1}} = \frac{ (a^2 v-v)^{q+1}}{(av-v)^{q^2+1}} = \frac{ v^{q+1}(a^2-1)^2}{v^{q^2+1}(a-1)^2} \\
		&=& \frac {(a+1)^2}{v^{q^2-q}} = \frac{(a+1)^2}{a^q}=\frac{(a+1)^2}a = \iota(\g).
	\end{eqnarray*}

	In case~(c), $\g = E_{\z,\l}$ for some $\l \in \F_{q^2} \setminus \F_q$ and $\z$ of order $\ell$. Here, $\iota(\gamma)=\iota\left(\textmatrix \z 001\right) = (\z+1)^2/\z$.
Let $v \in \cj\F_q$ such that 
$v^q = \g(v)$. Since $\g$ is rational as an element of $\PGL_2(\F_q)$, $v^{q^i} = \g^i(v) = E_{\z^i,\l}(v)$.
Let $u = C_\l^{-1}(v)$. Then $v^{q^i} = E_{\z^i,\l}(v) = C_\l(\z^i u)$. Using Lemma~\ref{lem:xyz},
\begin{eqnarray*} 
	Q(v) &=& \frac{(v^{q^2}-v)^{q+1}}{(v^q-v)^{q^2+1}}  = 
	 \left(\frac{v^{q^3}-v^q}{v^{q^3}-v^{q^2}}\right)  \left(\frac{v^{q^2}-v}{v^{q}-v}\right) \\
&=& \left(\frac{C_\l\left(\zeta u\right) - C_\l(\z^3 u)}{C_\l\left(\zeta^2 u \right)-C_\l(\z^3 u)}\right)
\left(\frac{C_\l(\z^2 u) - C_\l(u) } {C_\l(\z u) - C_\l(u) }\right) \\
&=& \frac{(\zeta u - \z^3 u)(\z^2 u - 1)}{(\z^2 u - \z^3u)(\z u - 1)} \frac{(\z^2 u - u)(\z u - 1)}{(\z u -  u)(\z^2 u - 1)} \\
	&=& \frac{ (\z - \z^3)(\z^2-1)}{(\z^2-\z^3)(\z-1)} = \frac{(\z+1)^2}{\z} = \iota(\g).
\end{eqnarray*} 

Combining the three cases, we see that if $v^q = \g(v)$ and $\circ(\g)\ge 3$ then $Q(v)=\iota(\g)$.  In each case, $\tau=\iota(\g) \in \F_q^\x$,
so it is regular and $\inv(\tau)$ is defined. Since $Q(v)=\tau$ and $v^q=\g(v)$, $\inv(\tau) = \calC_\g$.  We have shown $\inv\circ \iota$ is
the identity on $\{\calC_\g : o(\g) \ge 3 \}$.  

To prove that $\iota$ and $\inv_Q$ are bijections, it remains to prove that $\iota$ is surjective
from $\{\calC_\g : \circ(\gamma)\ge 3 \}$ onto $\F_q^\x$. Let $\tau \in \F_q^\x$. By (\ref{eq:zeta}), $\tau-2=\z + 1/\z$ where $\z^{q-1}=1$
or $\z^{q+1}=1$. If $\z=1$ then $\tau=4=\iota\left(\textmatrix 1101\right)$. In even characteristic, $4=0 \not \in \F_q^\x$. In odd characteristic, $\circ\textmatrix 1101
=p\ge 3$. If $\z=-1$ then $\tau=0 \not \in \F_q^\x$. If $\z^{q-1}=1$ and $\z\not\in \{1,-1\}$ then $\tau = \iota(\g)$ for $\g=\textmatrix \z001$, and $\circ(\gamma)
\ge 3$. Finally, if $\z^{q+1}=1$ and $\z \not \in \{1,-1\}$ then $\tau = \iota(E_{\z,\l})$ and $\circ(E_{\z,\l})\ge3$.  Thus, $\iota$ is surjective and the theorem is proved.
\end{proof}

\begin{corollary}  \label{cor:PGL2inv} 
	Let $\tau \in \F_q$.  By (\ref{eq:zeta}), $\tau-2 = \zeta+1/\zeta$, where $\zeta^{q-1}=1$ or $\zeta^{q+1}=1$.
Let $Q$ be the quotient map for $\PGL_2(\F_q)$ given in (\ref{eq:QGidentity}), and let $\inv = \inv_Q$. 
\begin{enumerate} 
\item[{\it (i)}] If $\zeta = -1$ (so $\tau = 0$) then $Q^{-1}(\tau) = \F_{q^2} \setminus \F_q$, and $\inv(\tau)$ is undefined since this orbit is short.
\item[{\it (ii)}] If $\zeta = 1$ and $1 \ne -1$ (so $\tau = 4$ and $q$ is odd), then $\inv(\tau) = \calC_{\textmatrix 1101}$.
\item[{\it (iii)}] If $\zeta \not \in \{1,-1\}$ and $\zeta^{q-1}=1$, then $\inv(\tau) = \calC_{\textmatrix \z 001}$.
\item[{\it (iv)}] If $\zeta \not \in \{1,-1\}$ and $\zeta^{q+1}=1$, then $\inv(\tau) = \calC_\gamma$, where $\gamma = E_{\z,\lambda}$ 
for any $\lambda \in \F_{q^2}\setminus \F_q$.
(All such matrices $E_{\z,\lambda}$ belong to the same conjugacy class in $\PGL_2(\F_q)$.)
\end{enumerate}
\end{corollary}
\begin{proof} {\it (i)} was shown in  (\ref{eq:deg2Orbit}), {\it (ii)}-{\it (iv)} follow from Theorem~\ref{thm:pgl2bijection}.
\end{proof}

\begin{corollary} \label{cor:order_gamma} 
Let $1\ne \g \in \PGL_2(\F_q)$, where $q=p^e$ and $p$ is any prime.   \\
	{\it (i)} $\iota(\g)=0$ if and only if $\circ(\g)=2$. \\
	{\it (ii)} $\iota(\g)=4$ if and only if $\circ(\g)=p$. \\
	{\it (iii)} Write $\iota(\g) = \z + 1/\z + 2$, where $\z \in \mu_{q-1}\cup \mu_{q+1}$. If $\z\ne 1$  (equivalently, $\iota(\g)\ne 4$), then $\circ(\g)=\circ(\z)$.
\end{corollary}
\begin{proof} {\it (i)}\ If $\g = \textmatrix abcd$ then $$\g^2=\begin{pmatrix} a^2+bc & b(a+d) \\ c(a+d) & d^2+bc \end{pmatrix}.$$
This matrix is scalar iff $b(a+d)=c(a+d)=a^2-d^2=0$.  These equations hold iff $a+d=0$ or $b=c=a-d=0$. The latter is excluded since we assume $\g \ne 1$.

{\it (ii)} and {\it (iii)}\  Assume first that $\iota(\g)\ne0$ and write $\iota(\g)=\z+1/\z+2$, where $\z\ne -1$. 
By {\it (i)}, $\circ(\g)\ge 3$. By Theorem~\ref{thm:pgl2bijection}, if $\iota(\a)=\iota(\g)$ and $\circ(\a)\ge 3$, then $\calC_\g=\calC_{\a}$, so $\circ(\g)=\circ(\a)$. 
If $\z=1\ne-1$, then $p$ is odd and $\iota(\g)=4=\iota(\textmatrix 1101)$. Since $\circ(\textmatrix 1101)=p\ge3$, $\circ(\g)=p$.
If $\z \in\mu_{q-1}\setminus \mu_2$ then $\iota(\g)=\iota(\textmatrix \z001)$, and $\circ(\textmatrix \z001)=\circ(\z)\ge3$, so $\circ(\g)=\circ(\z)$. Finally, if $\z \in \mu_{q+1}\setminus \mu_2$
then $\iota(\g)=\iota(E_{\z,\l})$ and  $\circ(E_{\z,\l})=\circ(\z)\ge3$, so $\circ(\g)=\circ(\z)$.

The proofs of {\it (ii)} and {\it (iii)} are complete when $\iota(\g)\ne 0$.  In ({\it ii}), $\iota(\g)=0$ iff $4=0$ iff $p=2$. In that case
{\it (ii)} follows from {\it (i)}. In {\it (iii)}, $\iota(\g)=0$ iff $\z=-1$. Since the case $\z=1$ is excluded, $q$ must be odd. Then  $\circ(\g)=2$ by {\it (i)}, but also $\circ(\z)=\circ(-1)=2$.
\end{proof}

\begin{corollary} \label{cor:PGL2conjugate}
	If $\g,\g' \in \PGL_2(\F_q) \setminus \{1\}$ and $\iota(\g)=\iota(\g')\ne 0$ then $\g$ and $\g'$ are conjugate. In particular, if $\g\ne1$ and $\iota(\g)\ne0$
	then $$\calC_\g = \{\a \in \PGL_2(\F_q) : \text{$\a\ne1$ and $\iota(\a)=\iota(\g)$}\}.$$
\end{corollary}
\begin{proof} By Corollary~\ref{cor:order_gamma}{\it (i)}, the hypothesis implies that $\circ(\g)\ge 3$ and $\circ(\g')\ge3$. Then $\iota(\g)=\iota(\g')$ implies $\calC_\g=\calC_{\g'}$
	by Theorem~\ref{thm:pgl2bijection}.
	\end{proof}

	\begin{proposition} \label{prop:factoring} Let $q=p^e$ where $p$ is prime, and let $\textmatrix abcd \in \GL_2(\F_q)$ be a nonscalar matrix.
		Let $f(x)=x^q(cx+d)-(ax+b)$, and write $(a+d)^2/(ad-bc)-2=\z+1/\z$ where $\z\in\mu_{q-1}\cup\mu_{q+1}$. \\ 
		{\it (i)}\ If $\z=1$ then $f$ has exactly $p^{e-1}$ irreducible factors of degree~$p$, and the remaining factors are linear.  \\
		{\it (ii)} If $\z=-1$ and $1\ne-1$ (so $q$ is odd), then $f$ has exactly $(q+\kappa)/2$ irreducible quadratic factors and the remaining factors are linear, where 
		$\kappa=-\jacobi{-(ad-bc)}q$. \\ 
		{\it (iii)}  If $\z\in\mu_{q+\kappa}\setminus \mu_2$, where $\kappa\in\{1,-1\}$, then $f$ has $(q+\kappa)/t$ irreducible factors of degree~$t$ and the remaining factors are linear, 
		where $t=\circ(\z)$. \\
		In each case, if $v$ is one irrational root of $f$, then the others are $\a(v)$ such that $\a \in \PGL_2(\F_q)$ and $\a\textmatrix abcd=\textmatrix abcd\a$.
	\end{proposition}

	\begin{proof} Let $\g=\textmatrix abcd$, considered as an element of $\PGL_2(\F_q)$, and let $1-\kappa$ denote the number of fixed points of $\g$ in $\F_q\cup\{\infty\}$.
		Note that $\iota(\g)=2+\z+1/\z=(\z+1)^2/\z$.
		By Corollary~\ref{cor:factoring1}, $f$ has $(q+\kappa)/\circ(\g)$ irreducible factors of degree~$\circ(\g)$, the remaining factors are linear, and the irrational
		roots comprise a $Z_\g$-orbit of full size.
		So to prove the proposition, it suffices to compute $\circ(\g)$ and $\kappa$ in each of the cases {\it (i)}--{\it (iii)}.

		{\it (i)}\ Given that $\z=1$, we must show $\circ(\g)=p$ and $\kappa=0$.  $\iota(\g)=4$ and so $\circ(\g)=p$ by Corollary~\ref{cor:order_gamma}{\it (ii)}.  
		If $p$ is odd, then $\g\sim\textmatrix 1b01$ by Corollary~\ref{cor:PGL2conjugate}, 
		so $\g$ has a unique fixed point in $\F_q\cup\{\infty\}$ and $\kappa=0$. If $p=2$, then $\iota(\g)=4=0$, so $\g=\textmatrix abc{-a} =\textmatrix abca$. Note that $b$, $c$
		cannot both be zero, as otherwise $\g$ would be scalar. It is easy to see that $\g$ has a unique fixed point $(b/c)^{1/2}$ if $c\ne0$, or $\infty$ if $c=0$. 
		Thus, $1-\kappa=1$ and $\kappa=0$.  We have shown $t=p$ and $\kappa=0$ for any $q$, even or odd, as required.

		{\it (ii)}\ Given that $\z=-1$ and $p$ is odd, we must show $\circ(\g)=2$ and $\kappa=-\jacobi{-\det(\g)}q$. Since $\iota(\g)=2+\z+1/\z=0$, $\g=\textmatrix abc{-a}$, and 
		$\circ(\g)=2$ by Corollary~\ref{cor:order_gamma}{\it (i)}. 
		The fixed points of $\g$ are the roots of $cz^2-2az-b$, and the number of rational roots is $1+\jacobi{4a^2+4bc}q=1+\jacobi{-\det(\g)}q$. Thus, there are $1-\kappa$ fixed points
		in $\F_q\cup\{\infty\}$, where $\kappa=-\jacobi{-\det(\g)}q$, as was to be shown. 

		{\it (iii)}\ The hypothesis is that $\iota(\g)=(\z+1)^2/\z$ where $\z\in\mu_{q-1}\cup\mu_{q+1}$ and $\z^2\ne 1$.
		By Corollary~\ref{cor:order_gamma}{\it (iii)}, $\circ(\g)=\circ(\z)$, and by Corollary~\ref{cor:PGL2conjugate},
		$\g \sim \textmatrix \z001$ if $\z^{q-1}=1$, and $\g\sim E_{\z,\l}$ if $\z^{q+1}=1$. In the former case, $\kappa=-1$, and in the latter case, $\kappa=1$.
		{\it (iii)} now follows from Corollary~\ref{cor:factoring1}.
	\end{proof}

We conclude this section by proving the second theorem from Example~\ref{example:pgl2}.   Let $G=\PGL_2(\F_q)$, $Q_G$ the quotient map given by~(\ref{eq:QGidentity}), $H$ a subgroup of $G$,
and $Q_H$ a quotient map for~$H$.
By Proposition~\ref{prop:H}, there is a unique function $h\in \F_q(x)$ such that $Q_G = h\circ Q_H$.

\begin{theorem} \label{thm:Fq} Let $H\subset \PGL_2(\F_q)$ and $Q_H,h$ be as above. Suppose $\tau \in \F_q$ is regular with respect to $Q_H$ and let $\inv_{Q_H}(\tau,q)=\calC_{\g,H}$. 
	If $\g=1$ then $h(\tau)=\infty$. If $\g\ne1$ then $h(\tau)=\iota(\g)$.
\end{theorem}

\begin{proof}  By Proposition~\ref{prop:Qorbit}, $V=Q_H^{-1}(\tau)$ is an $H$-orbit, and since $\tau$ is regular with respect to $Q_H$, the orbit has full size.
	Let $v\in V$. Since $Q_H(v^q)=\tau^q=\tau$, $v^q \in V$, and consequently $v^q=\d(v)$ for a unique $\d \in H$. By~(\ref{eq:invariantFq}), $\inv_{Q_H}(\tau)=\calC_{\d,H}$, therefore
	$\d$ is conjugate to $\g$. In particular, $\d$ has the same order as $\g$ and $\iota(\d)=\iota(\g)$. 
	Note that $Q_G(v)=h(Q_H(v))=h(\tau)$. By (\ref{eq:deg2Orbit}), $Q_G^{-1}(\infty)=\F_q\cup\{\infty\}$ and $Q_G^{-1}(0)=\F_q^2\setminus \F_q$. 

	First, $\g=1\Rightarrow \d=1\Rightarrow v\in\F_q\Rightarrow h(\tau)=Q_G(v)=\infty$.

	Next, suppose $\circ(\g)=2$. Since $v,\g(v)=v^q$ are distinct and $v^{q^2}=\g^2(v)=v$, $v$ belongs to $\F_{q^2}\setminus \F_q$. Then $h(\tau)=Q_G(v)=0$. On the other hand, 
	$\circ(\g)=2\iff \iota(\g)=0$. So $h(\tau)=\iota(\g)=0$ in this case.

	Finally, if $\circ(\g)\ge 3$ then $[\F_q(v):\F_q]=\circ(\d)=\circ(\g)\ge3$ by Lemma~\ref{lem:degree}, so $v \not\in\F_{q^2}$. Then $h(\tau)=Q_G(v)\in\F_q^\x$, so that $h(\tau)$ is regular with respect to $Q_G$. 
	Lemma~\ref{lem:Hinv} then implies that $\inv_{Q_G}(h(\tau))= 
	\calC_{\g,G}$. Finally, Theorem~\ref{thm:pgl2bijection} implies $h(\tau)=\iota(\gamma)$.
\end{proof}
	
\section{$G=\PSL_2(\F_q)$.} \label{sec:PSL2}

The projective special linear group is defined as $\SL_2(\F_q)$ modulo the scalar matrices $\textmatrix a00a\in \SL_2(\F_q)$.
If $\a \in GL_2(\F_q)$ and $\det(\a)=c^2$ with $c \in \F_q$, then $c^{-1}\a \in \SL_2(\F_q)$, so ($\a$ mod scalars) represents an element of $\PSL_2(\F_q)$.
On the other hand, if $\det(\a)$ is a nonsquare, then it has no scalar rational multiple in $\SL_2(\F_q)$. Thus, there is a short exact sequence
$$ 1 \longrightarrow \PSL_2(\F_q) \longrightarrow \PGL_2(\F_q) \longrightarrow \{\pm 1 \} \longrightarrow 1 $$
where the first map is inclusion and the second map is $\jacobi{\det(\a)}q$. The square-class of the determinant is well defined, because a
scalar matrix has square determinant. 

It follows that $[\PGL_2(\F_q):\PSL_2(\F_q)] = |\{\pm1\}|$, which is 1 if $q$ is even and 2 if $q$ is odd. In particular, $\PSL_2(\F_q)=\PGL_2(\F_q)$ when
$q$ is even, and we already studied this group in Section~\ref{sec:PGL2}.
For this reason, in this section we assume $q$ is odd.  Then $|\PSL_2(\F_q)| = (1/2)q(q-1)(q+1)$.

Usually we first find short orbits and then find the quotient map.  However, for this example it turns out to be easier to do these steps in reverse order.

\begin{proposition} \label{prop:PSL2Q} A quotient map for $\PSL_2(\F_q)$  is
	$$Q_S(x) = \frac{(x^{q^2}-x)^{(q+1)/2}}{(x^q-x)^{(q^2+1)/2}}.$$
If $\gamma \in \PGL_2(\F_q)$ then 
\begin{equation} Q_S\circ \gamma(x) = \jacobi{\det(\gamma)}q Q_S(x). \label{eq:QSgamma} \end{equation}
\end{proposition}

\begin{proof} 
First we prove (\ref{eq:QSgamma}). If the equation holds for $\gamma_1$ and for  $\gamma_2$, then it holds for $\gamma_1\circ \gamma_2$ as well, because
$$Q_S\circ (\gamma_1 \circ \gamma_2) = (Q_S \circ \gamma_1) \circ \gamma_2 = \jacobi{\det(\g_1)}q Q_S\circ \g_2 = \jacobi{\det(\g_1)}q \jacobi{\det(\g_2)}q Q_S,$$
and $\jacobi{\det(\g_1)}q \jacobi{\det(\g_2)}q = \jacobi{\det(\g_1) \det(\g_2)}q = \jacobi{\det(\g_1\g_2)}q$. 
Thus, it suffices to prove (\ref{eq:QSgamma}) for $\gamma =\textmatrix c001$,
$\gamma = \textmatrix 1b01$, and $\gamma = \textmatrix 0110$, as these generate $\PGL_2(\F_q)$.

If $\gamma = \textmatrix c001$ with $c \in \F_q^\x$ then
$$Q_S(\gamma (x) ) = Q_S(cx) = \frac{\left((cx)^{q^2}-cx\right)^{(q+1)/2}}{\left((cx)^q-cx\right)^{(q^2+1)/2}} 
= \frac{c^{(q+1)/2}(x^{q^2}-x)^{(q+1)/2}}{c^{(q^2+1)/2}(x^q-x)^{(q^2+1)/2}}.$$
The right side is $Q_S(x)$ times
\begin{eqnarray*} c^{(q-q^2)/2} &=&  (c^q)^{-(q-1)/2} = c^{-(q-1)/2} = c^{q-1} c^{-(q-1)/2} \\
&=& c^{(q-1)/2} = \jacobi cq = \jacobi {\det(\gamma)} q.
\end{eqnarray*}
 Thus, $Q_S(\gamma(x)) = \jacobi {\det(\gamma)} q Q_S(x)$.

If $\gamma = \textmatrix 1b01$ with $b\in\F_q$ then
$$Q_S(\gamma(x)) = Q_S(x+b) = \frac{ ((x+b)^{q^2} - (x+b))^{(q+1)/2}}{((x+b)^q-(x+b))^{(q^2+1)/2}} = Q_S(x)$$
and $\det(\gamma)=1$. 

Finally,  if $\gamma = \textmatrix 0110$ then
$$Q_S(\gamma(x)) = Q_S(1/x) = \frac{ (x^{-q^2} - x^{-1})^{(q+1)/2}}{(x^{-q}-x^{-1})^{(q^2+1)/2}}.$$
Multiply numerator and denominator by $x^{(q^2+1)(q+1)/2}$ to obtain 
$$ Q_S(1/x) = \frac{ (x - x^{q^2})^{(q+1)/2}}{(x-x^{q})^{(q^2+1)/2}} = (-1)^{(q-q^2)/2} Q_S(x).$$
Since $q^2 \equiv 1 \pmod 4$, $(-1)^{(q-q^2)/2}=(-1)^{(q-1)/2} = \jacobi{-1}q$.
Noting that $\det(\gamma)=-1$, the result follows.

Since $\det(\gamma)$ is a square for all $\gamma \in \PSL_2(\F_q)$, eq.~(\ref{eq:QSgamma}) shows
that $Q_S \circ \gamma = Q_S$ for all $\gamma \in \PSL_2(\F_q)$. 
Note that $Q_S^2=Q_G$, where
	$Q_G$ is the quotient map for $\PGL_2(\F_q)$ given by (\ref{eq:QGidentity}). Then
	$Q_S(\infty)=\infty$ and $\deg(Q_S)=(1/2)\deg(Q_G)=(1/2)|\PGL_2(\F_q)|=
	|\PSL_2(\F_q)|$.
	Thus, $Q_S$ is a quotient map for $\PSL_2(\F_q)$.

\end{proof}

\begin{lemma} \label{lem:PSL2short} Let $S= \PSL_2(\F_q)$. \\
{\it (i)}\ $\calO_\infty = \F_q \cup \{\infty\}$, and it has multiplicity $(1/2)(q^2-q)$.  \\
{\it (ii)}\ $\F_{q^2} \setminus \F_q$ is an $S$-orbit, and it has multiplicity $(1/2)(q+1)$. \\
{\it (iii)}\ $\calO_\infty$ and $\F_{q^2}\setminus \F_q$ are the only short orbits.  \\
{\it (iv)}\ The irregular elements of $\cj\F_q\cup\{\infty\}$ with respect to $Q_S$ are 0 and~$\infty$.\\
	{\it (v)}\ If $\tau \in \F_q^\x$ and $\inv_{Q_S}(\tau)=\calC_{\g,S}$ then $\circ(\g)\ge 3$ and $\iota(\g)=\tau^2$, where $\iota\left(\textmatrix abcd\right)=(a+d)^2/(ad-bc)$.
\end{lemma}
\begin{proof}
{\it (i)}\ Certainly $\calO_\infty \subset \F_q \cup \{\infty\}$. If $b\in \F_q$ then $b=\textmatrix 1b01 \textmatrix {0\,}{-1}{1\,}0 (\infty)$, so equality holds.
The multiplicity is $|S|/|\calO_\infty|=(1/2)(q^3-q)/(q+1)=(1/2)(q^2-q)$.

	{\it (ii)}\  Since $Q_S^2=Q_G$, $Q_S^{-1}(0)=Q_G^{-1}(0)$, which is $\F_{q^2}\setminus \F_q$ by~(\ref{eq:deg2Orbit}).
By Proposition~\ref{prop:Qorbit}, it is an $S$-orbit. The size of the orbit is $q^2-q$,
and the multiplicity is $|S|/(q^2-q)
= (1/2)(q^3-q)/(q^2-q)=(q+1)/2$.

{\it (iii)} Both $\calO_\infty$ and $\F_{q^2}\setminus \F_q$ are short as their multiplicities are greater than~1. There are no other short orbits
by Lemma~\ref{lem:short}.

{\it (iv)} holds because the images of the short orbits under $Q_S$ are $\infty$ and 0.

	{\it  (v)}  Write $\inv_{Q_S}(\tau)=\calC_{\g,S}$. Apply Theorem~\ref{thm:Fq}, with $h(x)=x^2$.  If $\g$ were the identity, then the theorem guarantees that
$h(\tau)=\infty$, and if $\circ(\g)=2$ then the theorem says $h(\tau)=0$.
However, $h(\tau)=\tau^2\in\F_q^\x$, so it must be that $\circ(\g)\ge 3$.
\end{proof}

\begin{theorem} \label{thm:QSinv} Let $\tau \in \F_q^\x$. Then $\inv_{Q_S}(\tau)=
	\calC_{\gamma,S}$, where $\gamma$ is as follows.
\begin{enumerate} 
	\item[{\it (i)}\ ] If $\tau = 2$ then $\gamma = \textmatrix 1201$.
	\item[{\it (ii)}\ ] If $\tau = -2$ then $\gamma = \textmatrix 1{2u}01$, where $u \in \F_q$ and $\jacobi uq=-1$.
	\item[{\it (iii)}\ ] If $\tau = \pm(a + 1/a)$ with $a \in \F_q^\x $ and $a^4 \ne 1$ then $\gamma = \textmatrix {a} 0 0 {a^{-1}}$.
	\item[{\it (iv)}\ ] If $\tau = \pm(\z+1/\z)$ with $\zeta^{q+1}=1$ and $\zeta^4\ne 1$ then $\gamma =E_{\zeta^2,\lambda}$, 
where $\l$ is any element of $\F_{q^2} \setminus \F_q$. 
Here $E_{\zeta,\l} $ is defined by (\ref{eq:Ezeta}), and $E_{\zeta^2,\l}=E_{\zeta,\l}^2 \in \PSL_2(\F_q)$.
\end{enumerate}
\end{theorem}
\begin{proof} All elements of $\F_q^\x$ are regular by Lemma~\ref{lem:PSL2short}, so $\inv_{Q_S}(\tau)=\calC_{\g,S}$ is defined.

	{\it (i)} and {\it (ii)}\ Let $v$ be a solution to $v^q=v+b$ where $b \in \F_q^\x$.  Then
$$Q_S(v) = \frac{(v^{q^2}-v)^{(q+1)/2}}{(v^q-v)^{(q^2+1)/2} } 
= \frac{(v + 2b - v)^{(q+1)/2}}{(v + b -v)^{(q^2+1)/2}}  = 2^{(q+1)/2} b^{(q-q^2)/2}.$$
Now $b^{(q-q^2)/2} = (b^q)^{(1-q)/2} = b^{(1-q)/2}=b^{q-1}b^{(1-q)/2}  = b^{(q-1)/2}$, so
$$Q_S(v)=  2^{(q+1)/2}b^{(q-1)/2} = 2 (2b)^{(q-1)/2} = 2 \jacobi{2b} q.$$
	If $b=2$ then $Q_S(v)=2 $ and if $b=2u$ then $Q_S(v) = -2$. Since $v^q=v+b=\textmatrix 1b01(v)$, we conclude that $\inv_{Q_S}(2) = \calC_{\textmatrix 1201}$ and
$\inv_{Q_S}(-2) = \calC_{\textmatrix 1{2u}01}$.

	{\it (iii)}\   By Lemma~\ref{lem:PSL2short}{\it (v)}, if $\inv_{Q_S}(\tau)=\calC_{\g,S}$ then  $\iota(\g)=\tau^2=(a+1/a)^2=\iota(\textmatrix a00{a^{-1}})$ and $\circ(\g)\ge3$.  By Corollary~\ref{cor:PGL2conjugate},
	$\g$ and $\textmatrix a00{a^{-1}}$ are conjugate in $\PGL_2(\F_q)$, 
	say $\g=\a \textmatrix a00{a^{-1}} \a^{-1}$. Let $\a' = \a \textmatrix d001$, where $d=1/\det(\a)$. Then $\a' \in \PSL_2(\F_q)$ and 
	$\g=\a' \textmatrix a00{a^{-1}} (\a')^{-1}$. Thus, $\inv_{Q_S}(\tau)=\calC_{\g,S}=\calC_{\textmatrix a00{a^{-1}},S}$. 

	{\it (iv)}\   By Lemma~\ref{lem:PSL2short}{\it (v)}, if $\inv_{Q_S}(\tau)=\calC_{\g,S}$ then  $\iota(\g)=\tau^2=(\z+1/\z)^2=(\z^2+1)^2/\z^2=\iota(E_{\z^2,\l})$ and $\circ(\g)\ge3$. By Corollary~\ref{cor:PGL2conjugate},
	$\g$ is conjugate to $E_{\z^2,\l}$ in $\PGL_2(\F_q)$,
	say $\g=\a E_{\z^2,\l} \a^{-1}$ where $\a \in \PGL_2(\F_q)$.  
	Here, $E_{\z^2,\l}=E_{\z,\l}^2\in \PSL_2(\F_q)$.

	If $\a\in \PSL_2(\F_q)$ then $\inv_{Q_S}(\tau)=\calC_{\g,S}=\calC_{E_{\z^2,\l},S}$ as required. 

	If $\a \not\in\PSL_2(\F_q)$ then $\det(\a)$ is
	a nonsquare.  Let $\d$ be a primitive element of $\F_{q^2}$ and let $D_{\d,\l}=C_\l \textmatrix {\d^q}00\delta C_\l^{-1}$. Then $\det(D_{\d,\l})=\d^{q+1}$ is a primitive
	element of $\F_q$, and in particular a nonsquare in $\F_q$.
	By Proposition~\ref{prop:Dickson}, its entries are rational, so it belongs to $\PGL_2(\F_q)$. It commutes with $E_{\z^2,\l}$, therefore 
	$\g=\a' E_{\z^2,\l} (\a')^{-1}$, where $\a'=\a D_{\d,\l}$. Since $\jacobi{\det(\a')}q=\jacobi{\det(\a)} q \jacobi {\det(D_{\d,\l})} q = (-1)\cdot(-1) = 1$, $\a'\in \PSL_2(\F_q)$,
	so that $\g$ and $E_{\z^2,\l}$ are in the same conjugacy class of $\PSL_2(\F_q)$.
	Then $\inv_{Q_S}(\tau) = \calC_{\g,S} = \calC_{E_{\z^2,\l},S}$ as required.
\end{proof}

\section{Acknowledgements}

The author thanks Xander Faber for reviewing this article and providing some very insightful comments. First, he observed that $\inv(\tau,q)$ is
essentially the Artin map, as explained in the introduction. This led to a change in emphasis, and even a change in the title.
Second, he greatly simplified Section~\ref{sec:cyclic} by finding
a more direct way to compute the quotient map for a cyclic group of order $\ell$ in the case where $\ell|q+1$.  
In addition, he gave many other suggestions that greatly improved the exposition.
His collegiality is immensely appreciated.


\begin{thebibliography}{10}
\bibitem{Artin} Emil Artin, \emph{Galois Theory}, Notre Dame Mathematical Lectures Number~2, University of Notre Dame Press, 1942.

\bibitem{Structure} Antonia W.\ Bluher, \emph{A structure theorem
	for finite fields}, {\sl Finite Fields and Their Applications}~{\bf 68} (2020), doi:10.1016/j.ffa.2020.101732.

\bibitem{Wilson-like} Antonia W. Bluher, \emph{New Wilson-like theorems arising from Dickson polynomials}, {\sl Finite Fields Appl.}~{\bf 72} (2021), 
	doi:10.1016/j.ffa.2021.101819.

\bibitem{Permutation} Antonia W.\ Bluher, \emph{Permutation properties of Dickson and Chebyshev polynomials and connections to number theory},
Finite Fields Appl.~76 (2021), 
doi:10.1016/j.ffa.2021.101899. 

\bibitem{Brewer} B. W. Brewer, \emph{On certain character sums}, Trans.\ Amer.\ Math.\ Soc.\ {\bf 99}, pp.~241--245, (1961), 
doi:10.1090/S0002-9947-1961-0120202-1.

\bibitem{CMPZ18} Bence Csajb\'{o}k, Giuseppe Marino, Olga Polverino, and Yue Zhou, 
\emph{MRD codes with maximum idealizers}, Discrete Mathematics~{\bf 343}(9), 
(2020), doi:10.1016/j.disc.2020.111985.

\bibitem{CMPZ19} Bence Csajb\'{o}k, Giuseppe Marino, Olga Polverino, and Ferdinando Zullo, \emph{A characterization of linearized polynomials with
maximum kernel}, \emph{Finite Fields Appl.}~{\bf 56}, 109--130, (2019),
doi:10.1016/j.ffa.2018.11.009

\bibitem{Dickson} Leonard Eugene Dickson, {\it Linear Groups with an
Exposition of the Galois Theory}, Dover Publications, 1958.

\bibitem{DD} John F. Dillon and Hans Dobbertin, New cyclic difference sets with Singer parameters,
Finite Fields and Their Applications 10 (2004), 342--389,
doi:10.1016/j.ffa.2003.09.003.

\bibitem{DF} David S.\ Dummit and Richard M.\ Foote, \emph{Abstract Algebra}, Third Ed., John Wiley \& Sons, 2004.

\bibitem{Goss} David Goss, \emph{Basic Structures of Function Field Arithmetic}, Springer, 1996.

\bibitem{Lang} Serge Lang, \emph{Algebra (Revised Third Edition)}, Graduate Textbooks in Mathematics 211, Springer-Verlag, 2002.

\bibitem{MM} Gary McGuire and Daniela Mueller, \emph{Some results on linearized trinomials that split completely}, arXiv:1905.11755v2 [math.NT], 8 June 2019.

\bibitem{MS} Gary McGuire and John Sheekey, \emph{A characterization of the number of roots of linearized and projective polynomials in the field of
coefficients}, \emph{Finite Fields and Appl.} 57, 68--91 (2019),
doi:10.1016/j.ffa.2019.02.003.

\bibitem{Mumford} David Mumford, \emph{Abelian Varieties}, with appendices by C.\ P.\ Ramanujam and Yuri Manin, published for the Tata Institute
of Fundamental Research by Hindustan Book Agency, with international distribution by the American Mathematical Society, 2012.

\bibitem{Neukirch} J.~Neukirch, \emph{Class Field Theory}, Springer, 1986.

\bibitem{Rosen} Michael Rosen, \emph{Number Theory in Function Fields}, Graduate Texts in Mathematics, Springer, 2002.

\bibitem{VdW} B.\ L.\ van der Waerden, \emph{Algebra, Volume~I}, Springer-Verlag, 1991.

\end{thebibliography}
\end{document}